\documentclass[11pt]{article} 
\usepackage{a4wide}
\usepackage{amsmath,amsfonts,amsthm}

\newtheorem{theorem}{Theorem}
\newtheorem{remark}{Remark}

\newtheorem{corollary}{Corollary}
\newtheorem{proposition}{Proposition}
\newtheorem{lemma}{Lemma} 
\newtheorem{claim}{Claim} 

\begin{document}

\title{Asymptotic stability of solitons for 
 the  Benjamin-Ono equation}
\author{C.E. Kenig$^{(1)}$ and Y. Martel$^{(2)}$}
\date{
(1) Department of Mathematics, University of Chicago,\\
5734 University ave., Chicago, Illinois 60637-1514,\\
cek@math.uchicago.edu\\\quad \\
(2) Universit\'e de Versailles Saint-Quentin-en-Yvelines,
 Math\'ematiques, \\
 45, av. des Etats-Unis,
 78035 Versailles cedex, France\\   
 martel@math.uvsq.fr}
\maketitle
\begin{abstract}
In this paper, we prove the asymptotic stability of the family of solitons of the 
Benjamin-Ono equation in the energy space. The proof is based on a Liouville 
property for solutions close to the solitons for this equation, in the spirit of
\cite{MM1}, \cite{MMas1}. As a corollary of the proofs, we obtain the asymptotic stability of exact multi-solitons.
\end{abstract}
 
\section{Introduction}
We consider the  Benjamin-Ono equation (BO)
\begin{equation}\label{BO}
    u_t + \mathcal{H} u_{xx} + u  u_x = 0, \quad
    (t,x)\in \mathbb{R}\times \mathbb{R},
\end{equation}
where   $\mathcal{H}$ denotes the Hilbert transform
\begin{equation}\label{hilbert}
\mathcal{H} u(x)=\frac 1\pi \, \mathrm{p.v.} \int_{-\infty}^{+\infty} \frac {u(y)}{y-x} dy
=\frac 1\pi\lim_{\varepsilon\to 0} \int_{|y-x|> \varepsilon} \frac {u(y)}{y-x} dy.
\end{equation}
Note that with this notation, 
$\int u_x \mathcal{H} u =\int |D^{\frac 12} u|^2 = \|u\|_{\dot H^{\frac 12}}^2$.

The Cauchy problem for \eqref{BO} is globally well-posed in $H^s$, for any $s\geq 0$
(see Tao \cite{Tao} for $s\geq 1$ and
Ionescu and Kenig \cite{IK} for the case $s\geq 0$, see also Burq and Planchon \cite{BP} for the case $s>\frac 14$).
Moreover, for solutions in the energy space $H^{\frac 12}$
the following quantities are invariant
\begin{equation}\label{invariants}
  \int u^2(t,x)dx =\int u^2(0,x)dx, \quad 
E(t)=\int \Big( u_x \mathcal{H} u - \tfrac 13 {u^{3}} \Big) (t,x) dx=
E(0).
\end{equation}
Recall the scaling and translation invariances of equation \eqref{BO}
\begin{equation}\label{inv}
\text{if $u(t,x)$ is solution then $\forall c>0$, $x_0\in \mathbb{R}$,
$v(t,x)=c \, u(c^2 t,c(x-x_0))$ is solution.}
\end{equation}

We call soliton any travelling wave solution $u(t,x)=Q_c(x-x_0-ct)$,
where $c>0$, $x_0\in \mathbb{R}$, and $Q_c(x)=c  Q(c x)$   solves:
\begin{equation}\label{stationary}
- \mathcal{H} Q' + Q - \tfrac 12 {Q^{2}} =0, \quad
Q\in H^{\frac 12},\quad Q>0.
\end{equation}
It is known that there is a unique (up to translations)   solution of \eqref{stationary}, which is
\begin{equation}\label{explicitQ}
Q(x)=\frac 4{1+x^2}.
\end{equation}
(see Benjamin \cite{B} and Amick and Toland \cite{AT} for the uniqueness statement).
This solution is stable (see Bennet et al. \cite{BB} and Weinstein \cite{W})
in the following sense.

\medskip

\noindent\textbf{Stability of soliton in the energy space (\cite{BB}, \cite{W}).}
\emph{
There exist  $C, \alpha_0>0$ such that
if $u_0\in H^{\frac 12}$ satisfies
$\|u_0-Q\|_{H^{\frac 12}}=\alpha\leq \alpha_0$ then 
the solution  $u(t)$  of \eqref{BO} with $u(0)=u_0$ satisfies}
$$\sup_{t\in \mathbb{R}} \inf_{y\in \mathbb{R}}
\|u(t)-Q(.-y)\|_{H^{\frac 12}}\leq C\alpha .$$
See a sketch of proof of this result  in Section \ref{sec:51}.

\medskip

The main result of this paper is  the asymptotic stability of the family of solitons
of \eqref{BO}. Then, we consider the multisoliton case (see Section \ref{sec:5}).

\begin{theorem}[Asymptotic stability of solitons in the energy space]\label{TH1}\quad \\
There exist $C, \alpha_0>0$ such that if $u_0\in H^{\frac 12}$ satisfies
$\|u_0-Q\|_{H^{\frac 12}}=\alpha\leq \alpha_0$, then there
exists $c^+>0$ with $|c^+-1|\leq C \alpha$ 
and a $C^1$ function $\rho(t)$ such that the solution $u(t)$ of
\eqref{BO} with $u(0)=u_0$ satisfies  
\begin{align}
& u(t,.+\rho(t)) \rightharpoonup Q_{c^+} \quad \text{in $H^{\frac 12}$ weak,}
\qquad \|u(t)-Q_{c^+}(.-\rho(t))\|_{L^2(x> \frac t{10})} \to 0,\\
& \rho'(t)\to c^+ \quad \text{as $t\to +\infty$.}
\end{align}
\end{theorem}

The proof of Theorem \ref{TH1} is based on the following rigidity result.

\begin{theorem}[Nonlinear Liouville property]\label{TH2}\quad \\
There exist $C,\alpha_0>0$ such that if $u_0\in H^{\frac 12}$ satisfies
$\|u_0-Q\|_{H^{\frac 12}}=\alpha\leq \alpha_0$ and if 
the solution $u(t)$ of
\eqref{BO} with $u(0)=u_0$ satisfies for some function $\rho(t)$  
\begin{equation}\label{th2:1}
\forall \varepsilon>0,\exists A_{\varepsilon}>0, \text{ s.t. }
\forall t\in \mathbb{R},\quad
\int_{|x|>A_\varepsilon} u^2(t,x+\rho(t)) dx <\varepsilon,
\end{equation}
then there exist $c_1>0$, $x_1\in \mathbb{R}$, such that
\begin{equation}
u(t,x)=Q_{c_1}(x-x_1-c_1t),
\quad 
|c_1-1|+|x_1|\leq C \alpha.
\end{equation}
\end{theorem}

\begin{remark}\label{re:1}
In Theorem \ref{TH1}, the convergence of $u(t)$ to
$Q_{c^+}$ as $t\to +\infty$ is obtained strongly in $L^2$ in the region $x>\frac t{10}$.
The value $\frac 1{10}$ is somewhat arbitrary, the result  holds for $x>\varepsilon t$,
for any $\varepsilon>0$, provided $\alpha_0=\alpha_0(\varepsilon)>0$ is small enough.
Note that this result is optimal in $L^2$ since $u(t)$ could contain other small 
(and then slow) solitons and since in general $u(t)$ does not go to $0$ in $L^2$ for $x<0$.
For example, 
if $\|u(t)-Q_{c^+}(.-\rho(t))\|_{H^{\frac 12}(\mathbb{R})}\to 0$ as $t\to +\infty$,
then $E(u)=E(Q_{c^+})$ and $\int u^2=\int Q_{c^+}^2$ and so 
by the variational characterization  of $Q(x)$ (see \cite{W}), $u(t)=Q_{c^+}(x-x_0-c^+t)$
is exactly a soliton.

Under the assumptions of Theorem \ref{TH1}, we expect strong convergence in $H^{\frac 12}$
to be true as well in the same local sense ($x>\varepsilon t$). This could 
require some more analysis.

By the methods of this paper, we are able to obtain the following weaker result (Section~\ref{sec:4.3})
\begin{equation}\label{inre:1}
\lim_{t\to +\infty} \int_t^{t+1} \|u(s,.+\rho(s))-Q_{c^+}\|_{H^{\frac 12}_{loc}}^2 ds =0.
\end{equation}
\end{remark}

The proof of Theorem \ref{TH1} follows the approach of \cite{MMjmpa}, \cite{MM1}, concerning the case of the generalized KdV equations, where the asymptotic stability of the family of solitons is deduced from
a Liouville type theorem such as Theorem~\ref{TH2}.
Moreover, similarly as in \cite{MM1}, the proof of  Theorem~\ref{TH2}
follows from a
Liouville property on the linearized equation around $Q$,
see Theorem~\ref{LINEARLIOUVILLE} in Section \ref{sec:3}.

With respect to the gKdV case, there are two main difficulties : (1)  
$L^2$ monotonicity type results, which are
similar to the ones for the gKdV equations (\cite{MM1}), but whose proof are more subtle due to
the nonlocal nature of the (BO) operator (see Section \ref{sec:2}). For this part, we use
a Kato type identity for \eqref{BO} (see \cite{GV} and \cite{Po}).
 
 (2) The proof of the linear Liouville theorem, which requires the analysis of
some linear operators related to $Q$. Note that for this part, we use the fact that
$Q(x)$ is explicit, and some known results about the linearized equation around $Q$
(\cite{BB}, \cite{W}). We point out that except for this part of the analysis, all the
arguments are quite flexible and could be applied to generalized versions of
the (BO)) equation. In particular, we do not use the integrability property of the
equation.

As a corollary of the proof of Theorem \ref{TH1} and of Theorem \ref{TH2}, we obtain 
stability and asymptotic stability of  multisoliton solutions.
See  Theorem 4 in Section \ref{sec:5}  for a precise statement. After the paper was finished and submitted, we learned that S. Gustafson, H.~Takaoka, and T-P. Tsai \cite{GTT} have obtained independently the stability part of Theorem 4. Note that the main result of the present paper, i.e. asymptotic stability of (single or multi-) solitons is not addressed in \cite{GTT}.

The rest of the paper is organized as follows. In Section \ref{sec:2}, we prove 
$L^2$ monotonicity type results in the context of Theorem \ref{TH1}.
In Section \ref{sec:3}, we state and prove the linear Liouville Theorem, which is
the main ingredient of the proof of Theorem \ref{TH2}. In Section \ref{sec:4},
we prove Theorems \ref{TH1} and \ref{TH2} using Sections \ref{sec:2} and \ref{sec:3}.
Section \ref{sec:5} is devoted to the multisoliton case.
In Section \ref{sec:6}, we prove some weak convergence and well-posedness results
used in the proofs. Finally, Appendix \ref{sec:A} contains the proof of some technical points. 

\medskip

\noindent\textbf{Acknowledgments.}
The first author is partly supported by the NSF grant DMS-0456583.
This work was initiated when the second author was visiting the University of Chicago.
He would like to thank the Department of Mathematics for its hospitality.
The second author is partly supported by the Agence Nationale de la Recherche (ANR {ONDENONLIN}).

\section{Monotonicity arguments for solutions close to $Q$}\label{sec:2}

\subsection{Modulation}
\begin{lemma}[Choice of translation parameter]\label{MODULATION}
There exist $C,\alpha_0>0$ such that for any $0<\alpha<\alpha_0$,
if $u(t)$ is an $H^\frac 12$ solution of \eqref{BO} such that
\begin{equation}
\forall t\in \mathbb{R},\quad \inf_{r\in \mathbb{R}}\|u(t)-Q(.-r)\|_{H^{\frac 12}}
< \alpha,
\end{equation}
then there exists $\rho(t)\in C^1(\mathbb{R})$ such that
$$
\eta(t,x)=u(t,x+\rho(t))-Q(x) 
$$
satisfies
\begin{equation}\label{ortho}\begin{split}
\forall  t\in \mathbb{R},\quad 
    & \int Q'(x) \eta(t,x)dx=0,\quad
	  \|\eta(t)\|_{H^{\frac 12}}\leq C \alpha,\\&
  	 |\rho'(t)-1|\leq  C \left(\int \frac {\eta^2(t,x)}{1+x^2}dx\right)^{\frac 12}\leq  C\|\eta(t)\|_{L^2}.
\end{split}\end{equation}
\end{lemma}
\noindent\emph{Proof of Lemma \ref{MODULATION}.}
This follows from standard arguments 
(see e.g. \cite{BSS}, Lemma 4.1, \cite{MMgafa}, Proposition 1 and Lemma 4).

\noindent\emph{Time independent arguments.}
For $u\in H^{\frac 12}$ and $y\in \mathbb{R}$, set
$$
I_{y}(u)=\int Q'(x) (u(x+y)-Q(x)) dx \quad \text{so that}\quad
{\frac {\partial I_{y}}{\partial y}}~_{| y=0, u=Q} = \int (Q')^2>0.
$$
Thus, by the implicit function theorem, there exists $\alpha_1>0$, $V$
a neighborhood of $0$ in $\mathbb{R}$  and a unique $C^1$ map:
$$
y:\{u\in H^{\frac 12}, ~ \|u-Q\|_{H^{\frac 12}} \leq \alpha_1\}\to V~
\text{such that $I_{y(u)}(u)=0$,
$|y(u)|\leq C\|u-Q\|_{H^{\frac 12}}$ .} 
$$
We uniquely extend the $C^1$ map $y(u)$ to 
$U_{\alpha_1}=\{u\in H^{\frac 12}, ~ \inf_r\|u(.+r)-Q\|_{H^{\frac 12}} \leq \alpha_1  \}$ so that for all $u$ and $r$,
$y(u)=y(u(.+r))+r$. Then, we set $\eta_u(x)=u(x+y(u))-Q(x)$, so that
$$\int \eta_u Q'=0\quad \text{and}\quad \|\eta_u\|_{H^\frac 12}\leq C \|u-Q\|_{H^{\frac 12}}.$$

\noindent\emph{Estimates depending on $t$.} 
For all $t$, we define $\rho(t)=y(u(t))$ and $\eta(t)=\eta_{u(t)}$.
To conclude the proof of the lemma, we just have to prove the estimate on
$\rho'(t)-1$.

We perform formal computations which can be justified for $H^{\frac 12}$ solutions
by density and continuous dependence arguments.
The function  $\eta(t,x)$ satisfies the following equation:
\begin{equation}\label{eqofeps}
\eta_t = (\mathcal{L} \eta - \tfrac 12 \eta^2)_x
+ (\rho'-1) (Q+\eta)_x\quad \text{where $\mathcal{L} \eta = - \mathcal{H} \eta_x
+ \eta - Q \eta$.}
\end{equation}
Thus, multiplying the equation of $\eta$ by $Q'$ and using $\int \eta Q'=0$, we obtain
\begin{equation}\label{nouveau}
(\rho'-1) \left[ \int (Q')^2 - \int \eta Q''\right]
= \int \eta \mathcal{L}(Q'') -\tfrac 12 \int \eta^2 Q'',
\end{equation}
which finishes the proof for $\alpha_0$ small enough.

\begin{remark}\label{rk:10}
By the proof of Lemma \ref{MODULATION}, $\rho(t)$ depends continuously on $u(t)$ in $H^{\frac 12}$.
In particular, let $u(t)$ satisfy the assumptions of  Lemma \ref{MODULATION} with $u(0)=u_0$. If $u_n(0)\to u_0$ in $H^{\frac 12}$ as $n\to +\infty$,  then by continuous dependence
(see \cite{IK}), we obtain for all $t\in \mathbb{R}$, $\rho_n(t)\to \rho(t)$ as $n\to +\infty$,
where $\rho_n(t)$ is defined from $u_n(t)$ ($u_n(t)$ is the solution of \eqref{BO} corresponding to $u_n(0)=u_{0n}$).

Note also that in the proof of Lemma \ref{MODULATION}, we can replace the space $H^{\frac 12}$ 
by $L^2$, so that in the same context if $u_n(0)\to u_0$ in $L^2$ as $n\to +\infty$ then
for all $t\in \mathbb{R}$, $\rho_n(t)\to \rho(t)$ as $n\to +\infty$ (see continuous dependence in $L^2$ also in \cite{IK}).

Finally, for future reference, we justify that if $u_n\rightharpoonup u$ in $H^{\frac 12}$ weak, then
$y(u_n)\to y(u)$, where $y(u)$ is defined in the proof of Lemma \ref{MODULATION}.
Indeed, in this proof, by the decay of $Q'(x)$, we can also replace $H^{\frac 12}$
by the weighted space $L^2(\frac 1{1+|x|} dx)$, so that 
if $u_n\to u$ in $L^2_{loc}$ and $\|u_n\|_{L^2}+\|u\|_{L^2}\leq C$, then $y(u_n)\to y(u)$
as $n\to +\infty$.
\end{remark}

In the rest of  this section, we present monotonicity arguments on $L^2$ quantities 
for both $u(t)$ and   $\eta(t)$, in the context of Lemma \ref{MODULATION}.
These results
  are reminiscent of similar  results for the gKdV equation in \cite{MM1}
  and \cite{MMas2}, but due to
the nonlocal nature of the operator $\mathcal{H}$, the proofs are more involved.

\subsection{Monotonicity results for $u(t)$}\label{SECmonou}

Let $A>1$ to be chosen later and set
\begin{equation}\label{defphi0}
\varphi(x)=\varphi_A(x)=\frac \pi 2 + \arctan\Big(\frac x A\Big)\quad \hbox{so that}\quad
\varphi'(x)=\frac {\frac 1A} {1+(\frac x A)^2}>0.
\end{equation}
\begin{proposition}\label{MONOTONICITY1}
Let $0<\lambda<1$.
Under the assumptions of Lemma \ref{MODULATION}, 
for $\alpha_0$ small enough and $A$ large enough, there exists $C>0$ such that for all $x_0>1$,  $t_1\leq t_2$,
\begin{enumerate}
\item Monotonicity on the right of the soliton: 
\begin{equation}\label{monotonicity1}
\int u^2(t_2,x)\varphi(x-\rho(t_2)-x_0) dx 
   \leq \int u^2(t_1,x)\varphi(x-\rho(t_1)-\lambda (t_2-t_1)-x_0) dx +\frac C {x_0}.
\end{equation}
\item Monotonicity on the left of the soliton:  \begin{equation}\label{monotonicity2}
\int u^2(t_2,x)\varphi(x-\rho(t_2)+\lambda (t_2-t_1)+x_0) dx 
   \leq \int u^2(t_1,x)\varphi(x-\rho(t_1)+x_0) dx +\frac C {x_0}.
\end{equation}
\end{enumerate}
\end{proposition}

\noindent\emph{Proof of Proposition \ref{MONOTONICITY1}.}\quad 
First, we note that \eqref{monotonicity2} is a consequence of \eqref{monotonicity1}
and the $L^2$ norm conservation. Indeed, let
$v(t,x)=u(-t,-x)$. Then $v(t)$ is a solution of \eqref{BO} satisfying the assumptions
of Lemma \ref{MODULATION} and $\rho_v(t)=-\rho(-t)$.
Thus, from \eqref{monotonicity1} applied on $v(t,x)$, we deduce
$$
 \int u^2(-t_2,x)\varphi(-x+\rho (-t_2)-x_0) dx 
   \leq  \int u^2(-t_1,x)\varphi(-x+\rho(-t_1)-\lambda (t_2-t_1)-x_0) dx +\frac C {x_0}.
$$
Since $\varphi(x)=\pi -\varphi(-x)$, from $\int u^2(-t_2)=\int u^2(-t_1)$, we obtain
$$
 \int u^2(-t_2,x)\varphi( x-\rho (-t_2)+x_0) dx +\frac C {x_0}
   \geq \int u^2(-t_1,x)\varphi(x-\rho (-t_1)+\lambda (t_2-t_1)+x_0) dx ,
$$
which is exactly formula \eqref{monotonicity2} for $t_2'=-t_1$, $t_1'=-t_2$.

\medskip

We are reduced to  prove \eqref{monotonicity1}.
We perform calculations on regular solutions
and then use density arguments and continuous dependence to obtain the result
in the framework of Lemma \ref{MODULATION}.

First, we recall a Kato type identity for solutions of the BO equation. 
By direct computations, we have
\begin{equation}\label{kato}\begin{split}
\frac 12 \,\frac d{dt}\int u^2(t,x) \varphi(x) dx  & =
\int u_t u \varphi(x) dx  = - \int (\mathcal{H} u_{xx} + u  u_x) u \varphi(x) dx\\
& = \int (\mathcal{H} u_x) ( u\varphi'(x)+u_x \varphi(x))dx + \frac 1{3} \int u^{3} \varphi'(x)dx.
\end{split}
\end{equation}

For the first term in \eqref{kato}, we prove the following result.

\begin{lemma}\label{FIRSTTERM}
 For all $u\in H^1(\mathbb{R})$,
\begin{equation}\label{firstterm}
\int (\mathcal{H} u_x)  u\varphi'(x) dx
\leq \frac C A \int u^2 \varphi'(x) dx.
\end{equation}
\end{lemma}

\noindent\emph{Proof of Lemma \ref{FIRSTTERM}}.
For $f\in L^2(\mathbb{R})$, we define the harmonic extension of $f$ on
$\mathbb{R}\times \mathbb{R}_+=\mathbb{R}^2_+$,
\begin{equation}
\forall x\in \mathbb{R},\quad
F(x,0)=f(x) \quad \mathrm{and}\quad 
F(x,y)=\frac 1\pi \int_{-\infty}^{\infty} \frac {y} {(x-x')^2 +  y^2} \, f(x') \, dx', \quad
\hbox{if $y>0$.}
\end{equation}
In particular, recall that $\mathcal{H} f'(x)=\partial_y F(x,0)$
(see Stein \cite{ST} Chapter III, and the Introduction of Toland \cite{To}).

We denote by $\Phi(x,y)$ the harmonic extension of $\varphi'(x)$ and $U(x,y)$
the harmonic extension of $u(x)$ on $\mathbb{R}\times \mathbb{R}_+$.
Note that $\Phi(x,y)$ is explicitly given by
\begin{equation}\label{defphicap}
\Phi(x,y)=\frac 1 A \frac {1+\frac y A}{(\frac x A)^2+ ( 1+\frac y A)^2}.
\end{equation}
Then, by the Green Formula on $\mathbb{R}^2_+$ 
(using decay properties of $\Phi(x,y)$ and $\Delta U^2=2|\nabla U|^2$), we obtain
formally
\begin{equation}\label{toprove}\begin{split}
\int (\mathcal{H}u_x) u \varphi' & = \int \partial_y U(t,x,0) U(t,x,0) \Phi(x,0) dx
=\frac 12 \int_{y=0} \partial_y (U^2) \Phi dx \\
& = -\frac 12 \iint_{\mathbb{R}^2_+} (\Delta U^2) \Phi
+ \frac 12 \iint_{\mathbb{R}^2_+}  U^2 \Delta \Phi + \frac 12 \int_{y=0} U^2 \partial_y \Phi \\
& = - \iint_{\mathbb{R}^2_+} |\nabla U|^2 \Phi + \frac 12 \int u^2 (\mathcal H \varphi'') dx.
\end{split}\end{equation}
See Appendix \ref{secTOPROVE}  for a rigorous proof of \eqref{toprove}.
Since  $\Phi\geq 0$ on $\mathbb{R}^2_+$, we obtain
\begin{equation} 
\int (\mathcal{H}u_x) u \varphi'  \leq \frac 12  \int u^2 (\mathcal H \varphi'').
\end{equation}
By explicit computations, since $
\mathcal H\big(\frac 1{1+x^2}\big) = -\frac x{1+x^2}$, we have
\begin{equation}\label{Hdephi}
\mathcal H \varphi' = - \frac 1  {A^2} \frac {x}{1+(\frac xA)^2},\quad
\mathcal H \varphi'' = \frac 1 A \varphi' - 2 (\varphi')^2
\quad\hbox{and}\quad 
\mathcal H \varphi'' \leq \frac 1A \varphi'.
\end{equation}
Lemma \ref{FIRSTTERM} follows.

\medskip

For the second term in \eqref{kato}, we have the following.

\begin{lemma}\label{SECONDTERM}
 For all $u\in H^1(\mathbb{R})$,
\begin{equation}\label{secondterm}
\left|\int (\mathcal{H} u_x) u_x \varphi dx\right|
\leq \frac C A \int u^2 \varphi'(x) dx.
\end{equation}
\end{lemma}

\noindent\emph{Proof of Lemma \ref{SECONDTERM}}.
We prove \eqref{secondterm} for $u$ smooth and compactly supported in $\mathbb{R}$,
the general case will follow by  a density argument.

Since the limit in \eqref{hilbert} holds in $L^2$ 
(see Stein \cite{ST}, Chapter II), we have
\begin{equation}\label{quinze}\begin{split}
&\int (\mathcal{H} u_x) u_x \varphi dx
 = \frac 1\pi \int \mathrm{p.v.}\bigg( \int \frac {u_x(y)}{y-x} dy \bigg) u_x(x) \varphi(x) dx
\\ &=\frac 1\pi\lim_{\varepsilon\to 0} \iint_{|y-x|>\varepsilon} u_x(y) u_x(x) \frac {\varphi(x)}{y-x} dydx  \\
&= \frac 1{2\pi} \iint   u_x(y)u_x(x)  \, \frac {\varphi(x)-\varphi(y)}{y-x} dxdy
=\frac 1{2\pi} \iint  u(y)u(x)  K_\varphi(x,y) dxdy,
\end{split}\end{equation}
by symmetry and then integration by parts, where 
\begin{equation}\label{defK}
 K_{\varphi}(x,y)=-  \frac {\partial^2}{\partial x\partial y}  \bigg(\frac {\varphi(x)-\varphi(y)}{x-y}\bigg)
 =\frac   {
 2 (\varphi(x)-\varphi(y))-(\varphi'(x)+\varphi'(y))(x-y)} {(x-y)^3}.
\end{equation}
Note that all the integrals in \eqref{quinze} make sense since $u(x)$ is compactly
supported, 
$(\varphi(x)-\varphi(y))/(x-y)$ is bounded and moreover,
by subtracting the following two Taylor formulas:
\begin{equation*}\begin{split}
& \varphi(x)=\varphi(y)+(x-y)\varphi'(y) + \frac 12 (x-y)^2 \varphi''(y)
+\frac 16 (x-y)^3 \varphi'''(x_1),\\
& \varphi(y)=\varphi(x)+(y-x)\varphi'(x) + \frac 12 (y-x)^2 \varphi''(x)
+\frac 16 (y-x)^3 \varphi'''(x_2),
\end{split}\end{equation*}
where $x_1,x_2\in (y,x)$, we find:
\begin{equation}\label{pf}
K_{\varphi}(x,y)
=\frac 12\frac {\varphi''(y)-\varphi''(x)}{x-y} 
+\frac 16 (\varphi'''(x_1)+\varphi'''(x_2)),
 \end{equation}
which is also bounded on $\mathbb{R}^2$.
Note also that by explicit computations, we have
\begin{equation}\label{bg}
\varphi'''(x)=
\frac {\varphi'(x)}{A^2} \left(\frac {-2}{1+\big(\frac x A\big)^2}
+ \frac {8 \big(\frac xA\big)^2} {\big(1+\big(\frac xA\big)^2\big)^2}\right)
=
\frac {\varphi'(x)}{A} \left(-2\varphi'(x) + \frac 8A {x^2} (\varphi')^2\right).
\end{equation}

We are reduced to prove the following estimate
\begin{equation}\label{cc}
\left|\iint  u(y)u(x)  K_\varphi(x,y) dxdy\right|\leq \frac C A\int u^2 \varphi'(x)dx.
\end{equation}
We consider only the case $|y|<|x|$ (by symmetry), and we divide
$\{(x,y),\,:\, |y|<|x|\}$ into the following regions:

\medskip

$\bullet$ $\Sigma_1=\{(x,y) \,:\, x>A,\, 0<y<\frac x 2\}.$
For $(x,y)\in \Sigma_1$,   by \eqref{defK} and the fact that
$\varphi'$ is decreasing on $\mathbb{R}^+$, we have
$$
|K_{\varphi}(x,y)|\leq \frac 4{ (x-y)^2} \sup_{[y,x]}\varphi'
\leq\frac {16}{x^2}  \varphi'(y)
=\frac {16} {A^2} \frac 1 {(\frac xA)^2} \varphi'(y)
\leq \frac {32} A \varphi'(x) \varphi'(y). 
$$
Thus, by Cauchy-Schwarz inequality, since $\int \varphi'(x)=\pi$, we obtain
\begin{equation*}\begin{split}
\left|\iint_{\Sigma_1}  u(y)u(x)  K_\varphi(x,y) dxdy\right|& \leq
\frac {C}{A} \int  |u(x)|Ê\varphi'(x)dx \int |u(y)| \varphi'(y) dy\\
& \leq \frac {C\pi }{A} \int u^2(x) \varphi'(x) dx.
\end{split}\end{equation*}

The case of the region
$\Sigma_1^-=\{(x,y) \,:\, x<-A,\, \frac x 2 <y<0\}$
is similar.

\medskip

$\bullet$ $\Sigma_2=\{(x,y) \,:\, x>A, \, -x<y<0\}.$
For $(x,y)\in \Sigma_2$, we have by \eqref{defK},
$|x-y|=x-y>x>\frac 12(x+A)$, $\varphi'(y)>\varphi'(x)$ and so
by \eqref{defK} and $\varphi$ bounded, we obtain
$$
|K_{\varphi}(x,y)|\leq \frac C {(x+A)^3}+\frac {C\varphi'(y)} {x^2}.
$$
For the term $\frac {C\varphi'(y)} {x^2}$, we argue as for $\Sigma_1$.
For the other term, 
by Cauchy-Schwarz' inequality  and the expression of $\varphi'$, we have
\begin{equation*}\begin{split}
& \iint_{\Sigma_2}  |u(y)| |u(x)| \frac 1{(x+A)^3} dxdy\leq 
C \left(\iint_{\Sigma_2} \frac {u^2(x) } {(x+A)^3} dxdy\right)^{\frac 12}
\left(\iint_{\Sigma_2} \frac {u^2(y)} {(x+A)^3} dxdy\right)^{\frac 12} \\
&\leq  C
  \left(\frac 12 \int_{x>A} \frac {u^2(x) } {(x+A)^2} dx\right)^{\frac 12}
\left(\frac 12 \int_{y<0} \frac {u^2(y)} {(-y+A)^2} dy\right)^{\frac 12}  
\leq \frac {C'} A \int u^2(x) \varphi'(x) dx.
\end{split}\end{equation*}

The case of $\Sigma_2^-=\{(x,y) \,:\, x<-A, \, 0<y<-x\}$ is similar
to $\Sigma_2$.

\medskip

$\bullet$ $\Sigma_3=\{(x,y) \,:\, |x|<A,\, |y|< |x|\}.$
For $(x,y)\in \Sigma_3$, and $ |s|<|x|$, we have $\frac 1{2A} \leq \varphi'(s)\leq \frac 1{A}$
and thus, from \eqref{pf} and \eqref{bg}, we obtain
$|K_{\varphi}(x,y)|\leq C \sup_{|s|<|x|}|\varphi'''(s)| 
\leq \frac C {A^3}\leq \frac {C} A \varphi'(x)\varphi'(y)$.
We finish as for $\Sigma_1$.

\medskip

$\bullet$ $\Sigma_4=\{(x,y) \,:\, x>A,\, \frac 12 x<y<x\}.$
For $(x,y)\in \Sigma_4$, and $y<s<x$, we have from \eqref{bg}:
$$
|\varphi'''(s)|\leq \frac {10} A (\varphi'(s))^2\leq 
\frac {10} A \varphi'(y)\varphi'(s)\leq \frac {40}A \varphi'(y)\varphi'(x)$$
  thus
$ 
|K_{\varphi}(x,y)|\leq   \frac C A \varphi'(x)\varphi'(y),$ 
and we conclude as for $\Sigma_1$.
The case of $\Sigma_4^-=\{(x,y) \,:\, x<-A,\, x<y<\frac x 2\}$ is similar

In conclusion, we have obtained \eqref{cc} and Lemma \ref{SECONDTERM} is proved. 

\medskip

From \eqref{kato},  Lemmas \ref{FIRSTTERM} and \ref{SECONDTERM}, there exists $C_0>0$ such that
\begin{equation}\label{bmbis}
\frac 12 \,\frac d{dt}\int u^2(t,x) \varphi(x) dx 
\leq 
\frac {C_0}A  \int u^2(t,x) \varphi'(x) dx+\frac 13 \int |u^3(t,x)| \varphi'(x)dx.
\end{equation}

Now, let $u(t)$ be a solution of \eqref{BO} satisfying the assumptions of
Lemma \ref{MODULATION} on $\mathbb{R}$. 
Let $\eta(t)$,  $\rho(t)$ be associated to the decomposition of $u(t)$ on $I$ as in Lemma \ref{MODULATION}.

Let $0<\lambda<1$, $t_0\in [t_1,t_2]$ and $x_0\geq 1$. For any $t\in [t_1,t_0]$, $x\in \mathbb{R}$,
we set 
\begin{equation} 
  \widetilde x=x-x_0 - \rho(t) - \lambda (t_0-t),\quad 
   M_{\varphi}(t)= \frac 12 \int u^2(t,x) \varphi(\widetilde x) dx.
\end{equation}
Then, by \eqref{bmbis}, we find
\begin{equation}
M_{\varphi}'(t) \leq -\frac 12\left(\rho'(t)-   \lambda-\frac {2 C_0}A\right) \int u^2(t) \varphi'(\widetilde x)
+\frac 13 \int |u(t)|^3 \varphi'(\widetilde x).
\end{equation}
Fix now $A>0$ large enough so that $\frac {2C_0}A \leq \frac 14 (1-\lambda)$.
Then, by \eqref{ortho}, we   choose $\alpha_0>0$ small enough so that 
$\forall t\in I$, $\rho'(t)-\lambda>\frac 12(1-\lambda)$.
Therefore, we obtain
\begin{equation}\label{step3}
M_{\varphi}'(t) \leq -\ \frac 18 (1-\lambda) \int u^2(t) \varphi'(\widetilde x)
+\frac 13 \int |u(t)|^3 \varphi'(\widetilde x).
\end{equation}
Finally, we estimate the nonlinear term $\int |u(t)|^3 \varphi'(\widetilde x)$.
We first observe:
\begin{equation}
\int |u(t)|^3 \varphi'(\widetilde x)
\leq C \int   Q^3(x{-}\rho(t)) 
\varphi'(\widetilde x) dx + C \int  |\eta(t,x)|^3 \varphi'(\widetilde x) dx.
\end{equation}

For the first term, we distinguish two regions in $x$:

\medskip

$\bullet$ $\Omega_1=\{x\,:\, x<\rho(t)+\frac 12 x_0+\frac 12 \lambda(t_0-t)\}$. For $x\in \Omega_1$, we have
$\widetilde x<-\frac 12 x_0-\frac 12 \lambda (t_0-t)$, and thus
$$\varphi'(\widetilde x)\leq \frac C {(x_0+\lambda (t_0-t) )^2}.$$
This implies
\begin{equation}\label{partI}
\int_{\Omega_1} Q^3(x{-}\rho(t)) \varphi'(\widetilde x)
\leq \frac C {(x_0+\lambda (t_0-t) )^2}\int Q^3 
\leq \frac C {(x_0+\lambda (t_0-t) )^2}.
\end{equation}

$\bullet$ $\Omega_2=\{x>\rho(t)+\frac 12 x_0+\frac 12 \lambda(t_0-t)\}$. For $x\in \Omega_2$, we have $x-\rho(t)>\frac 12 x_0+\frac 12 \lambda(t_0-t)$ 
and thus
$$Q^3(x{-}\rho(t))\leq \frac C {(x_0+\lambda (t_0-t) )^6},\quad
\int_{\Omega_2} Q^3(x{-}\rho(t))\varphi'(\widetilde x) dx\leq 
 \frac C {(x_0+\lambda (t_0-t) )^6}.$$
Now, we claim
\begin{equation}\label{appendix1}
\int |\eta(t,x-\rho(t))|^3 \varphi'(\widetilde x)  dx \leq C \alpha_0 \int \eta^2(t,x-\rho(t)) \varphi'(\widetilde x) dx,
\end{equation}
where $C$ is independent of $A$.
See proof of \eqref{appendix1} in Appendix \ref{secAPPENDIX1}.
Moreover, as before, we find
\begin{equation*}\begin{split}
\int \eta^2(t,x-\rho(t)) \varphi'(\widetilde x) dx & \leq
C \int (u^2(t,x) + Q^2(x{-}\rho(t))) \varphi'(\widetilde x) dx
\\
&\leq C \int  u^2(t,x)  \varphi'(\widetilde x) dx +\frac C {(x_0+\lambda (t_0-t) )^2}.
\end{split}\end{equation*}
Thus, it follows from \eqref{step3}--\eqref{appendix1} that for $\alpha_0>0$ small enough,
$\forall t\in [t_1,t_0],$
\begin{equation}\label{finstep3}\begin{split}
M_{\varphi}'(t)  
&\leq -\ \frac 1{8} (1-\lambda) \int u^2(t) \varphi'(\widetilde x)
+C\alpha_0 \int u^2(t)\varphi'(\widetilde x)
+ \frac  C {(x_0+\lambda (t_0-t) )^2}
\\ &\leq -\ \frac 1{16} (1-\lambda) \int u^2(t) \varphi'(\widetilde x)
+ \frac  C {(x_0+\lambda (t_0-t) )^2}.
\end{split}
\end{equation}
Let $t\in [t_1,t_0]$. By integration of \eqref{finstep3} on $[t,t_0]$,
since
$$
\int_{t}^{t_0} \frac {dt'} {(x_0+\lambda (t_0-t') )^2}
=\frac {1}{\lambda x_0} \int_0^{\frac {\lambda (t_0-t)}{x_0}} 
\frac {dt''}{(1+t'')^2} \leq \frac  C {x_0},
\quad (t''=\frac \lambda {x_0} (t_0-t'))
$$
we find:
\begin{equation}\label{conclusion1}\begin{split}
&\int u^2(t_0,x)\varphi(x-x_0-\rho(t_0)) dx 
+\frac 1C\int_{t}^{t_0} \int u^2(t',x) \varphi'(x-x_0-\rho(t') -\lambda (t-t')) dxdt'
\\ &\leq \int u^2(t,x)\varphi(x-x_0-\rho(t)-\lambda (t_0-t)) dx +\frac C  {x_0}.
\end{split}\end{equation}
By density and continuous dependence (\cite{IK}) estimate \eqref{conclusion1} also holds for $H^{\frac 12}$ solutions.

\subsection{Monotonicity results for $\eta(t)$}
Here, we present similar monotonicity arguments for $\eta(t)$. See \cite{MMas2} for similar results
in the case of  the gKdV equations.

\begin{proposition}\label{MONOTONICITY5}
Let $0<\lambda<1$.
Under the assumptions of Lemma \ref{MODULATION}, 
for $\alpha_0$ small enough and $A$ large enough, there exists $C>0$ such that for all $x_0>1$,  $t_1\leq t_2$,
\begin{equation*}
\begin{split}
&\int \eta^2(t_2,x)(\varphi(x-x_0)-\varphi(-x_0))\, dx 
\\&    \leq \int \eta^2(t_1,x)(\varphi(x-\lambda (t_2-t_1)-x_0)
   -\varphi(-x_0-\lambda (t_2-t_1))) dx +  C \int_{t_1}^{t_2} \frac {\|\eta(t)\|_{L^2}^2}{
   (x_0+\lambda (t_2-t))^2} dt.
   \end{split}
\end{equation*}
\end{proposition}

\begin{remark} 
With respect to Proposition \ref{MONOTONICITY1}, we need to modify slighty the function in the integral
($\varphi(x-x_0)-\varphi(-x_0)$ instead of $\varphi(x-x_0)$) to remove some terms
in the second member, see comments in the proof. This estimate
is clearly improving Proposition \ref{MONOTONICITY1} since the remainder
term can now be controlled by $\frac {C}{x_0} \sup_{t} \|\eta(t)\|_{L^2}^2$.
\end{remark}

As for $u(t)$ in the proof of Proposition \ref{MONOTONICITY1}, 
we have by direct computations using \eqref{eqofeps},
\begin{align}&
\frac 12 \frac d{dt}\int \eta^2(t,x) \varphi(x) dx   =
\int \eta_t \eta \varphi(x) dx \nonumber\\ &
= - \int (\mathcal{L} \eta) (\eta\varphi'  + \eta_x \varphi )  
+ \frac 13 \int \eta^3 \varphi'   
+ (\rho'-1)\left( \int Q' \eta \varphi  -\tfrac 12  \int \eta^2 \varphi'\right) \nonumber\\ &
= \int (\mathcal{H} \eta_x ) \eta \varphi' + 
\int (\mathcal{H} \eta_x) \eta_x \varphi
-\frac 12 \int \eta^2 \varphi'
+\frac 12 \int \eta^2 (Q \varphi' - Q' \varphi)
+ \frac 13 \int \eta^3 \varphi' \nonumber \\& 
+ (\rho'-1)\left( \int Q' \eta \varphi  -\tfrac 12  \int \eta^2 \varphi'\right).\label{page10bis}
\end{align}
Let $0<\lambda<1$ and $\overline x= x-x_0  - \lambda (t_0-t)$.
Then, by Lemmas \ref{FIRSTTERM} and \ref{SECONDTERM}, we get
\begin{align*}
\frac d{dt} \int \eta^2 \varphi(\overline x) & \leq 
-\left( \rho'(t)-\lambda - \frac {2C_0}A\right) \int \eta^2 \varphi'(\overline x)
+ \int \eta^2 (Q \varphi'(\overline x) - Q' \varphi(\overline x) )
 + \frac 23 \int |\eta|^3 \varphi'(\overline x)  \\&
+ 2 (\rho'-1) \int  Q' \eta \varphi(\overline x).
\end{align*}
Now, as in the proof of Proposition \ref{MONOTONICITY1}, we fix
$A>1$ such that $\frac {2C_0}A \leq \frac 14 (1-\lambda)$
and $\alpha_0$ small enough so that
$\rho'-\lambda >\frac 12 (1-\lambda)$ by \eqref{ortho}. Then, 
by  \eqref{appendix1}
and \eqref{ortho}, we can choose $\alpha_0>0$ small enough so that
$$
\frac 23 \int |\eta|^3 \varphi'(\overline x)  
\leq \frac 18 (1-\lambda) \int \eta^2 \varphi'(\overline x).
$$
Thus, we obtain
\begin{equation*}
\frac d{dt} \int \eta^2 \varphi(\overline x)  \leq 
-\frac 18 \left( 1-\lambda \right) \int \eta^2 \varphi'(\overline x)
+ \int \eta^2 (Q \varphi'(\overline x) -  Q' \varphi(\overline x) )
+ 2 (\rho'-1) \int Q' \eta \varphi(\overline x) .
\end{equation*}
At this point, note that the term $\int \eta^2 Q' \varphi(\overline x)$ has no
sign, and since $\varphi(y)\sim \frac C{|y|}$ as $y\to -\infty$, this term can only
be controlled by $\frac C{(x_0+\lambda (t_0-t))} \int \eta^2$, which is not sufficient
for our purposes.
We modify slightly  the functional to cancel the main order of this
term.

\medskip

Indeed,  since $\int \eta Q'=0$, using \eqref{eqofeps}, we have
$$
\frac d {dt}\int \eta^2 = 2\int Q \eta \eta_x = - \int Q' \eta^2.
$$
Therefore, using also $\int Q'\eta=0$, we get
\begin{align*}
&\frac d{dt} \int \eta^2 \left(\varphi(\overline x)-\varphi(-x_0-\lambda (t_0-t))\right)
\leq -\frac 18 \left( 1-\lambda \right) \int \eta^2 \varphi'(\overline x)\\ &
+ \int \eta^2 \left(Q \varphi'(\overline x) -  Q' (\varphi(\overline x) -\varphi(-x_0-\lambda (t_0-t)))\right)
\\Ê&+ 2 (\rho'-1) \int \eta  Q' (\varphi(\overline x)-\varphi(-x_0-\lambda (t_0-t))) 
-\lambda \varphi'(-(x_0+\lambda (t_0-t))) \int \eta^2.
\end{align*}
Now, we claim the following estimate
\begin{equation}\label{decroissance}
\forall x\in \mathbb{R},\quad 
Q(x) \varphi'(\overline x)+\left| Q(x) \left(\varphi(\overline x)-\varphi(-(x_0+\lambda (t_0-t)))\right)\right|
\leq \frac C{(x_0+\lambda (t_0-t))^2}.
\end{equation}
Since 
$$
Q(x) \varphi'(\overline x)\leq
\frac C {(1+x^2)(1+(x-x_0-\lambda (t_0-t))^2)}
$$
(recall that the value of $A$ has been fixed)
estimate \eqref{decroissance} is clear for $Q(x) \varphi'(\overline x)$
 by considering the two regions
$|x|>\frac 12 (x_0+\lambda(t_0-t))$ and $|x|<\frac 12 (x_0+\lambda(t_0-t))$.

For the other term, we first note that since $|Q(x)|\leq \frac C{1+x^2}$
and $\varphi$ is
bounded,  the estimate is clear for $|x|>\frac 12 (x_0+\lambda(t_0-t))$.
For $|x|<\frac 12 (x_0+\lambda(t_0-t))$, we have
$$
|\varphi(\overline x)-\varphi(-x_0-\lambda (t_0-t))|
\leq |x| \sup_{[\frac 12 (x_0+\lambda (t_0-t),\frac 32 (x_0+\lambda (t_0-t)]} \varphi'
\leq \frac {C |x|}{ (x_0+\lambda(t_0-t))^2};
$$
thus, for such $x$, we obtain the following estimate which finishes the proof of
\eqref{decroissance}:
$$
\left| Q(x) \left(\varphi(\overline x)-\varphi(-x_0-\lambda (t_0-t))\right)\right|
\leq \frac C{(x_0+\lambda (t_0-t))^2}.
$$

\medskip

By \eqref{ortho} and
\eqref{decroissance}, and since $|Q'(x)|\leq \frac C{1+|x|} Q(x)$, we obtain
\begin{equation}\label{conseq1}
\left|\int \eta^2 (Q \varphi'(\overline x) -  Q' (\varphi(\overline x) -\varphi(-x_0-\lambda (t_0-t))))\right|\leq 
\frac {C \|\eta(t)\|_{L^2}^2} {(x_0+\lambda (t_0-t))^2},
\end{equation}
\begin{equation}\label{conseq2}\begin{split}
\left|(\rho'-1) \int Q' \eta (\varphi(\overline x)-\varphi(-x_0-\lambda (t_0-t))) \right|
&\leq \frac {C \|\eta(t)\|_{L^2}} {(x_0+\lambda (t_0-t))^2} \int \frac {|\eta|} {1+|x|} 
\\& \leq \frac {C \|\eta(t)\|_{L^2}^2} {(x_0+\lambda (t_0-t))^2}.
\end{split}
\end{equation}

The conclusion is thus:
\begin{align*}
&\frac d{dt} \int \eta^2 \left(\varphi(\overline x)-\varphi(-(x_0+\lambda (t_0-t)))\right)
\leq -\frac 18 \left( 1-\lambda \right) \int \eta^2 \varphi'(\overline x)
+\frac {C \|\eta(t)\|_{L^2}^2} {(x_0+\lambda (t_0-t))^2}.
\end{align*}
By integration on $[t,t_0]$, we get
\begin{equation*} \begin{split}
&\int \eta^2(t_0,x) \left(\varphi(x-x_0)-\varphi(-x_0)\right)dx
+ \frac 1C \int_t^{t_0} \int \eta^2(t',x) \varphi'(x-x_0-\lambda(t_0-t')) dxdt'\\
& \leq \int \eta^2(t,x) \left(\varphi( x-x_0-\lambda(t_0-t))-\varphi(-x_0-\lambda(t_0-t))\right) dx + 
C \int_t^{t_0} \frac { \|\eta(t')\|_{L^2}^2 dt'} {(x_0+\lambda (t_0-t'))^2}.
\end{split}
\end{equation*}

\section{Linear Liouville property}\label{sec:3}

In this section, we prove the following result.

\begin{theorem}\label{LINEARLIOUVILLE}
Let $w\in C(\mathbb{R},L^2(\mathbb{R}))\cap L^\infty(\mathbb{R},L^2(\mathbb{R}))$
be a solution of
\begin{equation}\label{lineareq}
w_t=(\mathcal{L} w)_x+\beta(t) Q',\quad (t,x)\in \mathbb{R}^2,
\quad \text{where $\beta$ is continuous,}
\end{equation}
satisfying
\begin{equation}\label{orthow}
\forall t\in \mathbb{R},\quad\int w(t,x)Q(x)dx=\int w(t,x)Q'(x)dx=0,
\end{equation}
\begin{equation}\label{decayw}
\forall t\in \mathbb{R},~\forall x_0>1,\quad
\int_{|x|>x_0} w^2(t,x) dx \leq \frac C{x_0}.
\end{equation}
Then
\begin{equation}\label{conclusionw}
w\equiv 0 \quad \text{on $\mathbb{R}^2$.}
\end{equation}
\end{theorem}

This result is similar to Theorem 3 in \cite{MMjmpa}. For the proof, we follow
the strategy of \cite{yvanSIAM}, \cite{MMas1}, introducing a dual problem
whose operator 
has better spectral properties. Since $w(t)$ is only $L^2$ and has a weak decay
at infinity in space, we will need to regularize and localize the dual solution.

For the sake of clarity, we now present the formal argument. The complete
justification will be  presented in Sections \ref{subsec:31} and \ref{secVIRIEL}.

Multiplying the equation of $w(t)$ by $x w(t)$, we get
$$
\frac d{dt} \int xw^2 = - 2 \int (\mathcal{H} w)w_x - \int w^2 
+\int w^2 (Q-xQ') + 2 \beta(t) \int x Q'w,
$$
where $(\int (Q')^2)\beta(t)=Ê\int w \mathcal{L}(Q'')$ (multiply the equation of $w$ by $Q'$ and use $\int w Q'=0$).
But it is not clear how to study the spectral properties of
the operator $$2 \int (\mathcal{H} w)w_x + \int w^2 
-\int w^2 (Q-xQ') +\frac 2{\int (Q')^2} \left(\int w \mathcal{L}Q''\right)\left( \int x Q'w\right).$$
Moreover, the decay estimate \eqref{decayw} is not quite enough to control
$\int x w^2$.

Therefore, we instead rely on the dual problem, setting $v=\mathcal{L} w$.
Since $\mathcal{L} Q'=0$ (direct calculation), we obtain the following equation for $v(t)$:
$v_t=\mathcal{L}(v_x)$. Multiplying the equation by $xv$, we obtain
$$
- \frac d{dt} \int xv^2 =  2 \int (\mathcal{H} v)v_x + \int v^2 
-\int v^2 (Q+xQ').
$$
Note that the operator in $v$ is much easier to study since now the potential $xQ'$ has a positive contribution ($xQ'\leq 0$), moreover, there is no scalar product.
In fact, we will obtain (see Proposition \ref{pureQUADRA}) the positivity of this operator
 under the orthogonality
condition $\int v (xQ)'=0$.
Observe that $\int v(xQ)'=\int (\mathcal{L} w) (xQ)'=- \int w Q=0$ since
$\mathcal{L}((xQ)')=-Q$ (see \eqref{surT}).

Provided that $\int |x|v^2(t)\leq C$, we would obtain from the above identity
$$
\int_{-\infty}^{+\infty} \|v(t)\|_{H^{\frac 12}}^2 dt \leq C,
$$
which says that for a subsequence $t_n\to +\infty$, $v(t_n)\to 0$, 
$w(t_n)\to 0$.
Combined with  energy conservation ($(\mathcal{L} w(t),w(t))=C$)
and Lemma \ref{76} below, this gives $w\equiv 0$.
But \eqref{decayw} is not enough to obtain the estimate $\int |x|v^2(t)\leq C$
In fact, since $w(t)$ is only in $L^2$, we both need
to  localize and regularize the dual problem.

\subsection{Proof of Theorem \ref{LINEARLIOUVILLE} assuming positivity of a quadratic form}\label{subsec:31}

\begin{lemma}[Regularized dual problem]\label{DUAL}
There exists $\gamma_0>0$ such that for any $0<\gamma<\gamma_0$, the
following is true.
Let $v=(1-\gamma \partial_x^2)^{-1}(\mathcal{L} w)$.
Then, $v\in C(\mathbb{R},H^1(\mathbb{R}))\cap L^\infty(\mathbb{R}, H^1(\mathbb{R}))$ and \begin{enumerate}
\item Equation of $v$. 
\begin{equation}\label{eqofv}
v_t 
=\mathcal{L}(v_x)  - \gamma (1-\gamma \partial_x^2)^{-1}(2 v_{xx} Q' + v_x Q'') .
\end{equation}
\item Decay of $v$.
\begin{equation}\label{lemmdec}
\forall t\in \mathbb{R},x_0>1,\quad
\int_{|x|>x_0}  (v_x^2(t,x)+v^2(t,x)) dx \leq \frac {C_{\gamma}} {x_0^{\frac 34}}.
\end{equation}
\item Virial type estimate.
\begin{equation}\label{finalviriel}
\int_{-\infty}^{+\infty}
\frac 1{(1+t^2)^{\frac 25}} \|v(t)\|_{H^1}^2 dt <C.
\end{equation}
\end{enumerate}
\end{lemma}
\noindent\emph{Proof of Lemma \ref{DUAL}.}
First, since
 $\sup_t \|w(t)\|_{L^2} \leq C$, we obtain
 $\sup_t \|v(t)\|_{H^1}\leq C_{\gamma}$ (see Claim~\ref{TECH} below).

\medskip

\noindent\emph{ 1. Equation of $v$.}
Let $\widetilde v=\mathcal{L} w$ so that $w_t=\widetilde v_x + \beta Q'$. Since $\mathcal{L} Q'=0$, 
the function $\widetilde v$ satisfies $\widetilde v_t = \mathcal{L} w_t= \mathcal{L}(\widetilde v_x)$.
Now, we introduce a regularization of the function $\widetilde v$. For $0<\gamma<\frac 12$ to be chosen
later small enough, we set:
\begin{equation}\label{defofv}
v(t,x)=(1-\gamma \partial_x^2)^{-1} \widetilde v(t,x) \quad
\text{or equivalently}\quad
v-\gamma v_{xx} = \widetilde v = \mathcal{L} w.
\end{equation}
Then, $v(t,x)$ satisfies the following equation
\begin{equation*}
v_t=(1-\gamma \partial_x^2)^{-1} \widetilde v_t = (1-\gamma \partial_x^2)^{-1} \mathcal{L}(\widetilde v_x)
=\mathcal{L}(v_x) - (1-\gamma \partial_x^2)^{-1} (\widetilde v_x Q) + v_x Q.
\end{equation*}
But
$-(1-\gamma \partial_x^2)^{-1} (\widetilde v_x Q) + v_x Q
=(1-\gamma \partial_x^2)^{-1}(-2 \gamma v_{xx} Q' - \gamma v_x Q'')$, and so
\begin{equation}
v_t 
=\mathcal{L}(v_x)  - \gamma (1-\gamma \partial_x^2)^{-1}(2 v_{xx} Q' + v_x Q'') .
\end{equation}

\smallskip

\noindent\emph{2. Decay estimate on $v$.}
By using the decay on $w(t)$, we claim
\begin{equation}\label{decayv}
\forall x_0>1,\forall t,\quad 
\int_{|x|\geq x_0} (v_x^2(t,x)+v^2(t,x))dx \leq \frac {C_{\gamma}}
{x_0^{\frac 34}}.
\end{equation}
Indeed, let ($x_0>1$) 
$$h(x)=h_{x_0}(x)=
\varphi_{\sqrt{x_0}}^2(x-x_0)=
\left(\frac \pi 2 + \arctan\Big(\frac {x-x_0}{\sqrt{x_0}}\Big)\right)^2.
$$ 
Note that $0\leq |h'|+|h''|\leq C h$.
Since $v-\gamma v_{xx}=\mathcal{L} w$, multiplying by $vh$, we have
\begin{equation}\label{brrr}
 \int v^2 h+\gamma \int v_x^2h-\frac 12 \gamma \int v^2 h''
=\int w \mathcal{L} (vh)
=\int w D (vh) + \int wvh - \int Q wvh.
\end{equation}

First, from 
$$
\left| \int wvh\right| + \left|  \int Q wvh\right|\leq C \|w \sqrt{h}\|_{L^2} \|v \sqrt{h}\|_{L^2}
$$
and
$\int w^2 h\leq \int_{x<\frac {x_0}2} w^2 h + \int_{x>\frac {x_0}2} w^2
\leq C \frac 1 {x_0}$
(using the definition of $h$ and \eqref{decayw}) it follows that
$$
\left| \int wvh \right| + \left| \int Q wvh\right|
\leq \frac C{x_0^{\frac 12}} \|v \sqrt{h}\|_{L^2}.
$$

Second, by Lemma \ref{COMMUTATOR}, we have
\begin{align*}
& \left| \int w D(vh) - \int D(v\sqrt{h}) \sqrt{h} w\right|
\leq \|w\|_{L^2} \|v \sqrt{h}\|_{L^4} \|D(\sqrt{h})\|_{L^4} \leq
 \frac {C}{x_0^{\frac 38}}  \|v \sqrt{h}\|_{H^1}.
\end{align*}

Since 
$\left| \int  D(v\sqrt{h}) \sqrt{h} w\right|\leq \|w \sqrt{h}\|_{L^2}
\|v \sqrt{h}\|_{H^1}$ and 
$\|v \sqrt{h}\|_{H^1}^2 = \int (v_x^2  + v^2) h + O(\frac 1{x_0})$,
we obtain from \eqref{brrr}
$$
\int_{x>x_0} (v_x^2+v^2)\leq 
\int (v_x^2  + v^2) h \leq \frac { C_{\gamma}} {x_0^{\frac 34}}.
$$

\smallskip

\noindent\emph{3. Virial type estimate on $v(t)$.}
Let $\frac 13 < \theta <\frac 12$, $B>1$ to be chosen later and set
$$
I(t)=\frac 12 \int g\bigg(\frac {x} {(B+t^2)^\theta}\bigg) v^2(t,x) dx,
\quad z=v \sqrt{g' \bigg(\frac {x} {(B+t^2)^\theta}\bigg)}\quad \text{where $g(x)=\arctan(x)$,}
$$
\begin{equation}\label{defLtilde}
(\widetilde{\mathcal{L}}z,z)=-2(\mathcal{L} (z_x),xz)=2 \int |D^{\frac 12} z|^2 + \int z^2 - \int(xQ'+Q) z^2.
\end{equation}
For any $0<\sigma_0<1$, we claim
\begin{equation}\label{examen}
\begin{split}
& \left|2 I'(t) + \frac 1{(B+t^2)^\theta} (\widetilde{\mathcal{L}}z,z)\right | \leq 
\frac {\sigma_0}{(B+t^2)^\theta} \|z\|_{L^2}^2+ \frac C{\sigma_0(B+t^2)^{1-\theta}} \|v\|_{L^2}^2\\
& + \frac {C}{(B+t^2)^{\frac 74\theta}}  \|z\|_{H^{\frac 12}}\|v\|_{L^2}
+ \frac {C}{(B+t^2)^\theta}  \gamma^{\frac 14}  \|z\|_{H^{\frac 12}}\|v\|_{L^2}
+\frac {C}{(B+t^2)^{2\theta}} \|z\|_{L^2}^2.
\end{split}\end{equation}

\noindent\emph{Proof of \eqref{examen}.}
We compute $I'(t)$:
\begin{align*}
I'(t)=-\frac {\theta t}{(B+t^2)^{\theta+1}} \int x g'\bigg(\frac {x} {(B+t^2)^\theta}\bigg) v^2
+\int g\bigg(\frac {x} {(B+t^2)^\theta}\bigg)  v v_t.
\end{align*}
First, note that by Cauchy-Schwarz' inequality, for any $\sigma_0>0$,
\begin{align*}
  \left|\frac {\theta t}{(B+t^2)^{\theta+1}} \int x g'\bigg(\frac {x} {(B+t^2)^\theta}\bigg) v^2\right|
&\leq \frac {\sigma_0}{(B+t^2)^\theta} \int   g'\bigg(\frac {x} {(B+t^2)^\theta}\bigg) v^2
\\& + \frac {\theta^2 t^2} {4 \sigma_0 (B+t^2)^{2-\theta}} \int \bigg(\frac {x}{(B+t^2)^\theta}\bigg)^2 g'\bigg(\frac {x} {(B+t^2)^\theta}\bigg) 
v^2.
\end{align*}
Since $s^2 g'(s)\leq 1$, we obtain
\begin{align*}
& \left|\frac {\theta t}{(B+t^2)^{\theta+1}} \int x g'\bigg(\frac {x} {(B+t^2)^\theta}\bigg) v^2\right|
\leq \frac {\sigma_0}{(B+t^2)^\theta} \int   z^2
+ \frac {C\theta^2} { \sigma_0 (B+t^2)^{1-\theta}} \int v^2.
\end{align*}
Second, we use the equation of $v$ to compute the term
$\int g\big(\frac {x} {(B+t^2)^\theta}\big)  v v_t$.
\begin{align*}
&\int g\bigg(\frac {x} {(B+t^2)^\theta}\bigg)  v v_t
 \\& = \int g\bigg(\frac {x} {(B+t^2)^\theta}\bigg)  v \mathcal{L}v_x
-\gamma \int g\bigg(\frac {x} {(B+t^2)^\theta}\bigg) v (1-\gamma \partial_x^2)^{-1}
(2v_{xx} Q' + v_x Q'') = \mathbf{A}+\mathbf{B}.
\end{align*}
\textit{Estimate on $\mathbf{A}$.}
\begin{align*}
\mathbf{A}&=
\int g\bigg(\frac {x} {(B+t^2)^\theta}\bigg) v (-\mathcal{H} v_{xx} + v_x)
-\int g\bigg(\frac {x} {(B+t^2)^\theta}\bigg)  Q vv_x\\
&=-\int  \bigg(\frac 1 {(B+t^2)^\theta}  g'\bigg(\frac {x} {(B+t^2)^\theta}\bigg)
v + g \bigg(\frac {x} {(B+t^2)^\theta}\bigg) v_x\bigg) (-\mathcal{H} v_{x} + v)
\\ &+ \frac 12 \int \bigg(\frac 1 {(B+t^2)^\theta}  g'\bigg(\frac {x} {(B+t^2)^\theta}\bigg)Q+g \bigg(\frac {x} {(B+t^2)^\theta}\bigg) Q'\bigg) v^2.
\end{align*}
Next,
\begin{align*}
\mathbf{A}&=- \frac 1 {(B+t^2)^\theta} \int |D^{\frac 12}z|^2
+  v \Big(D \big(vg'\big(\tfrac {x} {(B+t^2)^\theta} \big)
-D \big(v\sqrt{g' \big(\tfrac {x} {(B+t^2)^\theta}\big)}\big)\sqrt{g' \big(\tfrac {x} {(B+t^2)^\theta}\big)}\Big)\\
& + \int (\mathcal{H} v_x) v_x  g \bigg(\frac {x} {(B+t^2)^\theta}\bigg) 
-\frac 12 \frac 1 {(B+t^2)^\theta} \int z^2 \\
&+   \frac 12 \frac 1 {(B+t^2)^\theta} \int (xQ'+Q)z^2
+\frac 12  \int \bigg(g \bigg(\frac {x} {(B+t^2)^\theta}\bigg)-\frac {x} {(B+t^2)^\theta}
g'\bigg(\frac {x} {(B+t^2)^\theta}\bigg)\bigg) Q'v^2\\
&=-\frac 12 \frac 1 {(B+t^2)^\theta} (\widetilde{\mathcal{L}}z,z) +\mathbf{A_1}+\mathbf{A_2}+\mathbf{A_3},
\end{align*}
where
$$
\mathbf{A_1}=
- \frac 1 {(B+t^2)^\theta}\int v \Big(D \big(z\sqrt{g'\big(\tfrac {x} {(B+t^2)^\theta}} \big)
\big)
-(D z) \sqrt{g' \big(\tfrac {x} {(B+t^2)^\theta}\big)}\big),
$$ $$
\mathbf{A_2}=
\int (\mathcal{H} v_x) v_x  g \bigg(\frac {x} {(B+t^2)^\theta}\bigg),$$
$$\mathbf{A_3}=
\frac 12  \int \bigg(g \bigg(\frac {x} {(B+t^2)^\theta}\bigg)-\frac {x} {(B+t^2)^\theta}g' \bigg(\frac {x} {(B+t^2)^\theta}\bigg)\bigg) Q'v^2.
$$ 

\noindent \textit{Estimate on $\mathbf{A_1}$.}
By Lemma \ref{COMMUTATOR}, we have
$$
|\mathbf{A}_1|\leq
\frac C{(B+t^2)^\theta} \|v\|_{L^2} \|z\|_{L^4} \big\|D\sqrt{g'
\big(\tfrac {x} {(B+t^2)^\theta}\big)}\big\|_{L^4} 
\leq \frac C{(B+t^2)^{\frac {7\theta}4 }} \|v\|_{L^2} \|z\|_{H^{\frac 12}}.
$$

\noindent\textit{Estimate  on $\mathbf{A_2}$.}
Since $\int (\mathcal{H}v_x)v_x=0$, Lemma \ref{SECONDTERM}, applied to $A=(B+t^2)^{\theta}$
gives
$$
|\mathbf{A}_2|\leq \frac C{(B+t^2)^{2\theta}}
\|z\|_{L^2}^2.
$$

\noindent\textit{Estimate on $\mathbf{A_3}$.}
Since for all $y\in \mathbb{R}$, $|\arctan y - \frac {y}{1+y^2}|\leq C y^2$,
we have, for all $x\in \mathbb{R}$,
$$
\left|\bigg(g\bigg(\frac {x} {(B+t^2)^\theta}\bigg)-\frac {x} {(B+t^2)^\theta}
g'\bigg(\frac {x} {(B+t^2)^\theta}\bigg)\bigg) Q'(x)
\right|\leq \frac {x^2 |Q'(x)|}{(B+t^2)^{2\theta}}\leq
\frac C{(B+t^2)^{2\theta}}\frac 1{1+|x|}.
$$
Thus,
$$
|\mathbf{A}_3|\leq \frac C{(B+t^2)^{2\theta}}
\|v\|_{L^2}\|z\|_{L^2}.
$$

\noindent\textit{Estimate on $\mathbf{B}$.}
First, we claim the following.

\begin{claim}\label{TECH}
\begin{description}
\item{\rm (i)}
$x(1-\gamma \partial_x^2 )^{-1} f =
(1-\gamma \partial_x^2 )^{-1} (xf)- 2 \gamma (1-\gamma \partial_x^2)^{-2} (f')$.
\item{\rm (ii)}
$\|(1-\gamma \partial_x^2 )^{-1}f\|_{L^2}
+ \gamma^{\frac 12}\|(1-\gamma \partial_x^2 )^{-1}(f')\|_{L^2}
+ \gamma\|(1-\gamma \partial_x^2 )^{-1}(f'')\|_{L^2}
\leq C \|f\|_{L^2},$

$\|(1-\gamma \partial_x^2 )^{-1}(f'')\|_{L^2}
\leq C \gamma^{- \frac 34}\|f\|_{\dot H^{\frac 12}}.$
\end{description}
\end{claim}

\noindent\emph{Proof of Claim \ref{TECH}.}
(i) Let $h=(1-\gamma \partial_x^2 )^{-1} f$. Then,
$xh-\gamma (xh)''=xf - 2 \gamma h'$ and so
$xh = (1-\gamma \partial_x^2 )^{-1} (xf-2 \gamma(1-\gamma \partial_x^2 )^{-1}f')$.

(ii) $\int |f|^2= \int |h-\gamma h''|^2 
=\int h^2 + 2 \gamma \int (h')^2 + \gamma^2 \int (h'')^2$, which proves the first estimate.

Next, 
$\|(1-\gamma \partial_x^2 )^{-1}f''\|_{L^2}
\leq C \|(\frac {\xi^2}{1+\gamma \xi^2}) \hat f\|_{L^2}
\leq   C {\gamma^{-\frac 34}} \| |\xi|^{\frac 12} \hat f\|_{L^2}$,
since 
$\forall \xi\in \mathbb{R}$, $\forall \gamma>0$,
$\frac {\xi^2}{1+\gamma \xi^2} \leq   {\gamma^{-\frac 34}} {|\xi|^{\frac 12}}$.
The claim is proved.

\medskip

Using (i) of Claim \ref{TECH}, we obtain
$$\mathbf{B}=-\gamma \int \frac 1x g\big(\frac {x} {(B+t^2)^\theta}\big) v \, (1-\gamma \partial_x^2)^{-1} H,$$
where
$$H=2 x v_{xx} Q' + xv_x Q'' - 2 \gamma (1-\gamma \partial_x^2)^{-1}( 2 v_{xx} Q' + v_x Q'')_x.$$
Since $|g(y)|\leq C |y|$, for all $y$, we have
$|\mathbf{B}|\leq \frac {C\gamma}{(B+t^2)^{\theta}} \|v\|_{L^2}\| (1-\gamma \partial_x^2)^{-1}H\|_{L^2}$.
Now, we use Claim \ref{TECH} (ii) to estimate $\| (1-\gamma \partial_x^2)^{-1}H\|_{L^2}$.
We can rewrite $H$ under the form:
$$
H=(2 v x Q')''+ (vF_1)'+ vF_2- 2 \gamma (1-\gamma \partial_x^2)^{-1} ((2vQ')''+(vF_3)'+vF_4)_x,
$$
where for $j=1,\ldots,4$, $|F_j(x)|\leq C \frac 1{1+x^2}$. Thus,
\begin{align*}
&\|(1-\gamma \partial^2)^{-1}H\|_{L^2}\\
&\leq C \gamma^{-\frac 34} \| v x Q'\|_{\dot H^{\frac 12}}
+ C \gamma^{-\frac 12} \|v\tfrac 1{1+x^2}\|_{L^2} + \gamma^{\frac 12}
\|(1-\gamma \partial_x^2)^{-1}((2vQ')''+(vF_3)'+vF_4)\|_{L^2}\\
&\leq C \gamma^{-\frac 34} \| v x Q'\|_{\dot H^{\frac 12}}
+ C \gamma^{-\frac 12} \|v\tfrac 1{1+x^2}\|_{L^2}.
\end{align*}
Now, we claim
\begin{equation}\label{compvz}
\|v \tfrac 1{1+x^2}\|_{L^2}\leq C \|z\|_{L^2},\quad
\|v x Q'\|_{\dot H^{\frac 12}}\leq C \|z\|_{H^{\frac 12}}.
\end{equation}
The first estimate is clear since
$\tfrac 1{1+x^2}\leq C \sqrt{g'}$. 
Let $f(x)=\frac {xQ'(x)}{\sqrt{g'(\frac {x}{(B+t^2)^{\theta}})}}$. Then,
by Lemma \ref{COMMUTATOR},
$$
\|D^{\frac 12} (vxQ')\|_{L^2}
=\|D^{\frac 12} (zf)\|_{L^2}
\leq \|(D^{\frac 12} z) f\|_{L^2} + C \|z\|_{L^4} \|D^{\frac 12}f\|_{L^4}
\leq C \|z\|_{H^{\frac 12}},
$$
since $\|f\|_{L^\infty}+\|D^{\frac 12} f\|_{L^4} \leq \|f\|_{H^1}  \leq C.$

Thus,
$\|(1-\gamma \partial^2)^{-1}H\|_{L^2}\leq C \gamma^{-\frac 34} \|z\|_{H^{\frac 12}}$
and
 in conclusion for the term $\mathbf{B}$:
$$
|\mathbf{B}|\leq \frac {C \gamma^{\frac 14}}{(B+t^2)^{\theta}} \|z\|_{H^{\frac 12}}
\|v\|_{L^2}.
$$
Putting together the above estimates, we obtain \eqref{examen}.

\medskip

We now claim the following (see proof in Section \ref{secVIRIEL}):

\begin{proposition}\label{QUADRA}
There exist $\lambda>0$, $\gamma_0>0$ and $B_0>1$ such that, 
 for $0<\gamma <\gamma_0$, $B\geq B_0$,
$$\forall t,\quad 
(\widetilde{\mathcal{L}}z(t),z(t))\geq \lambda \|z(t)\|_{H^{\frac 12}}^2,\quad
\text{where $z$ is as above}.
$$
\end{proposition}

\begin{remark} The operator $\widetilde{\mathcal{L}}$ does not depend on 
$\gamma$ and $B$, but the   orthogonality conditions on $w$ imply almost orthogonality conditions on
$z$ that depend on $\gamma$, $B$, see proof of Proposition~\ref{QUADRA}.
\end{remark}

Choose $\theta=\frac 25$
and fix $\sigma_0=\frac \lambda 4$.
Then,
$$
-2 I'(t)\geq \frac {\lambda}{2(B+t^2)^\theta} \|z(t)\|_{H^{\frac 12}}^2
-\frac {C}{(B+t^2)^\theta}\bigg(\frac 1{(B+t^2)^{\frac 15}}   
+  \gamma^{\frac 12} \bigg) \|v\|^2_{L^2}.
$$
By the decay property \eqref{lemmdec},
$$
\int v^2(t)\leq \int_{|x|\leq \frac 12 (B+t^2)^\theta} v^2(t)
+ \frac {C_{\gamma}} {(B+t^2)^{\frac 34\theta}}
\leq C \int z^2(t)+\frac {C_{\gamma}} {(B+t^2)^{\frac 3{10}}}.
$$
For $\gamma>0$ small enough and $B$ large enough,
and by $\|v\|_{H^{\frac 12}}\leq C$,  we get
$$
-2 I'(t)\geq \frac {\lambda}{4(B+t^2)^\theta} \|z(t)\|_{H^{\frac 12}}^2
-\frac {C_{\gamma}}{(B+t^2)^{\frac 35}}.
$$
Since $I(t)$ is bounded, we obtain by integration
\begin{equation}\label{mach}
\int_{-\infty}^{+\infty} \frac {1}{(B+t^2)^\theta} \|z(t)\|_{H^{\frac 12}}^2 dt
<C_{\gamma}.
\end{equation}
We claim that \eqref{mach} and \eqref{lemmdec} imply
\begin{equation}\label{survvv}
\int_{-\infty}^{+\infty} \frac {1}{(B+t^2)^\theta} \|v(t)\|_{H^{\frac 12}}^2 dt
<C.
\end{equation}
Indeed,
by \eqref{lemmdec} and the expression of $g'$,
and considering the two regions
$x>\frac 1{(B+t^2)^{\frac \theta 2}}$, $x<\frac 1{(B+t^2)^{\frac \theta 2}}$,
we have
\begin{equation}\label{truc}
\|v-z\|_{H^1}^2=
\|v(1-\sqrt{g'})\|_{H^1}^2 \leq  \frac C{(B+t^2)^{\frac 3{8} \theta}}
=\frac C{(B+t^2)^{\frac 3{20}}}.
\end{equation}
Thus, by $\|v\|_{H^{\frac 12}}\leq \|z\|_{H^{\frac 12}} + \|v-z\|_{H^{\frac 12}}$,
and \eqref{mach}
$$
\int_{-\infty}^{+\infty} \frac {1}{(B+t^2)^\theta} \|v(t)\|_{H^{\frac 12}}^2 dt
\leq 2 C_\gamma+\int_{-\infty}^{+\infty} \frac {1}{(B+t^2)^{\frac {11}{20}}} dt
\leq C.
$$

Using another virial argument, we claim
\begin{equation}\label{otherviriel}
\int_{-\infty}^{+\infty} \frac {1}{(B+t^2)^\theta} \|z(t)\|_{H^{\frac 32}}^2 dt
<C.
\end{equation}
Proof of \eqref{otherviriel}.
We set
$$
J(t)=\frac 12 \int g\bigg(\frac {x}{(B+t^2)^\theta}\bigg) v_x^2(t).
$$
Proceeding as in the proof of \eqref{examen} (the equation for $v_x$
is very similar to the one for $v$), we obtain
$$
\left| J'(t) + \frac 1 {(B+t^2)^\theta} \int (D^{\frac 32}z)^2\right| \leq
\frac C{(B+t^2)^\theta} \|v\|_{H^1}
(\|v\|_{H^1}+ \|z\|_{H^{\frac 32}}).
$$
Using 
$
\|v\|_{H^1} \leq \|z\|_{H^1} + \|v-z\|_{H^1},
$ \eqref{truc} 
and the following estimate 
$$
\|z\|_{H^1} \leq \varepsilon \|D^{\frac 32} z\|_{L^2}
+ C_{\varepsilon} \|z\|_{L^2},
$$
we obtain, for $\varepsilon>0$ small enough,
$$
-J'(t)\geq \frac 12 \frac 1 {(B+t^2)^\theta} \|D^{\frac 32}z\|_{L^2}^2
+ C \frac 1 {(B+t^2)^\theta} \|z\|_{L^2}^2.
$$
Since $J(t)$ is bounded and using \eqref{mach},
we obtain \eqref{otherviriel}.
 
Finally, by    \eqref{mach}, \eqref{truc} and \eqref{otherviriel}, we get  \eqref{finalviriel}. Lemma \ref{DUAL} is proved.

\begin{lemma}[Decay estimate on $w(t)$]\label{BACKTOw}
The following hold
\begin{equation}\label{firststa}
\int_{-\infty}^{+\infty} \frac {1}{(1+t^2)^{\frac 25}} \|w(t)\|_{L^2}^2 dt
<C,
\end{equation}
\begin{equation}\label{secondsta}
\sup_{t\in \mathbb{R}} \int |x| w^2(t,x) dx \leq C.
\end{equation}
\end{lemma}
\noindent\emph{Proof of Lemma \ref{BACKTOw}.}
Estimate \eqref{firststa} is a consequence of Lemma \ref{DUAL} by comparing $v$ and $w$.
Let $\gamma>0$ small. We have by the definition of $v$:
$ 
(1-\gamma \partial_x^2)  v = \mathcal{L} w.
$ 
Let $\widetilde w= (1-\gamma \partial_x^2)^{-\frac 14} w$.
Then,
$$
\int  w (1-\gamma \partial_x^2)^{\frac 12}  v = 
\int  \widetilde w (1-\gamma \partial_x^2)^{-\frac 14}(\mathcal{L} w).
$$
On the one hand, we have
\begin{equation*} 
\left|\int   w (1-\gamma \partial_x^2)^{\frac 12} v\right|
\leq C \|w\|_{L^2} \|v\|_{H^{1}}.
\end{equation*}
On the other hand, as in the proof of Claim \ref{TECH}
$$
\|(1-\gamma \partial_x^2)^{-\frac 14}(\mathcal{L} w)
- \mathcal{L} \widetilde w\|_{L^2}
\leq  \|(1 - \gamma \partial_x^2)^{-\frac 14}(Q w) - Q w\|_{L^2}
+\| Q(w-\widetilde w)\|_{L^2}
\leq \gamma^{\frac 14} \|w\|_{L^2}.
$$
Thus,
$$
\left|\int  \widetilde w (1-\gamma \partial_x^2)^{-\frac 14}(\mathcal{L} w)
- (\mathcal{L} \widetilde w,\widetilde w)\right|
\leq C \gamma^{\frac 14} \| w\|_{L^2}^2
$$
and since $(\mathcal{L} \widetilde w, \widetilde w) \geq \frac 12 \lambda \|\widetilde w\|_{H^{\frac 12}}$
for $\gamma>0$ small enough
(this is a consequence of Lemma \ref{76} and the orthogonality
conditions on $w$ -- see Section \ref{secVIRIEL},
in particular the proof of Proposition~\ref{QUADRA}), we obtain
$$
\int  \widetilde w (1-\gamma \partial_x^2)^{-\frac 14}(\mathcal{L} w)
\geq \frac \lambda 2 \|\widetilde w\|_{H^{\frac 12}}^2 - C  \gamma^{\frac 14} \| w\|_{L^2}^2
\geq \lambda_1 \|  w\|_{L^2}^2.
$$
In conclusion, we have obtained
$$
\|w\|_{L^2}\leq C \|v\|_{H^1},
$$
and Lemma \ref{DUAL} then implies \eqref{firststa}.

\medskip

Now, we prove \eqref{secondsta}.
Indeed, the integrability property \eqref{firststa} allows us to obtain the decay
on $w(t,x)$ by monotonicity properties.

By the proof of Proposition \ref{MONOTONICITY5}, we have, for any
$\lambda\in (0,1)$, for any $t_0$, $t\in (-\infty,t_0]$, $x_0>1$,
\begin{equation}\label{bonjbis}\begin{split}
&\int w^2(t_0,x) \left(\varphi(x{-}x_0)-\varphi(-x_0)\right)dx
\\&\leq \int w^2(t,x) \left(\varphi( x{-}x_0{-}\lambda(t_0{-}t))-\varphi(-x_0-\lambda(t_0{-}t))\right) dx + 
C \int_t^{t_0} \frac { \|w(t')\|_{L^2}^2 dt'} {(x_0+\lambda (t_0-t'))^2}.
\end{split}
\end{equation}
The last term in \eqref{bonjbis} is treated as follows $(x_0>1)$
$$
\int_t^{t_0} \frac { \|w(t')\|_{L^2}^2 dt'} {(x_0+\lambda (t_0-t'))^2}
\leq 
C x_0^{-\frac 65}
\int_{-\infty}^{+\infty} \frac { \|w(t')\|_{L^2}^2 dt'} {(1+(t_0-t'))^{\frac 45}}
$$
Thus, by \eqref{firststa}
(applied to $w(t+t_0)$) and \eqref{decayw},
 letting $t\to -\infty$ in \eqref{bonjbis}, we obtain
$$
\int w^2(t_0)\left(\varphi(x-x_0)-\varphi(-x_0)\right)dx
\leq C{x_0^{-\frac 65}}.
$$
By the change of  variable $x\to -x$, $t\to -t$, which leaves the equation invariant,
we get:
$$
\int w^2(t_0)\left(\varphi(x_0)-\varphi(x+x_0)\right)dx
\leq C{x_0^{-\frac 65}},
$$
and thus, summing up the two estimates,
\begin{equation*}
\int w^2(t_0) \left(\varphi(x-x_0)-\varphi(x+x_0)+\varphi(x_0)-\varphi(-x_0)\right)dx
\leq C{x_0^{-\frac 65}}.
\end{equation*}
We verify easily that for all $|x|>x_0\geq 1$,
\begin{equation}\label{calcul}\begin{split}
\varphi(x-x_0)-\varphi(x+x_0)+\varphi(x_0)-\varphi(-x_0)& \geq \varphi(0)- \varphi(2x_0)+\varphi(x_0)-\varphi(-x_0)\\
& \geq \tfrac \pi 2 -\arctan(2)>0.
\end{split}
\end{equation}
Thus, for all $x_0>1$,
\begin{equation}\label{voiture}
\int_{|x|\geq x_0} w^2(t_0)\leq C{x_0^{-\frac 65}}.
\end{equation}
By integrating in $x_0$, we obtain the following estimate
\begin{equation}\label{uniform}
\forall t\in \mathbb{R},\quad \int |x| w^2(t)\leq C.
\end{equation}
Thus Lemma \ref{BACKTOw} is proved.

\medskip

Now, we claim that estimate \eqref{secondsta} implies a gain of regularity
on $w(t)$.

\begin{lemma}[Gain of regularity on $w(t)$]\label{RIGUEUR}
Let $w\in C(\mathbb{R},L^2(\mathbb{R}))\cap L^\infty(\mathbb{R},L^2(\mathbb{R}))$
be a solution of \eqref{lineareq} 
satisfying \eqref{secondsta}. Then, $w(t)\in  C(\mathbb{R},H^{\frac 12}(\mathbb{R}))$ and the following
identity holds
\begin{equation}\label{virielw}
\int x w^2(t_2)- \int x w^2(t_1)
= -\int_{t_1}^{t_2} \int \big(2 |D^{\frac 12} w|^2+w^2 + w^2 (xQ' -  Q)\big)
+2 \int_{t_1}^{t_2} \beta(t) \int x Q ' w.
\end{equation}
\end{lemma}

\noindent\emph{End of  the proof of Theorem \ref{LINEARLIOUVILLE} assuming Lemma \ref{RIGUEUR}.}\quad 
Note first that multiplying the equation of $w(t)$ by $Q'$ and using $\int w Q'=0$,
we find $\left(\int (Q')^2 \right) \beta(t)= \int w \mathcal{L}(Q'')$, so that
\begin{equation}\label{controlbeta}
|\beta(t)|\leq C \|w\|_{L^2}.
\end{equation}
Multiplying the equation of $w(t)$ by
$\mathcal{L} w$ and using $\mathcal{L} Q'=0$, we also have
$$
\forall t\in \mathbb{R}, Ê\quad
(\mathcal{L}w(t),w(t))=(\mathcal{L}w(0),w(0)).
$$

By \eqref{virielw}, 
the estimates on $\int |x| w^2(t)$ and on $\beta(t)$, and Lemma \ref{BACKTOw},  we have
$$
\int_{-\infty}^{+\infty} \frac 1{(1+t^2)^{\frac 25}}
\|w(t)\|_{H^{\frac 12}}^2 dt<C.
$$
This implies that for a sequence $t_n\to +\infty$, we have
$\|w(t_n)\|_{H^{\frac 12}}\to 0$ as $n\to +\infty$.

Since $(\mathcal{L}w(t),w(t))=\lim_{t_n \to \infty}(\mathcal{L}w(t_n),w(t_n))$,
we obtain $(\mathcal{L}w(t),w(t))=0$ and so by the orthogonality conditions on
$w(t)$ and Lemma \ref{76}, 
we finally obtain $\forall t$, $w(t)=0$.

\bigskip

\noindent\emph{Proof of Lemma \ref{RIGUEUR}.} Formally, identity \eqref{virielw} follows from multiplying equation \eqref{lineareq}
by $xw$, integration by parts and properties of the Hilbert transform.
To justify \eqref{virielw}, we use a regularization of $w(t)$.

We set
$w_n=(1-\frac 1n \partial_x^2)^{-1} w$, so 
that for all $t$, $w_n(t)\to w(t)$ in $L^2(\mathbb{R})$
as $n\to +\infty$.
Then, $w_n$ satisfies the following equation
\begin{equation}\label{eqwn}
w_{nt}=(\mathcal{L}w_n)_x - \tfrac 1n (1-\tfrac 1n \partial_x^2)^{-1}(2Q'w_{nx}+w_n Q'')_x
+ \beta (1-\tfrac 1n \partial_x^2)^{-1}Q'.
\end{equation}
Let $h:\mathbb{R}\to \mathbb{R}$ be a smooth nondecreasing function such that
$h(x)=x$ if $x>1$ and $h(x)=0$ if $x<0$.
Then,
\begin{equation*}\begin{split}
\int h(x) w^2& =
\int h(x) (w_n-\tfrac 1n w_{nxx})^2
= \int h(x) w_n^2 - \tfrac 2n \int w_{nxx} w_n h(x) + \frac 1{n^2} \int w_{nxx}^2 h(x)
\\
&=\int h(x) w_n^2 + \tfrac 2n \int w_{nx}^2 h(x)
-\tfrac 1n \int w_n^2 h''(x)
 + \frac 1{n^2} \int w_{nxx}^2 h(x)
\end{split}\end{equation*}
implies that 
\begin{equation}\label{troisbarres}
\int_{x>0} x w_n^2\leq C \quad \text{and}\quad
\int_{x>0} x(w-w_n)^2\to 0 \quad \text{as $n\to +\infty$.}
\end{equation}
The same holds true in the region $x<0$.

For the functions $w_n$, we have the following identity,
for any $t_1<t_2$:
\begin{equation}\label{purviriel}\begin{split}
&\int x w_n^2(t_2)- \int x w_n^2(t_1)
= -\int_{t_1}^{t_2} \int \Big(2 |D^{\frac 12} w_n|^2+w_n^2 +
 w_n^2 (xQ' -  Q)\Big)dx dt\\ &
 +\int_{t_1}^{t_2} \int \big(
- \tfrac 2n x (1-\tfrac 1n \partial_x^2)^{-1} (2Q'w_{nx} + w_n Q'')_x w_n\big)
+ 2\int_{t_1}^{t_2}  \beta \int x((1-\tfrac 1n \partial_x^2)^{-1} Q')w_n  dx
 dt.
\end{split}\end{equation}
Indeed, 
multiplying the equation of $w_n$ by $Ag(\tfrac xA) w_n$ 
where $g(x)=\arctan(x)$, we find
\begin{equation*}\begin{split}
&\int Ag(\tfrac xA) w_n^2(t_2)- \int Ag(\tfrac xA) w_n^2(t_1)
\\&= -\int_{t_1}^{t_2} \int \Big(2 |D^{\frac 12} w_n|^2 g'(\tfrac xA) 
+ 2 D^{\frac 12} w_n (D^{\frac 12} (w_n g'(\tfrac x A)) - D^{\frac 12}(w_n) g'(\tfrac xA))+ 2 D w_n w_{nx} A g(\tfrac xA)
\\ &+w_n^2 g'(\tfrac xA) 
+ w_n^2 (Ag(\tfrac xA)Q' - g'(\tfrac xA) Q)\Big)dx dt -2
\int_{t_1}^{t_2} \beta(t) \int x((1-\tfrac 1n \partial_x^2)^{-1} Q )  A g(\tfrac xA) w_n\\
&-\frac 2n \int_{t_1}^{t_2} \int ((1-\tfrac 1n \partial_x^2)^{-1} Q )(2 Q' w_{nx} + w_n Q'')_x
A g(\tfrac xA) w_n.
\end{split}
\end{equation*}
Then, \eqref{purviriel} is proved 
using Lemmas \ref{SECONDTERM} and \ref{COMMUTATOR} (see the proof of Lemma \ref{DUAL}
for similar arguments) 
and then passing to the limit as $A\to +\infty$
applying the Lebesgue convergence theorem.

\medskip

From \eqref{purviriel}, we claim that for any $t_1,t_2$,
\begin{equation}\label{fini}
\limsup_{n\to +\infty} \int_{t_1}^{t_2} \|w_n(t)\|_{H^{\frac 12}}^2 dt
<+\infty.
\end{equation}
Proof of \eqref{fini}.
By Claim \ref{TECH} (i), we have
\begin{equation*}\begin{split}
& \tfrac 1n \int x (1-\tfrac 1n \partial_x^2)^{-1} (2Q'w_{nx} + w_n Q'')_x w_n
\\ &= \tfrac 1n \int  w_n  (1-\tfrac 1n \partial_x^2)^{-1} (2xQ'w_{nxx} + 3 x Q''  w_{nx}+  xQ^{(3)}w_n ) 
\\
&- \tfrac 2{n^2} \int w_n (1-\tfrac 1n \partial_x^2)^{-2} (2Q'w_{nx} + w_n Q'')_{xx} 
=\mathbf{I}+ \mathbf{II}.
\end{split}\end{equation*}
As in the proof of Lemma \ref{DUAL} (control of $\mathbf{B}$), 
we have
\begin{equation}\label{doublebarre}
|\mathbf{I}|\leq \frac C{n^{\frac 14}} \|w_n\|_{H^{\frac 12}}\|w_n\|_{L^2},\quad
|\mathbf{II}|\leq \frac C{n^{\frac 12}}  \|w_n\|_{L^2}^2.
\end{equation}
From \eqref{purviriel}, \eqref{controlbeta}, the $L^2$ bounds on $w(t)$
and $w_n(t)$ and \eqref{doublebarre} we obtain
$$
\int_{t_1}^{t_2} \|w_n(t)\|_{H^{\frac 12}}^2 dt
\leq C |t_2-t_1|+\sup_t \int |x| w_n^2(t)+
\frac C{n^{\frac 14}} \int_{t_1}^{t_2} \|w_n\|_{H^{\frac 12}}^2  dt.
$$
For $n$ large enough, we get
$
\int_{t_1}^{t_2} \|w_n(t)\|_{H^{\frac 12}}^2 dt
\leq C.
$
Thus \eqref{fini} is proved. 

By the well-posedness of the equation of 
$w(t)$ in $H^{\frac 12}$, we obtain $\forall t$, $w(t)\in H^{\frac 12}$
and $w_n\to w$ in $H^{\frac 12}$.
Finally, from \eqref{troisbarres} and \eqref{doublebarre},
we obtain \eqref{virielw} by passing to the limit as $n\to \infty$ in
\eqref{purviriel}.

\subsection{Positivity of a quadratic form related to the dual problem}\label{secVIRIEL}
 
In this section, we prove Proposition \ref{QUADRA}. The main ingredient is the following
result.

\begin{proposition}\label{pureQUADRA}
There exists $\lambda_0>0$ such that for all $z\in H^{\frac 12}$,
\begin{equation}\label{proofQUADRA}
\int z (xQ)' =0\quad \Rightarrow\quad 
(\widetilde{\mathcal{L}}z,z)=2 \int |D^{\frac 12} z|^2 + \int z^2 - \int(xQ'+Q) z^2
\geq \lambda_0 \|z\|_{H^{\frac 12}}^2.
\end{equation}
\end{proposition}

\noindent\emph{Proof of Proposition \ref{pureQUADRA}.}
First, we introduce some notation. Recall that 
\begin{equation}\label{defL}
\mathcal L f=-\mathcal{H} f' + f - Q f.
\end{equation}
We define
$S=(xQ)'$. Note that $S=\frac d{dc}{Q_c}_{|c=1}$
and thus by differentiating the equation of $Q_c$ with respect to 
$c$, and taking $c=1$, we find $\mathcal{L} S=-Q.$ 
Observe also that $\mathcal{L}Q=-\mathcal{H}Q'+ Q -Q^2= -\frac 12 Q^2$,
by the equation of $Q$.
Now, we set $T=S-Q$. Then, $\mathcal{L}T=-Q+\frac 12 Q^2=(xQ)'=S$,
by using the explicit expression $Q(x)= \frac 4{1+x^2}$.
We compute $\int TS=\int S^2 - \int QS$. Since $S=\mathcal{H} Q'=\frac 12 Q^2-Q$
(explicit computation), we have $\int S^2=\int (Q')^2$
and $(Q')^2=\frac {64x^2}{(1+x^2)^4}=Q^3-\frac 14 Q^4$,
thus $\int (Q')^2=\int Q^3 -\frac 14 \int Q^4=\int S^2=\int (\frac 12 Q^2-Q)^2
=\frac 14 \int Q^4 -\int Q^3 +\int Q^2$, we find
$\int S^2= \frac 12 \int Q^2$. Moreover, $\int S Q= -\int xQQ'=\frac 12 \int Q^2$,
and so $\int TS=0$. 
Finally, $\int TQ=-\int T \mathcal{L} S=-\int S^2$.

In conclusion, we have proved ($(.,.)$ denotes the $L^2$ scalar product):
\begin{equation}\label{surT}\begin{split}
& S=\tfrac 12 Q^2- Q=(xQ)',~ T=S-Q,~ \mathcal{L}Q=-\tfrac 12 Q^2, ~ \mathcal{L} S=-Q, ~ \mathcal{L} T = S,\\
& (S,Q)= \tfrac 12 \int Q^2, ~ (S,T)=0, ~(T,Q)=-\int S^2.
\end{split}
\end{equation}
Now, we claim the following.

\begin{lemma}\label{TRAV}
There exists $\lambda>0$ such that, for all $\varepsilon>0$,
if $\int w S_\varepsilon=0$, where $S_\varepsilon=S+\varepsilon Q$,
then $(\mathcal{L} w,w)\geq 0$ and $(\widetilde{\mathcal{L}}w,w)\geq \lambda \|w\|_{H^{\frac 12}}^2$.
\end{lemma}

\noindent\emph{Proof of Lemma  \ref{TRAV}.} 
Let $T_\varepsilon=T-\varepsilon S$ and $S_\varepsilon = S+ \varepsilon Q$,
then by \eqref{surT} : $\mathcal{L} T_\varepsilon = S_\varepsilon$ and
$$
(\mathcal{L} T_\varepsilon,T_\varepsilon)=
(S_\varepsilon,T_\varepsilon) = (S,T)+ \varepsilon (-(S,S)+(T,Q)) - \varepsilon^2 (S,Q)
\leq - 2 \varepsilon (S,S)<0.
$$
Moreover, it is clear that if $f_0$, $\lambda_0$ denote respectively the first eigenfunction and first eigenvalue of $\mathcal{L}$ (see Lemma \ref{76})  we have $(S,f_0)=(\mathcal{L} T,f_0)=(T,\mathcal{L} f_0)
=\lambda_0 (xQ',f_0)\neq 0$, since $f_0>0$.
Thus, by Lemma E.1 in \cite{We85}, we obtain the first part of Lemma \ref{TRAV}.

Now, we note that since $xQ'>0$,
\begin{equation}\label{trucdec}\begin{split}
(\widetilde{\mathcal{L}}w,w)& =2 \int |D^{\frac 12} w|^2 + \int w^2 - \int(xQ'+Q) w^2\\
& \geq 
2 \int |D^{\frac 12} w|^2 + \int w^2 - \int Q w^2
= 
\int |D^{\frac 12} w|^2 + (\mathcal{L}w,w).
\end{split}\end{equation}
Using the inequality $\|w\|_{L^4}^2 \leq C \|w\|_{L^2}\|D^{\frac 12} w\|_{L^2}$
(see \eqref{gn})
 and Cauchy-Schwarz' inequality, we have, for some constant $C_0>0$,
$$
\int Q w^2 \leq C \|w^2\|_{L^2}
\leq C_0 \int |D^{\frac 12} w|^2 + \frac 12 \int w^2.
$$
Thus, for $\delta_0>0$ such that $2-C_0\delta_0> 1-\frac {\delta_0}2$, we have
\begin{equation*}\begin{split}
(\widetilde{\mathcal{L}}w,w) & \geq (2-C_0\delta_0) \int |D^{\frac 12} w|^2 + (1-\tfrac 12 \delta_0) \int w^2 - 
(1-\delta_0) \int Q w^2\\
& \geq (1-\delta_0) (\mathcal{L} w,w) + \tfrac {\delta_0}2 \|w\|_{H^{\frac 12}}^2
\geq \tfrac {\delta_0}2 \|w\|_{H^{\frac 12}}^2,
\end{split}\end{equation*}
provided $\int w S_\varepsilon=0$.

\medskip

Now, we finish the proof of Proposition \ref{pureQUADRA}. Let $z\in H^{\frac 12}$ be such that
$\int z S=\int z(xQ)'=0$. Let $w=z+aQ$, where $\int w S_\varepsilon=0$, $0<\varepsilon<\varepsilon_0$,
where $\varepsilon_0$ is to be chosen small enough.
In particular, we have
$$
\int w S_\varepsilon=\int z S_\varepsilon + a \int Q S_\varepsilon=
\varepsilon \int zQ + a \int SQ + a\varepsilon \int Q^2=
\varepsilon \int zQ + a \left(\tfrac 12 + \varepsilon\right) \int Q^2=0,
$$
and so $|a|\leq \frac 2 {\|Q\|_{L^2}} \varepsilon \|z\|_{L^2}$, and $\|w\|_{L^2}\leq 2 \|z\|_{L^2}$
for $\varepsilon_0$ small enough. Similarly, we have $\|z\|_{L^2}\leq 2 \|w\|_{L^2}$, by possibly choosing a smaller $\varepsilon_0$.
By Lemma \ref{TRAV},
we obtain
\begin{equation*}\begin{split}
\frac \lambda 2 \|z\|_{H^{\frac 12}}^2\leq 
\lambda \|w\|_{H^{\frac 12}}^2 \leq (\widetilde{\mathcal{L}} w,w)
=(\widetilde{\mathcal{L}} z,z)+a^2 (\widetilde{\mathcal{L}} Q,Q) + 2 a (\widetilde{\mathcal{L}}Q,z).
\end{split}\end{equation*}
For $\varepsilon_0$ small, we get $(\widetilde{\mathcal{L}} z,z)\geq \frac \lambda 4 \|z\|_{H^{\frac 12}}^2$.

\medskip

Now, we are in a position to prove Proposition \ref{QUADRA}.

\medskip

\noindent\emph{Proof of Proposition \ref{QUADRA}.}
In Proposition \ref{QUADRA}, we want to prove that for $B$ large and $\gamma$ small, and for some $\lambda_1>0$, for all $t$,
$(\widetilde{\mathcal{L}} z(t),z(t))\geq \lambda_1 \|z(t)\|_{H^{\frac 12}}^2$, for
$z(t)=v(t)\sqrt{g'(\frac x {(B+t^2)^{\alpha}})}$ , where
$v=(1-\gamma \partial_x^2)^{-1}(\mathcal{L} w).$
Formally, if $B=+\infty$ and $\gamma=0$, we have $z(t)=v(t)=\mathcal{L} w$ and
$0=\int w Q=- \int w \mathcal{L} S=-\int z S$, and the result follows from Proposition \ref{pureQUADRA}.
Now, we justify that the result persists for large values of $B$ and small values of $\gamma$.

Let $S_{B,\gamma}(t)=(g'(\frac x {(B+t^2)^{\alpha}}))^{-\frac 12} (S-\gamma S'')$.
Then, $\mathcal{L}((1-\gamma \partial_x^2)^{-1}(\sqrt{g'(\frac x {(B+t^2)^{\alpha}})} S_{B,\gamma}(t))=-Q$
and so $\int S_{B,\gamma}(t) z=-\int wQ=0.$
Now, we control $S_{B,\gamma}(t)-S$:
$$
S_{B,\gamma}(t)-S=
\sqrt{1{+}\tfrac {x^2}{(B+t^2)^\alpha}} (S-\gamma S'') -S
=\Big(\sqrt{1{+}\tfrac {x^2}{(B+t^2)^\alpha}}-1\Big) S -\gamma \sqrt{1{+}\tfrac {x^2}{(B+t^2)^\alpha}} S''.
$$
Thus, by elementary estimates and the expression of $S$, we obtain:
$$|S_{B,\gamma}(t,x)-S(x)|\leq \big(B^{-\frac \alpha 2} + \gamma\big) \frac 1 {1+|x|}.$$
It follows that
$$
\left|\int S z(t)\right| = \left| \int (S-S_{B,\gamma}(t)) z(t)\right| 
\leq \big(B^{-\frac \alpha 2} + \gamma\big)\|z\|_{L^2}.
$$
Setting $z=z_1+aQ$, where $\int z_1 S=0$ and $|a|\leq \big(B^{-\frac \alpha 2} + \gamma\big)\|z\|_{L^2}$,
we conclude the proof of Proposition \ref{QUADRA} as at the end of the proof of Proposition \ref{pureQUADRA},
for $B$ large enough and $\gamma$ small enough.

\section{Proof of asymptotic stability - Theorem \ref{TH1}}\label{sec:4}
In this section, we first prove that Theorem \ref{TH2} implies Theorem \ref{TH1}.
Then, we prove that Theorem \ref{LINEARLIOUVILLE} (proved in Section \ref{sec:3})
implies Theorem \ref{TH2}.

\subsection{Proof of Theorem \ref{TH1} assuming Theorem \ref{TH2}}

We follow the strategy of \cite{MMjmpa}, \cite{MM1}, the main idea being to use monotonicity type arguments
(such as Proposition \ref{MONOTONICITY1}) to prove that a limiting solution of \eqref{BO} has uniform
decay in space. See also \cite{MMnonlinearity} for similar use of monotonicity arguments.

\medskip

We consider a solution $u(t)$ of \eqref{BO} in $H^{\frac 12}$ which satisfies
$\|u_0-Q\|_{H^{\frac 12}}=\alpha<\alpha_0$, for $\alpha_0>0$ small enough.
By the stability property, for all $t\in \mathbb{R}$,
$\inf_{y}\|u(t)-Q(.-y)\|_{H^{\frac 12}}\leq C \alpha$.

\medskip

\noindent\emph{1. Decomposition of $u(t)$ around the asymptotic soliton.}
First, we determine the parameter $c^+>0$. It is given by the amount of $L^2$
norm that remains on the region $x>\frac t{10}$ asymptotically as $t\to +\infty$.
Let $\varphi$ be as in \eqref{defphi0}, with $A>1$ so that Proposition \ref{MONOTONICITY1} holds.
Let 
\begin{equation}\label{defcplus}
c^+=\frac 1{\pi \int Q^2} \, \limsup_{t\to +\infty} \int u^2(t,x) \varphi(x-\tfrac t{10}¿) dx.
\end{equation}
From the stability property, $|c^+-1|\leq C \alpha_0$ ($\lim_{+\infty} \varphi=\pi$). Using Lemma \ref{MODULATION}
to decompose $u(t)$ around $Q_{c^+}$, we consider the following decomposition of $u(t)$
\begin{equation}\label{decompcplus}\begin{split}
& u(t,x)=Q_{c^+}(x-\rho(t))+\eta(t,x-\rho(t)),\\
& \int Q_{c^+}' \eta(t,x) dx =0,\quad \sup_t\|\eta(t)\|_{H^{\frac 12}}\leq K \alpha_0.
\end{split}\end{equation}
In what follows, we consider $\alpha_0>0$ small enough, so that the following
holds (by \eqref{ortho}):
\begin{equation}\label{centieme}
\forall t,\quad \frac {99}{100} \leq \rho'(t)\leq \frac {101}{100},\quad
\frac {99}{100} \leq c^+\leq \frac {101}{100}.
\end{equation}

\noindent\emph{2. Monotonicity arguments.} We claim the following estimates:

\begin{lemma}[Asymptotics on $u(t)$]\label{DECAY}
\begin{equation}\label{droite}
\forall y_0>1,\quad
\limsup_{t\to +\infty}\int u^2(t,x) \varphi(x-y_0-\rho(t)) dx\leq \frac {C}{y_0},
\end{equation}
\begin{equation}\label{gauche}
\forall y_0>1,\quad 
\limsup_{t \to +\infty} \int u^2( t  , x ) (\varphi(x-\rho(t)+\tfrac t{10})
-\varphi(x-\rho(t)+y_0))dx \leq  \frac C{y_0},
\end{equation}
\begin{equation}\label{tresgauche}
\lim_{t\to +\infty}
\int u^2(t,x)(\varphi(x-\rho(t)+\tfrac {19}{20} t)-\varphi(x-\rho(t)+\tfrac t{10})) dx
=0,
\end{equation}
\begin{equation}\label{autregauche}
\lim_{t\to +\infty} \int u^2(t,x) \varphi(x-\rho(t)+\tfrac {t}{10} ) dx = c^+ \pi \int Q^2.
\end{equation}
\end{lemma}

\noindent\emph{Proof of Lemma \ref{DECAY}.}
Monotonicity property on the right of the soliton. By \eqref{monotonicity1},
with $\lambda=\frac 12$, we have, for all $y_0>1$,
$$
\int u^2(t,x) \varphi(x-y_0-\rho(t)) dx \leq
\int u^2(0,x) \varphi(x-y_0-\rho(0)-\tfrac 12 t) dx + \frac {C}{y_0}.
$$
Since $\lim_{t\to +\infty} \int u^2(0,x) \varphi(x-y_0-\rho(0)-\tfrac 12 t) dx=0$,
we obtain \eqref{droite}.

\medskip

Monotonicity property on the left of the soliton. By \eqref{monotonicity2}, with
$\lambda=\frac {19}{20}$ and $x_0=\frac {19}{20} t'$, we have
for all $0\leq t'\leq t$,
$$
\int u^2(t,x) \varphi(x-\rho(t)+\tfrac {19}{20} t) dx \leq
\int u^2(t',x) \varphi(x-\rho(t')+\tfrac {19}{20} t') dx + \frac {C} {t'}.
$$
It follows that $\int u^2(t,x) \varphi(x-\rho(t)+\tfrac {19}{20} t) dx$
has a limit as $t\to +\infty$. Set
$$\ell = \lim_{t\to +\infty} \int u^2(t,x) \varphi(x-\rho(t)+\tfrac {19}{20} t) dx,
\qquad 
\ell\geq c^+\pi \int Q^2.$$

Applying \eqref{monotonicity2} with $\lambda=\frac {100}{99}(\frac {19}{20}-\frac 1{1000})<1$ and $x_0=\frac t{1000}$, we find
$$
\int u^2(t,x)\varphi(x-\rho(t)+\tfrac {19}{20} t) dx \leq
\int u^2(\tfrac t{100},x) \varphi(x-\rho(\tfrac t{100}) + \tfrac t{1000}) dx + \frac C{t}.
$$
Since
$$
\limsup_{t\to +\infty} 
\int u^2(\tfrac t{100},x) \varphi(x-\rho(\tfrac t{100}) + \tfrac 1{10} \tfrac t{100}) dx 
\leq c^+ \pi \int Q^2 ,$$
we obtain $c^+ \pi \int Q^2 =\ell$ and  \eqref{tresgauche}.

Fix $y_0>1$, pick $\lambda=\frac 12$. Consider $t_2>t$ and define $t_1=\frac 45 t_2 + 2 y_0$, so
that for $t$ large, $t_1<t_2$. But then, by  \eqref{monotonicity2},
\begin{equation*}\begin{split}
 \int u^2(t_2,x) \varphi(x-\rho(t_2)+ \tfrac {t_2} {10}) dx &= 
\int u^2(t_2,x) \varphi(x-\rho(t_2) +Ê\lambda (t_2-t_2) + y_0) dx \\
&\leq \int u^2(t_1,x) \varphi(x-\rho(t_1)+ y_0) dx +\frac {C} {y_0}.
\end{split}\end{equation*} 
 In light of \eqref{tresgauche} and the existence of $\ell$,  \eqref{gauche} follows.
Thus Lemma \ref{DECAY} is proved.

\medskip

\noindent\emph{3. Construction of a compact limit object.}
Let $t_n\to +\infty$. By the uniform bound on $u(t)$ in $H^{\frac 12}$, there exist
$\widetilde u_0\in H^{\frac 12}$ and 
a subsequence, still denoted by $(t_n)$, such that
$$ u(t_n,.+\rho(t_n))\rightharpoonup \widetilde u_0 \quad
\text{in $H^{\frac 12}$ weak as $n\to +\infty$.}$$
Consider $\widetilde u(t)$ the global $H^{\frac 12}$ solution of \eqref{BO}
such that  $\widetilde u(0)=\widetilde u_0$. 
By \eqref{decompcplus}, $\|\widetilde u_0-Q_{c^+}\|\leq C \alpha_0$
 and thus by the stability property, $\sup_t \inf_y \|\widetilde u(t)-Q(.-y)\|_{H^{\frac 12}}\leq C \alpha_0$.
Let $\widetilde \rho(t)$, $\widetilde \eta(t)$
correspond to the decomposition of $\widetilde u(t)$ around $Q_{c^+}$ given by
Lemma \ref{MODULATION}.

By Theorem \ref{weaku} below and Remark \ref{rk:10}, for all $t\in \mathbb{R}$, we have
\begin{align*}
& u(t_n+t,.+\rho(t_n)) \rightharpoonup \widetilde u(t) 
\quad \text{in $H^{\frac 12}$ weak,}\\
& \rho(t_n+t)-\rho(t_n)\to \widetilde \rho(t)
\quad \text{as $n\to +\infty$}.
\end{align*}
From weak convergence and
Lemma \ref{DECAY}, we claim the following decay estimate on $\widetilde u(t)$:
\begin{equation}\label{decaytilde}
\forall y_0>1,\forall t\in \mathbb{R},\quad
\int_{|x|>y_0} \widetilde u^2(t,x+\widetilde \rho(t)) dx \leq  \frac C{y_0}.
\end{equation}
Indeed, first, from \eqref{droite}, for any fixed $y_0>1$, $t\in \mathbb{R}$, we have
$$
\limsup_{n\to +\infty} \int u^2(t+t_n,x+\rho(t_n))\varphi(x-\rho(t_n+t)+\rho(t_n)-y_0) dx\leq
\frac C{y_0},
$$
and so by weak convergence
$$
\int \widetilde u^2(t,x) \varphi(x -\widetilde \rho(t) -y_0) dx \leq \frac C{y_0}.
$$
Second, from \eqref{gauche}, for fixed $t\in \mathbb{R}$,
$$
\limsup_{n\to +\infty}
\int u^2(t+t_n,x+\rho(t_n))
(\varphi(x-\rho(t_n+t)+\rho(t_n) + \tfrac {t+t_n}{10})
-\varphi(x-\rho(t_n+t)+\rho(t_n)+y_0)) dx \leq \frac C{y_0}.
$$
Note that for fixed $t$, $y_0$, we have
\begin{align*}
\lim_{n\to +\infty} \varphi(x-\rho(t_n+t)+\rho(t_n) + \tfrac {t+t_n}{10})
-\varphi(x-\rho(t_n+t)+\rho(t_n)+y_0)& = \pi-\varphi(x-\widetilde \rho(t)+y_0) \\ &=
\varphi(-x+\widetilde \rho(t) -y_0).
\end{align*}
Thus, we obtain
$$
\int \widetilde u^2(t,x) \varphi(-x+\widetilde \rho(t)-y_0) dx \leq \frac C{y_0}.
$$

Finally, from  \eqref{droite}--\eqref{autregauche}, for any $y_0>1$, we have
$$
\lim_{n\to +\infty} \bigg| \int u^2(t_n,x) (\varphi(x-\rho(t_n)+y_0)
-\varphi(x-\rho(t_n)-y_0)) dx - c^+ \pi \int Q^2 \bigg| \leq \frac C{y_0}.
$$
Thus, by $L^2_{loc}$ convergence, for any $y_0>1$,
$$
\bigg| \int \widetilde u^2_0(x) (\varphi(x+y_0)
-\varphi(x-y_0)) dx - c^+ \pi \int Q^2 \bigg| \leq \frac C{y_0}.
$$
Passing to the limit $y_0\to +\infty$, we obtain
$ \|\widetilde u_0\|_{L^2}=\|\widetilde u(t)\|_{L^2} = \sqrt{c^+} \|Q\|_{L^2}=\|Q_{c^+}\|_{L^2}.$

\medskip

\noindent\emph{4. Conclusion by Theorem \ref{TH2}.}
From Theorem \ref{TH2}, it follows that for some $c_1$ close to $c^+$ and
$x_1$ close to $0$, we have
$$
\widetilde u(t,x)\equiv Q_{c_1}(x-x_1-c_1 t).
$$
But $\|\widetilde u(t)\|_{L^2}  =\|Q_{c^+}\|_{L^2}
$ implies that $c_1=c^+$. Moreover, $\widetilde \rho(0)=0$ and
$\widetilde u(0)=Q_{c^+}(x-x_1)=Q_{c^+}(x)+\widetilde \eta(0,x)$
where $x_1$ is small and $\int \widetilde \eta(0) Q'_{c^+}=0$ imply
$x_1=0$. In conclusion,
$
\widetilde u_0= Q_{c^+}.
$

By a standard argument and \eqref{droite}, \eqref{autregauche}, we have obtained
$$
u(t,.+\rho(t))\rightharpoonup Q_{c^+} 
\quad\text{in $H^{\frac 12}$ weak as $t\to +\infty$,}
$$
\begin{equation}\label{page26}
\lim_{t\to +\infty} 
\int_{x>\frac t{10}} |u(t,x) - Q_{c^+}(x-\rho(t))|^2 dx =0.
\end{equation}
Thus Theorem \ref{TH1} is a consequence of Theorem \ref{TH2}.

\subsection{Proof of Theorem \ref{TH2}}
First, we note that it is sufficient to prove Theorem \ref{TH2} in the case $\int u_0^2=\int Q^2$.
Indeed, for $u_0$ satisfying the assumptions of Theorem \ref{TH2}, 
set $c_1=\int u_0^2/\int Q^2$ and
$\widetilde u(t)=\frac 1{c_1} u(\frac 1{c_1^2} t, \frac 1{c_1} x).$
Then, $|c_1-1|\leq C \alpha_0$ and $\widetilde u$ satisfies \eqref{BO}, $\int \widetilde u^2=\int Q^2$ 
and  $\|\widetilde u_0-Q \|_{H^\frac 12}
\leq C \alpha_0$. Thus, by the stability property -- see Introduction -- for all $t$, there exists $y(t)$ such that 
$\sup_{t} \|\widetilde u(t)- Q(.-y(t))\|_{H^{\frac 12}}\le C'\alpha_0$.
Moreover, $\widetilde u(t)$ also satisfies \eqref{th2:1}.
If we prove $\widetilde u(t,x)=Q(x-t-x_0)$, with $|x_0|\leq C \alpha_0$, the result follows
for $u(t)$.

The proof of Theorem \ref{TH2} is by contradiction. Assume that there exists
a sequence $u_n(t)$ of $H^{\frac 12}$ solutions of \eqref{BO} such that
\begin{align}
&
\sup_{t\in \mathbb{R}} \|u_n(t)-Q(.-\rho_n(t))\|_{L^2}\to 0 \quad
\text{as $n\to +\infty$,}\label{close}\\
& \int u_n^2(0)=\int Q^2,\qquad \eta_n\not \equiv 0,\\
& \forall n,\forall \varepsilon>0,\exists A_{n,\varepsilon}>0, \text{ s.t. }
\forall t\in \mathbb{R},\quad
\int_{|x|>A_{n,\varepsilon}} u_n^2(t,x+\rho_n(t)) dx <\varepsilon,
\label{compacite}
\end{align}
where $\rho_n(t)$ and  $\eta_n(t)$ are defined from $u_n(t)$ by Lemma \ref{MODULATION}.
 Note that    $\int u_n^2(0)=\int Q^2$ implies
\begin{equation}\label{Ltwo}
\forall n,\forall t,\quad 
\int \eta_n^2(t)=-2 \int \eta_n(t) Q.
\end{equation}

Define
\begin{equation}\label{defbn}
0\not \equiv b_n=\sup_t \|\eta_n(t)\|_{L^2}\to 0 \quad \text{as $n\to +\infty$.}
\end{equation}
Then, there exists $t_n$ such that $\|\eta_n(t_n)\|_{L^2}\geq \frac 12 b_n$.
We set
$$
w_n(t,x)=\frac {\eta_n(t_n+t,x)}{b_n}.
$$

For such a sequence $w_n$, we claim the following result.

\begin{proposition}[Weak convergence of the sequence of renormalized solutions]\label{RN}\quad \\
There exists $(w_{n'})$ a subsequence of $(w_n)$ and 
$w\in C(\mathbb{R},L^2(\mathbb{R}))\cap L^{\infty}(\mathbb{R},L^2(\mathbb{R}))$ such that
$$
\forall t\in \mathbb{R},\quad w_{n'}(t)\rightharpoonup w(t) 
\quad \text{in $L^2$ weak as $n\to +\infty$.}
$$
Moreover, $w(t)$ satisfies for some continuous function $\beta(t)$:
\begin{align*}
& w_t=(\mathcal{L} w)_x + \beta(t) Q'\quad \text{on $\mathbb{R}\times \mathbb{R}$},\\
& w(0)\neq 0,\quad  \int wQ=\int wQ'=0,\\
& \forall t\in \mathbb{R},\forall x_0>1,\quad 
\int_{|x|>x_0} w^2(t,x) dx \leq \frac C{x_0}.
\end{align*}
\end{proposition}

Proposition \ref{RN} is in contradiction with Theorem \ref{LINEARLIOUVILLE}. 
Thus, for $\alpha_0>0$ small, for $u(t)$ satisfying  the assumptions of Theorem \ref{TH2}, 
we have  $\eta\equiv 0$ so that
$\rho'(t)=1$ (by Lemma \ref{MODULATION}) and $u(t,x)=Q(x-\rho(0)-t)$, with
$|\rho(0)|\leq C \alpha_0$.

Therefore, we  are reduced to prove Proposition~\ref{RN}.

\medskip

\noindent\emph{Proof of Proposition \ref{RN}.}
One can actually prove a strong $L^2$ convergence result.
See  the end of the proof.

\medskip

Note that the main point in Proposition \ref{RN} is the fact $w\neq 0$.
For this, we need to obtain a strong convergence in $L^2$
for some suitable $t$.

\medskip

\noindent\emph{Decay estimate.}
From Proposition \ref{MONOTONICITY5}, we have 
\begin{equation*}\begin{split}
& \int \eta_n^2(t_0,x) (\varphi(x-x_0)-\varphi(-x_0)) dx
\\ & 
\leq \int \eta_n^2(t,x) (\varphi(x-x_0-\lambda(t_0-t))-\varphi(-x_0-\lambda(t_0-t))) dx
+ \frac {C b_n^2} {x_0}.
\end{split}\end{equation*}
Letting $t\to -\infty$ and using \eqref{compacite}, we obtain, for any
$x_0>1$,
$$
\int \eta_n^2(t_0,x) (\varphi(x-x_0)-\varphi(-x_0)) dx\leq
\frac {C b_n^2} {x_0}.
$$
Similarly, arguing on $\eta_n(-t,-x)$, for any $x_0>1$,
$$
\int \eta_n^2(t_0,x) (\varphi(x_0)-\varphi(x+x_0)) dx\leq
\frac {C b_n^2} {x_0},
$$
which gives, by \eqref{calcul}, similarly as in the proof of \eqref{voiture}:
\begin{equation}\label{UNIFD}
\forall x_0>1,\quad
\int_{|x|>x_0} \eta_n^2(t,x) dx \leq 
\frac {C b_n^2} {x_0}\quad\text{and}\quad  
\int_{|x|>x_0} w_n^2(t,x) dx \leq 
\frac {C} {x_0}.
\end{equation}

\noindent\emph{Local smoothing estimate on $w_n$.} Let
$\varphi$ be defined in \eqref{defphi0} for a fixed value of $A$ ($A=1$ for example).
Then,
\begin{equation}\label{KATOwn}
\int_0^1 \int |D^{\frac 12} (w_n(t,x) \sqrt{\varphi'(x)}) |^2 dx dt\leq C.
\end{equation}
Proof of \eqref{KATOwn}.
First, we claim the following estimate:
\begin{equation}\label{result}
\frac 12 \frac d{dt}\int \eta_n^2 \varphi 
\leq - \frac 12\int |D^{\frac 12} (\eta_n \sqrt{\varphi'})|^2  
+ C \int \eta_n^2
\leq - \frac 12 \int |D^{\frac 12} (\eta_n \sqrt{\varphi'})|^2  + C b_n^2.
\end{equation}
Thus, by integration,
\begin{equation}\label{dix}
\forall t\in \mathbb{R},\quad
\int_{t}^{t+1} \int |D^{\frac 12} (\eta_n \sqrt{\varphi'})|^2  dx dt
\leq C b_n^2\quad \text{and}\quad 
\int_{t}^{t+1} \int |D^{\frac 12} (w_n \sqrt{\varphi'})|^2  dx dt
\leq C.
\end{equation}
Now, we justify \eqref{result}. Using direct computations, Lemma \ref{SECONDTERM},
\eqref{ortho} and then $|\int \eta_n^3 \varphi'|\leq C \int |\eta|^3 \le C \int \eta^2$
(by \eqref{gn}), we get
\begin{align}
\frac 12 \frac d{dt}\int \eta_n^2 \varphi 
&  = - \int (\mathcal{L} \eta_n - \tfrac 12 \eta_n^2)
(\eta_{nx} \varphi + \eta_n \varphi')
  +(\rho_n'-1) \int Q' \eta_n \varphi 
-\frac 12 (\rho_n'-1) \int \eta_n^2 \varphi'\nonumber\\
&= \int \left((\mathcal{H} \eta_{nx}) \eta_{nx} \varphi + (\mathcal{H}\eta_{nx})\eta_n \varphi'\right) - \frac 12 \int \eta_n^2 \varphi' + 
\frac 12 \int \eta_n^2 (-Q'\varphi + Q \varphi')\nonumber
\\ &\quad +\frac 13 \int \eta_n^3 \varphi' +(\rho_n'-1) \int Q' \eta_n \varphi 
-\frac 12 (\rho_n'-1) \int \eta_n^2 \varphi'\nonumber\\
&\leq  \int (\mathcal{H}\eta_{nx})\eta_n \varphi' + C \int \eta_n^2.\nonumber
\end{align}
Using \eqref{COMM1} and then \eqref{gn}, we have
\begin{equation}\begin{split}\label{page28bis}
& -\int (\mathcal{H}\eta_{nx})\eta_n \varphi' =
\int \eta_n D(\eta_n \varphi') =
\int \eta_n  \sqrt{\varphi'} D(\eta_n \sqrt{\varphi'})
+ \int  \eta_n \left(D(\eta_n \varphi') -\sqrt{\varphi'} D(\eta_n \sqrt{\varphi'})  \right)\\ &
\geq \int |D^{\frac 12} (\eta_n \sqrt{\varphi'})|^2 
- C \|\eta_n\|_{L^2} \|D(\eta_n \varphi') -\sqrt{\varphi'} D(\eta_n \sqrt{\varphi'})\|_{L^2}
\\ &\geq \int |D^{\frac 12} (\eta_n \sqrt{\varphi'})|^2 
- C \|\eta_n\|_{L^2} \|\eta_n \sqrt{\varphi'}\|_{L^4}\|D\sqrt{\varphi'}\|_{L^4}
\geq \frac 12 \int |D^{\frac 12} (\eta_n \sqrt{\varphi'})|^2 
- C \|\eta_n\|_{L^2}^2.
\end{split}\end{equation}
(Note that we have used $\|D\sqrt{\varphi'}\|_{L^4}<+\infty$.) 
Thus \eqref{result} is proved.

\medbreak

\noindent\emph{Compactness in $L^2$ for some time.}
From the equation of $\eta_n$ and \eqref{ortho}, it follows that
$$
\frac d{dt} \int \eta_n^2 = - \frac 12 \int Q' \eta_n^2 + (\rho_n'-1) \int Q'\eta_n
\quad \text{and so}\quad
\left|\frac d{dt} \int \eta_n^2\right|\leq C_0 \int \eta_n^2.
$$
In particular, by the definition of $t_n$, 
$\forall t\in [0,1],$ $\int \eta_n^2(t+t_n) \geq
e^{-C_0} b_n^2$ and so
\begin{equation}\label{nonvanish}
\forall t\in [0,1],\quad \|w_n(t)\|_{L^2}\geq e^{-\frac 12 C_0}=\delta >0.
\end{equation}

It follows from \eqref{KATOwn} that for all $n$, there exists $\tau_n\in [0,1]$
such that $\int |D^{\frac 12} (w_n(\tau_n) \sqrt{\varphi'}) |^2 \leq C$.
Thus, there exists a subsequence of $(w_n)$ (still denoted by $(w_n)$) and
$s_0\in [0,1]$, $W\in H^{\frac 12}$ such that
$$
w_n(\tau_n) \sqrt{\varphi'}\rightharpoonup W \quad \text{in $H^{\frac 12}$ weak,}
\qquad \tau_n\to s_0\quad \text{as $n\to +\infty$.}
$$
But (by possibly extracting a further subsequence), there exists $w_{s_0}\in L^2$
such that
$$
\tau_n\to s_0,\quad 
w_n(\tau_n)  \rightharpoonup w_{s_0} \quad \text{in $L^2$ weak
as $n\to +\infty$.}
$$
It follows that $W= w_{s_0} \sqrt{\varphi'}$.
Since $\sqrt{\varphi'}>0$ on $\mathbb{R}$, we get
$$
w_n(\tau_n)  \to w_{s_0} \quad \text{in $L^2_{loc}$ 
as $n\to +\infty$.}
$$
By \eqref{UNIFD} and  \eqref{nonvanish}, we finally get
\begin{equation}\label{strongonetime}
w_n(\tau_n)\to w_{s_0} \quad \text{in $L^2$ 
as $n\to +\infty$},\quad \int w_{s_0} Q'=0,\quad 
w_{s_0}\neq 0.
\end{equation}
Note also that from \eqref{Ltwo} and $\int \eta_n Q'=0$, we have
\begin{equation}\label{strongortho}
\int w_{s_0} Q=\int w_{s_0} Q'=0.
\end{equation}

\noindent\emph{Weak convergence for all time.}
Consider $\widetilde w(t)\in C(\mathbb{R},L^2(\mathbb{R}))$ the unique
solution of
$$
\widetilde w_t =(\mathcal{L} \widetilde w)_x\quad \text{on $\mathbb{R}\times \mathbb{R}$},\quad
\widetilde w(s_0)=w_{s_0}, \quad \text{on $\mathbb{R}$.}
$$
(It is clear by a standard energy estimate and regularization arguments that the corresponding Cauchy 
problem is well-posed in $L^2$).

\medskip

Now, to obtain weak convergence, we need to remove  some terms from the equation of $w_n$, following some arguments in \cite{MMjmpa}, Lemma 8 and beginning of proof of Lemma 11.
We write
\begin{align*}
w_{nt} & = (\mathcal{L} w_n - \tfrac {b_n} 2 w_n^2)_x +\frac 1 b_n (\rho_n'-1) 
(Q+b_n w_n)_x \\
& =  
(\mathcal{L} w_n)_x  - \tfrac {b_n} 2(w_n^2)_x + \beta_n Q' + b_n F_n' 
+ b_n \widetilde \beta_n (w_n)_x,
\end{align*}
where
$$
\beta_n=\frac 1{\int (Q')^2}\int w_n \mathcal{L}(Q''),\quad \widetilde \beta_n =Ê\frac 1{b_n} (\rho_n'-1),
\quad F_n=\frac 1{b_n} (\widetilde \beta_n - \beta_n) Q.
$$

Set $\widetilde w_n(t)=w_n(t)-Q' \int_{\tau_n}^t \beta_n(s) ds.$ Then, the equation 
of $\widetilde w_n(t)$ writes
$$
\widetilde w_{nt} = (\mathcal{L} \widetilde w_n)_x -\tfrac {b_n} 2 (w_n^2)_x 
+ b_n F_n' + b_n \widetilde \beta_n (\widetilde w_n)_x + b_n \widetilde \beta_n Q'' \int_{\tau_n}^t
\beta_n(s) ds.
$$
We claim the following weak convergence result.

\begin{lemma}\label{WEAKwn}
For all $t\in \mathbb{R}$,
$$
\widetilde w_n(t)\rightharpoonup \widetilde w(t) \quad
\text{in $L^2$ weak.}
$$
\end{lemma}

Assuming this lemma, from \eqref{nouveau}, we have, for all $t$,
$$
\widetilde \beta_n(t) \to \widetilde \beta(t)=\frac 1{\int (Q')^2} \int \widetilde w \mathcal{L}(Q''),\quad
\int_{\tau_n}^t \widetilde \beta_n(s)ds \to \int_{s_0}^t \widetilde \beta(s) ds.
$$
Set $w(t)=\widetilde w(t) + Q' \int_{s_0}^t \widetilde \beta(s) ds$. Then,
$w(t)$ solves 
$$w_t = (\mathcal L w)_x + \widetilde \beta Q',$$
and $w(s_0)=w_{s_0}\neq 0$.
Moreover, for all $t\in \mathbb{R}$,
$$
w_n(t) \rightharpoonup w(t) \quad \text{in $L^2$ weak.}
$$
Finally, from \eqref{Ltwo} and $\int \eta_n Q'=0$, we have
$\int w(t) Q=\int w(t)Q'=0,
$
and by weak convergence and \eqref{UNIFD}, we have
$$
\forall x_0>1,\forall t, \quad
\int_{|x|>x_0} w^2(t,x)dx \leq \frac C{x_0}.
$$
Thus, we are reduced to prove Lemma \ref{WEAKwn}.

\medskip

\noindent\emph{Proof of Lemma \ref{WEAKwn}.}
Set
$$ G_{1,n}= -\frac 12 w_n^2,\quad 
G_{2,n}=
F_n + \widetilde \beta_n \widetilde w_n +   \widetilde \beta_n Q' \int_{\tau_n}^t
\beta_n(s) ds,\quad
G_n=G_{1,n}+G_{2,n}.
$$
Observe that
$$
\|G_{1,n}\|_{L^1} + \|G_{2,n}\|_{L^2} \leq C(t),\quad
\text{with $C(t)$ bounded on bounded intervals.}
$$
Let $T\in \mathbb{R}$. By $\sup_{t} \|w_n(t)\|_{L^2}\leq C$ and the expression of
$\widetilde w_n$, we have $\sup_{[-T,T]} \|\widetilde w_n(t)\|_{L^2}\leq C_T$.

Let $g\in C_0^\infty(\mathbb{R})$ and let $v$ solve the problem
\begin{equation*}
\left\{\begin{aligned}
&  \partial_t v = \mathcal{L} (v_x)\\
& v_{|t=T} = g.
\end{aligned}\right.\end{equation*}
Then
\begin{align*}
& \int (\widetilde w_n-\widetilde w)(T) g(x) dx -
\int (w_n(\tau_n)-w(\tau_n)) v(\tau_n) dx
=\int_{\tau_n}^T \int \partial_t ((\widetilde w_n-\widetilde w)(t) v(t,x)) dx dt \\
& \int_{\tau_n}^T \int ((\mathcal{L} \widetilde w_n)_x - (\mathcal{L} \widetilde w)_x
+ b_n (G_n)_x) v(t,x) +
(\widetilde w_n-\widetilde w) (\mathcal{L} v)_x dx dt\\
&=-b_n \int_{\tau_n}^T \int G_n v_x(t,x) dx dt.
\end{align*}
The energy method gives
$$
\|v\|_{L^\infty([\tau_n,T],L^2(\mathbb{R}))}+
\|v_x\|_{L^\infty([\tau_n,T]\times \mathbb{R})}+
\|v_x\|_{L^\infty([\tau_n,T],L^2(\mathbb{R}))}\leq C.
$$
Moreover, by continuity of $t\mapsto w(t)$ in $L^2$,
\begin{equation*}
\begin{split}
&\lim_{n\to +\infty} \int (w_n(\tau_n)-w(\tau_n)) v(\tau_n) dx
\\ &= \lim_{n\to +\infty} \int (w_n(\tau_n)-w(s_0)) v(\tau_n)dx +
\lim_{n\to +\infty} \int (w(s_0)-w(\tau_n)) v(\tau_n)dx=0.
\end{split}
\end{equation*}
Thus,
$$
\widetilde w_n(T)\rightharpoonup \widetilde w(T)
\quad \text{as $n\to +\infty$.}
$$
and the proof of Lemma \ref{WEAKwn} is concluded.

\medskip

\noindent\emph{Alternate proof by strong $L^2$ convergence for all time.}
Now, we use   Theorem \ref{THA} in Section \ref{sec:6} to prove 
strong $L^2$ convergence of the sequence $(w_n(t))$ for all $t$.

Let $T>0$.
Set 
\begin{equation}\label{defzeta}
\zeta_n(t,x)=w_n(t,x-\rho_n(t)+\rho_n(0))
+\frac 1{b_n} [Q(x-(\rho_n(t)-\rho_n(0)))-Q(x-t)],
\end{equation}
 so that
$$
u_n(t,x+\rho_n(0))
= Q(x-\rho_n(t)+\rho_n(0)) + b_n w_n(t,x-\rho_n(t)+\rho_n(0))
= Q(x-t) + b_n \zeta_n(t,x),
$$
and $\zeta_n$ satisfies 
\begin{align*}
& (\zeta_n)_t = (-\mathcal{H}(\zeta_n)_x - Q(x{-}t) \zeta_n)_x - \frac {b_n}2 (\zeta_n^2)_x,\\
& \|\zeta_n(t)\|_{L^2}\leq C_T, \quad \forall  t\in [-T,T]. 
\end{align*}
Indeed, since $|\rho'_n(t)-1|\leq C \|\eta_n\|_{L^2} \leq C b_n$,
we have  
\begin{equation}\label{surrhon}
|\rho_n(t)-\rho_n(0)-t|\leq C b_n |t|,
\end{equation}
and the estimate on $\zeta_n$ follows.

On the one hand, Theorem \ref{THA} applied to $\frac 1{C_T} \zeta_n$ for $n$ large enough
(so that $b_n$ is small enough) implies that $t\in [-T,T] \mapsto \zeta_n(t)\in L^2$ is equicontinuous
in $n$.

On the other hand, from \eqref{KATOwn}, we have
$$
\int_{[-T,T]} \int \left|D^{\frac 12} ( \zeta_n (t,x) \sqrt{\varphi'(x)})\right|^2 dx dt \leq C_T,
$$
and the decay property \eqref{UNIFD} also holds for $\zeta(t)$ on $[-T,T]$ with constant
depending on $T$.
In particular, there exists $N\subset [-T,T]$ of zero Lebesgue measure such that
for all $t\in [-T,T]\setminus N$,
$ \int |D^{\frac 12} ( \zeta_n (t,x) \sqrt{\varphi'(x)})|^2 dx dt <+\infty.
$
Now, we choose a dense and countable subset $I$ of $[-T,T]$ such that
for all  $t\in I$,
$ \int |D^{\frac 12} ( \zeta_n (t,x) \sqrt{\varphi'(x)})|^2 dx dt <+\infty.
$
Arguing as in the proof of \eqref{strongonetime}, and using a diagonal argument,
there exists a subsequence of $(\zeta_n)$ which we will still denote by
$(\zeta_n)$ such that for any 
$t\in I$,
$\zeta_n(t)\to \zeta(t)$  in $L^2$ strong as $n\to +\infty$.
Using the equicontinuity, we obtain
\begin{equation}\label{strongL2zeta}
\forall t\in [-T,T],\quad 
\zeta_n(t)\to \zeta(t) \quad \text{in $L^2$ strong as $n\to +\infty$.}
\end{equation}
By \eqref{surrhon} and $|\rho'_n-1|\leq C b_n$, 
we may also assume that for the same subsequence
\begin{equation}\label{convrhon}
\forall t\in [-T,T],\quad 
\frac 1{b_n} (\rho_n(t)-\rho_n(0)-t) \to \kappa(t).
\end{equation}

Now, we deduce from \eqref{defzeta}, \eqref{strongL2zeta} and \eqref{convrhon} that
$$
\forall t\in [-T,T],\quad w_n(t)\to w(t)=\eta(t,.+t) + \kappa(t) Q'\quad
\hbox{in $L^2$ strong as $n\to +\infty$.}
$$

\subsection{Proof of Remark \ref{re:1}}\label{sec:4.3}
Let $u(t)$ be a solution satisfying the assumptions of Theorem \ref{TH1}. Let $c^+$,
$\rho(t)$ and $\eta(t)$ be defined as in the proof of Theorem \ref{TH1}. In particular,
by \eqref{page26}, we have
\begin{equation}
\lim_{t\to +\infty} \int_{x>\frac t{10} -\rho(t)} |\eta(t,x)|^2 dx=0.
\end{equation}
To prove \eqref{inre:1}, we use the identity \eqref{page10bis} on $\eta$, where
$\varphi=\frac \pi 2 + \arctan(\frac x{A})$, $A>1$ large enough be to defined later:
\begin{align*}
\frac 12 \frac d{dt}\int \eta^2 \varphi 
 &= \int (\mathcal{H}\eta_{x})\eta \varphi'+ \int  (\mathcal{H} \eta_{x}) \eta_{x} \varphi   - \frac 12 \int \eta^2 \varphi' + 
\frac 12 \int \eta^2 (-Q'\varphi + Q \varphi')
\\ &\quad +\frac 13 \int \eta^3 \varphi' +(\rho'-1) \int Q' \eta \varphi 
-\frac 12 (\rho'-1) \int \eta^2 \varphi'.
\end{align*}
We claim that for $A$ large enough and $\alpha_0$ small enough, for  $C>0$
independent of $A$,
\begin{equation}\label{claimre1}
\frac 12 \frac d{dt}\int \eta^2 \varphi 
 \leq   \int (\mathcal{H}\eta_{x})\eta \varphi'  + C \int \frac {\eta^2}{1+x^2}.
\end{equation}
Indeed, by Lemma \ref{SECONDTERM}, we have $\int (\mathcal{H} \eta_{x}) \eta_{x} \varphi
\leq \frac C A \int \eta^2 \varphi'$. By the definition of $Q$,
$\int \eta^2 (-Q'\varphi + Q \varphi')\leq  C \int \frac {\eta^2}{1+x^2}$.
By \eqref{appendix1} (note that the constant in \eqref{appendix1} is independent of $A$)
$|\int \eta^3 \varphi'|\leq C \alpha_0 \int \eta^2 \varphi'$. 
Finally, the last two terms are controlled using \eqref{ortho}, so that \eqref{claimre1} is proved
for $A$ large enough, $\alpha_0$ small enough.

Now, we use \eqref{page28bis} on $\eta$.
We obtain
\begin{equation}
\frac 12 \frac d{dt}\int \eta^2 \varphi \leq 
- \int |D^{\frac 12} (\eta \sqrt{\varphi'})|^2 
+ C \|\eta\|_{L^2} \|\eta \sqrt{\varphi'}\|_{L^4}\|D\sqrt{\varphi'}\|_{L^4}
+ C \int \frac {\eta^2}{1+x^2}.
\end{equation}
Note that for $A>1$, we have $\frac 1{1+x^2}\leq C \varphi$ on $\mathbb{R}$.
Let $t_0>0$. Integrating the above estimate on $[t_0,t_0+1]$, we get
$$
\int_{t_0}^{t_0+1} \int |D^{\frac 12} (\eta \sqrt{\varphi'})|^2dt
\leq C \sup_{t\in [t_0,t_0+1]} \left(\int \eta^2(t)\varphi + C  \|\eta \sqrt{\varphi'}\|_{L^4}\|D\sqrt{\varphi'}\|_{L^4}\right).
$$

On the other hand, by \eqref{COMM1}, we have
\begin{equation*}\begin{split}
\int |D^{\frac 12} \eta|^2 \varphi'&\leq 2 \int |D^{\frac 12} (\eta \sqrt{\varphi'})|^2 +
2 \int |(D^{\frac 12}\eta) \sqrt{\varphi'} - D^{\frac 12}(\eta \sqrt{\varphi'})|^2 \\
&\leq  2 \int |D^{\frac 12} (\eta \sqrt{\varphi'})|^2  + C \|\eta\|_{L^4}^2 \|D^{\frac 12} \sqrt{\varphi'}\|_{L^4}^2\leq 2 \int |D^{\frac 12} (\eta \sqrt{\varphi'})|^2  + C \|D^{\frac 12} \sqrt{\varphi'}\|_{L^4}^2.
\end{split}
\end{equation*}
Thus, we obtain
\begin{equation*}
\int_{t_0}^{t_0+1} \int |D^{\frac 12} \eta|^2 \varphi'dt
\leq C\sup_{t\in [t_0,t_0+1]} \left(\int \eta^2(t)\varphi + C  \|\eta \sqrt{\varphi'}\|_{L^4}\|D\sqrt{\varphi'}\|_{L^4}\right)+ C  \|D^{\frac 12} \sqrt{\varphi'}\|_{L^4}^2.
\end{equation*}
We have $\|D^{\frac 12} \sqrt{\varphi'}\|_{L^4}^2\leq CA^{-\frac 32}$,
$\|D\sqrt{\varphi'}\|_{L^4}\leq C A^{-\frac 54}$ and
$ \|\eta \sqrt{\varphi'}\|_{L^4} \leq \|\eta\|_{L^8}\|\sqrt{\varphi'}\|_{L^8}\leq C A^{-\frac 38}$.
Therefore,
\begin{equation*}
\int_{t_0}^{t_0+1} \int |D^{\frac 12} \eta(t,x)|^2 \frac {dxdt}{1+(\frac xA)^2}
\leq A \sup_{t\in [t_0,t_0+1]} \left(\int \eta^2(t)\varphi\right) + C A^{-\frac 12}.
\end{equation*}
We now choose $A$ depending on $t_0$:
$$
A=A_{t_0}=\min\left(\frac {\sqrt{t_0}}2,\left(\sup_{t\in [t_0,t_0+1]} \int_{x\geq \frac t{10} - \rho(t)} \eta^2(t,x)dx\right)^{-\frac 12}\right).
$$
For this choice of $A_{t_0}$, we have $\lim_{t_0\to +\infty}A_{t_0}=+\infty$ and,
since $\frac t{10}-\rho(t)\leq - \frac t2$,  
$$
A \sup_{t\in [t_0,t_0+1]} \left(\int \eta^2(t)\varphi\right) 
\leq C A \sup_{t\in [t_0,t_0+1]} \left( \int_{x\geq \frac t{10} - \rho(t)} \eta^2(t) \right) +  \frac {CA} { t_0}
\leq C A^{-1}.
$$
so that $\lim_{t_0\to +\infty}  A \sup_{t\in [t_0,t_0+1]} \left(\int \eta^2(t)\varphi\right) =0$.
It follows that
\begin{equation*}
\lim_{t_0\to +\infty}
\int_{t_0}^{t_0+1} \int |D^{\frac 12} \eta(t,x)|^2 \frac {dxdt}{1+x^2}=0.
\end{equation*}

\section{Multi-soliton case}\label{sec:5}

Using the previous arguments and the strategy of \cite{MMT} for the gKdV equation,
we obtain the following result concerning multi-soliton solutions of \eqref{BO}.

\begin{theorem}[Asymptotic stability of a sum of decoupled solitons]\label{MULTI}\quad\\
Let $N\geq 1$ and $0<c_1^0<\ldots<c_N^0$. There exist $L_0>0$, $A_0>0$
and $\alpha_0>0$ such that if $u_0\in H^{\frac 12}$ satisfies
for some $0\leq \alpha< \alpha_0$, $L\geq L_0$,
\begin{equation}\label{MULTI:1}
\bigg\|u_0- \sum_{j=1}^N Q_{c_j^0} (.-y_j^0)\bigg\|_{H^{\frac 12}}
\leq \alpha  \quad \text{where} \quad
\forall j\in \{2,\ldots, N\},\quad 
 y_j^0-y_{j-1}^0\geq L,
\end{equation}
and if  $u(t)$ is the solution of \eqref{BO} corresponding to $u(0)=u_0$,
then there exist $\rho_1(t), \ldots, \rho_N(t)$ such that
the following hold
\begin{description}
\item{{\rm (a)}} Stability of the sum of $N$ decoupled solitons.
\begin{equation}\label{}
	\forall t\geq 0,\quad
	\bigg\| u(t)- \sum_{j=1}^N Q_{c^0_j}(x-\rho_j(t)) \bigg\|_{H^{\frac 12}}
	\leq A_0 \left(\alpha+ \frac 1 L\right).
\end{equation}
\item{\rm (b)} Asymptotic stability of the sum of $N$ solitons. There exist
$c_1^+, \ldots,c_N^+$, with $|c_j^+-c_j^0|\leq  A_0 \left(\alpha+ \frac 1 {L}\right)$,
such that
\begin{equation}\label{}
\forall j,\quad 
u(t,.+\rho_j(t))   \rightharpoonup  Q_{c^+_j} \quad 
\text{in $H^{\frac 12}$ weak as $t\to +\infty$},
\end{equation}
\begin{equation}\label{}
	\bigg\| u(t)- \sum_{j=1}^N Q_{c^+_j}(.-\rho_j(t)) \bigg\|_{L^2(x\geq \frac 1{10} {c_1^0} t)}\to 0,\quad
	\rho_j'(t)\to c_j^+ \quad \text{as $t\to +\infty$}.
\end{equation}
\end{description}
\end{theorem}

Recall that the Benjamin-Ono equation admits explicit multi-soliton solution. We denote by
  $U_N(x;{c_j},{y_j})$ the explicit family of $N$-soliton profiles,
  see e.g. \cite{MAT} formula (1.7) and Appendix A
(see also references in \cite{MAT}).
We obtain the following corollary of the above Theorem and
the continuous dependence of the solution in $H^{\frac 12}$. 

Let $N\geq 1$,  $0<c_1^0<\ldots<c_N^0$ and set
$$d_N(u)=\inf\big\{ \|u-U_N(.;c_j^0,y_j)\|_{H^{\frac 12}},\ y_j\in \mathbb{R}\big\}.$$

\begin{corollary}[Asymptotic stability in $H^{\frac 12}$ of multi-solitons]\label{cor:1}\quad \\
For all $\delta>0$, there exists $\alpha>0$ such that
if $d_N(u_0)\leq \alpha$ then for all $t\in \mathbb{R},$ $d_N(u(t))\leq \delta$.
\end{corollary}

Recall that a result of stability in $H^1$ of double solitons  for the BO equation
was proved by variational methods in \cite{NL}. See also \cite{MAT} for stability related results.

\subsection{Sketch of the stability argument \cite{W}}\label{sec:51}

For the reader's convenience, we now sketch the proof of the stability argument for one soliton
(see  statement in the Introduction).
Let $u(t)$ be an $H^{\frac 12}$ solution of \eqref{BO} such that $u(0)$ is close to $Q$ in $H^{\frac 12}$.
Let $c^+>0$ be close to $1$ such that $\int u^2(0)=c^+ \int Q^2$. We use Lemma \ref{MODULATION}
on $u(t)$ around $Q_{c^+}$ so that 
$\eta(t,x)=u(t,x+\rho(t))-Q_{c^+}(x)$ satisfies $\int \eta(t)Q_{c^+}'=0$ and by $L^2$ 
conservation $\int \eta(t) Q_{c^+}= -\frac12 \int \eta^2(t)$.

We define the functional 
\begin{equation}\label{defG}
\mathcal{G}(u(t))=E(u(t)) + {c^+}  \int u^2(t).
\end{equation}
Observing that  $\mathcal{G}(u(t))=\mathcal{G}(u(0))$ and so expanding $u(t)$ in $\mathcal{G}(u(t))$, we obtain
$$
(\mathcal{L}_{c^+} \eta(t),\eta(t))+O(\eta^3(t))=(\mathcal{L}_{c^+} \eta(0),\eta(0))+O(\eta^3(0))
$$
where $\mathcal{L}_{c^+}\eta=-\mathcal{H} \eta_x + c^+ \eta - Q_{c^+} \eta.$
By the positivity property of $\mathcal{L}_{c^+}$, (property \eqref{weinstein}
of $\mathcal{L}$ and a scaling argument), we then obtain
$$
\|\eta(t)\|_{H^{\frac 12}} \leq C \|\eta(0)\|_{H^{\frac 12}}.
$$
Note that $\left|\int \eta Q_{c^+}\right|\leq C \|\eta\|_{L^2}^2$ replaces the orthogonality condition
$\int \eta Q_{c^+}=0$.

\subsection{Sketch of the proof of Theorem \ref{MULTI}}\label{sec:52}
The proof is the same as the proof of Theorem 1 in \cite{MMT}.

First, we recall four lemmas (corresponding to Lemmas 1--4 in \cite{MMT})
which are the main tools in proving Theorem \ref{MULTI}.

\begin{lemma}
 [Decomposition of the solution]\label{lem:1}
 There exist $L_1, \alpha_1, K_1>0$ such that the
following is true. If for $L>L_1$, $0<\alpha<\alpha_1$,
$t_0>0$, 
\begin{equation*}
\sup_{0\le t\le t_0} \Big( \inf_{y_j>y_{j-1}+L}\Big\{
\Big\|u(t,.)-\sum_{j=1}^N Q_{c^0_j}(.-y_j)\Big\|_{H^{\frac 12}}
\Big\}\Big)<\alpha,
\end{equation*}
then there exist unique $C^1$ functions
$c_j:[0,t_0]\to (0,+\infty),$ $\rho_j:[0,t_0]\to \mathbb{R}$, such that
\begin{equation*}
\eta(t,x)=u(t,x)-R(t,x)\quad  \hbox{where} \quad R(t,x)=\sum_{j=1}^N R_j(t,x),
\quad R_{j}(t,x)=Q_{c_j(t)}(x-\rho_j(t)),
\end{equation*}
satisfies the following orthogonality conditions
\begin{equation*}
\forall j, \forall t\in [0,t_0],\quad \int R_j(t) \eta(t)=
\int (R_j(t))_x \eta(t)=0.
\end{equation*}
Moreover,  there exists $C>0$ such that 
$\forall t\in [0,t_0],$
\begin{equation*}
\|\eta(t)\|_{H^{\frac 12}}+
\sum_{j=1}^N |c_j(t)-c^0_j|
\le C \alpha,
\quad 
\forall j,~
\left|c_j'(t)\right|
+\left|\rho_j'(t)-c_j(t)\right|
\le C\left(\|\eta(t)\|_{L^2}+\frac 1{L}\right).
\end{equation*}
\end{lemma}

\begin{remark}\label{rk:modulation}
In the rest of the argument, the modulation in the scaling
parameter  for all time 
(i.e. the introduction of $c_j(t)$) is not  necessary. Indeed, modulation at
$t=0$ would be sufficient since we deal with the subcritical case. However, we have
preferred to introduce this modulation to match the strategy of \cite{MMT}.
\end{remark}

Expanding $u(t)$ in  the energy conservation and using $E(Q_c)=c^2 E(Q)$, we have

\begin{lemma} \label{lem:2}
There exists $C>0$ such that in the context of Lemma \ref{lem:1},
 $\forall t\in [0,t_0]$,
\begin{equation*}
\bigg| 
E(Q) \sum_{j=1}^N \left[c_j^2(t)) -c_j^2(0)\right] 
 +\frac 12 \int \, (\eta_x \mathcal{H}\eta - R\eta^2)(t)  \bigg |
  \leq
C \left( \|\eta(0)\|_{H^{\frac 12}}^2 
+ \|\eta(t)\|_{H^{\frac 12}}^3
+ \frac 1{L}\right) .
\end{equation*}
\end{lemma}

We consider $\varphi$ defined as in \eqref{defphi0}, with $A$ large enough, and we set
$$
\forall j\in \{2,\ldots,N\},\quad
\mathcal{I}_j(t)=\int u^2(t,x)\varphi(x-m_j(t)) dx ,\quad
m_j(t)=\frac 12 (\rho_{j-1}(t)+\rho_j(t)).
$$
Then, proceeding as in the proof of Proposition \ref{MONOTONICITY1}, we obtain the following.

\begin{lemma}\label{lem:3}
There exists $C>0$ such that 
in the context of Lemma \ref{lem:1},
$$
\forall j \in \{2,\ldots,N\},\ \forall t\in [0,t_0], \quad
\mathcal{I}_j(t)-\mathcal{I}_j(0)\leq \frac C{L}.
$$
\end{lemma}

Finally, setting
$
c(t,x)= c_1(t)+\sum_{j=2}^N (c_j(t)-c_{j-1}(t)) \varphi(x-m_j(t)),
$
and proceeding as in the proof of Propositions \ref{QUADRA} and \ref{pureQUADRA},
we have

\begin{lemma}\label{lem:4}
There exists $\lambda>0$ such that 
in the context of Lemma \ref{lem:1}, 
$$
\forall t\in [0,t_0],\quad 
\mathcal{G}_N(t):=\int \eta_x \mathcal{H}\eta + c(t,x) \eta^2 - Q \eta^2 \geq \lambda \|\eta(t)\|_{H^\frac 12}^2.
$$
\end{lemma}
Recall that the introduction of the functional $\mathcal{G}_N(t)$ for the problem of stability
of multi-soliton solutions is justified as follows. 
For the stability of one soliton, the suitable functional is $\mathcal{G}(u(t))$ defined in \eqref{defG}.
For the case of $N$ solitons, we introduce the functional $\mathcal{G}_N(t)$ which  is approximately  $E(u(t))+c_j(0) \int  u^2(t)$ around the soliton $Q_{c_j}$.
Then, we observe (using the energy conservation and  Lemma \ref{lem:2})
that this quantity is almost decreasing. This is sufficient to conclude the stability 
argument for several solitons. We now sketch the argument. We refer to \cite{MMT}, Section 3
for more details in the stability proof.

\medskip

\noindent\emph{Sketch of the proof of the stability.}
Let
$$
\mathcal{V}_{A_0}(L,\alpha)=
\bigg\{ u\in H^{\frac 12}; \inf_{y_j-y_{j-1}\geq L}
\bigg\|u-\sum_{j=1}^N Q_{c_j^0}(.-y_j)\bigg\|_{H^{\frac 12}}\leq
A_0\left(\alpha+\frac 1{L}\right)\bigg\}.
$$
Part (a) of Theorem \ref{MULTI}
is a consequence of the following proposition and   continuity arguments.
\begin{proposition}\label{prop:6}
There exist $A_0>0$, $L_0>0$ and $\alpha_0>0$ such that, for all $u_0\in H^{\frac 12}$,
if
$$
\bigg\|u_0-\sum_{j=1}^N Q_{c_j^0}(.-y_j^0)\bigg\|_{H^{\frac 12}}\leq \alpha,
$$
where $L\geq L_0$, $0<\alpha<\alpha_0$, $y_j^0-y_{j-1}^0+L$, and if for $t^*>0$,
$$
\forall t\in [0,T^*],\quad 
u(t)\in \mathcal{V}_{A_0}(L,\alpha),
$$
where $u(t)$ is the solution of \eqref{BO}, then
 $$
\forall t\in [0,T^*],\quad 
u(t)\in \mathcal{V}_{\frac 12 A_0}(L,\alpha).
$$
\end{proposition}

The proof of Proposition \ref{prop:6} is exactly the same as the proof of Proposition 1 in \cite{MMT}, using Lemmas \ref{lem:1}--\ref{lem:4}. 
In particular, we first prove
\begin{equation}\label{eq:lem5}
\forall t\in [0,t^*],\quad
\sum_{j=1}^N |c_j(t)-c_j(0)|\leq C_1 \left(\|\eta(t)\|_{H^{\frac 12}}^2 +\|\eta(0)\|_{H^{\frac 12}}^2+\frac 1{L}\right),
\end{equation}
and then
\begin{equation}\label{eq:lem6}
\|\eta(t)\|_{H^{\frac 12}}^2\leq C_2 \left(\|\eta(0)\|_{H^{\frac 12}}^2+\frac 1{L}\right),
\end{equation}
where $C_1,$ $C_2>0$ are independent of $A_0$,
and we then conclude by using the decomposition of $u(t)$ is terms of $\eta(t)$ and
$R(t)$.

Note that in proving \eqref{eq:lem5}, we make use of the following algebraic fact:
$$E(Q_c)=c^2 E(Q),\quad \int Q_c^2 = c \int Q^2,\quad E(Q)=-\frac 12 \int Q^2.$$
The last formula is easily obtained from the equation of $Q$ multiplying by $Q$
and then by $xQ'$ and using $\int (\mathcal{H}Q') (xQ')=0$.
This allows  us to prove the following estimate
\begin{equation*}
\bigg|E(Q) \sum_{j=1}^N \left(c_j(t)-c_j(0)\right) 
+ \int Q^2 \sum_{j=1}^N \left\{c_j(0) \left(c_j(t)-c_j(0)\right)\right\} \bigg|\leq
C \sum_{j=1}^N \left|c_j(t)-c_j(0)\right|^2.
\end{equation*}
which is the analogue of (44) in \cite{MMT}.

\medskip

The proof of part (b) of Theorem \ref{MULTI} is exactly the same as in \cite{MMT}, Section 4,
using Theorem~\ref{TH2}, the monotonicity arguments (Proposition \ref{MONOTONICITY1}) and
Theorem \ref{weaku}. It follows closely the proof of Theorem \ref{TH1} in the present paper.

\medskip

The proof of Corollary \ref{cor:1} is omitted since it is the same as the proof of Corollary
1 in \cite{MMT}.

\section{Weak convergence and well-posedness results}\label{sec:6}
 
 \subsection{Weak convergence}
\begin{theorem}[Weak  continuity of the BO flow map]\label{weaku}
Let $(u_n)$ be a sequence of global $H^{\frac 12}$ solutions of equation \eqref{BO}.
 Assume that $u_n(0) \rightharpoonup u_0$ in 
$H^{\frac 12}$ weak and let $u(t)$ be the
solution of \eqref{BO} corresponding to $u(0)=u_0$. Then, for all $t\in \mathbb{R}$,
$u_n(t)\rightharpoonup u(t)$ in $H^{\frac 12}$ weak.
\end{theorem}

\noindent\emph{Proof of Theorem \ref{weaku}.}
Let $u_{0,n}=u_n(0)$. It is sufficient to prove the result for $T\in [0,1]$.

\medskip

\noindent\emph{Step 1. $H^2 case$.} Here, we assume $u_{0,n}\rightharpoonup u_0$ in 
$H^2$. Let $w_n=u_n-u$. The equation for $w_n$ is 
\begin{equation}\left\{
\begin{aligned}
&w_{nt} + \mathcal{H}(w_n)_{xx}+ u_n w_{nx} + u_x w_n=0\\
&w_n(0)=\psi_n,\quad \psi_n=u_{0,n}-u_0.
\end{aligned}
\right. \end{equation}
Fix $t=T$, $g\in C^{\infty}_0(\mathbb{R})$. 
For a function $\widetilde u$ to be determined, we consider the solution 
$v(t)$ of
\begin{equation*}
\left\{
\begin{aligned}
&		v_t+ \mathcal{H} v_{xx} + (\widetilde u v)_x - u_x v = 0,
\\ &	 	v(T)=g.
\end{aligned}
\right.
\end{equation*}
Then
$$
\int w_n(T,x)g(x)dx - \int \psi_n(x) v(0,x)dx 
= \int_0^T \int w_{nt}(t)v(t) + \int_0^T \int w_n(t) v_t =\mathbf{I}+\mathbf{II}.
$$
$$
\mathbf{I}=
\int_0^T \int w_n (\mathcal{H} v_{xx} + (u_n v)_x - u_x v),
\qquad 
\mathbf{II}=
- \int_0^T \int w_n (\mathcal{H} v_{xx} + (\widetilde u v)_x - u_x v)
$$
so that
$$
\int w_n(T,x)g(x)dx - \int \psi_n(x) v(0,x)dx =
\int_0^T \int w_n ((u_n-\widetilde u)v)_x =
-\int_0^T \int w_{nx} (u_n-\widetilde u) v.
$$
We can assume, after passing to a subsequence, that 
$u_n-\widetilde u \to 0$ in $L^2_{loc}(\mathbb{R}\times [0,T])$.
Next, we will show that given $\varepsilon>0$, there exists $R>0$ such that
$$
\bigg|\int_0^T \int_{|x|>R} w_{nx} (u_n-\widetilde u) v \bigg|\leq \varepsilon,
\quad \text{uniformly in $n$.}
$$
In fact, since $\|w_{nx}\|_{L^\infty}\leq C$, $\sup_t \|v\|_{L^2}\leq C$ and
$\sup_t \|u_n-\widetilde u\|_{L^2}\leq C$, the claim is clear.

But then, $\mathbf{I}+\mathbf{II}\to 0$ as $n\to +\infty$.
We only needed $\sup_t \|v\|_{L^2}\leq C$,
which needs $\widetilde u_x \in L^\infty$, $u_x \in L^\infty$, which are both clear.
(We use the energy method to bound $v$.)

\medskip

\noindent\emph{Step 2. General case.} 
Fix $N$ large, define $u_{0,n}^N$ such that $\widehat {u_{0,n}^N}(\xi)=\mathbf{1}_{[-N,N]}(\xi) \widehat u_{0,n}(\xi)$, where $\mathbf{1}_I$ is the characteristic function of $I$.
Note that
$$
\|u_{0,n}^N - u_{0,n}\|_{L^2}^2 
= \int_{|\xi|\geq N} \big|\widehat {u_{0,n}^N}(\xi)\big|^2
\leq \frac 1N \|u_{0,n}\|_{H^{\frac 12}}^2 \leq \frac C N,
$$
so that $u_{0,n}^N\to u_{0,n}$ in $L^2$ as $N\to +\infty$, uniformly in $n$.

Fix $g \in C^\infty_0$, $T\in \mathbb{R}$, $\varepsilon>0$.
The proof of the $L^2$ continuity of the flow map (see \cite{IK}) shows that
$$
\sup_{t\in [0,1]} \|u^N(t) - u(t)\|_{L^2} \leq C \|u_0^N-u_0\|_{L^2},
\quad
\sup_{t\in [0,1]} \|u_n^N(t)-u_n(t)\|_{L^2} \leq C \|u_{0,n}^N-u_{0,n}\|_{L^2}
$$
for some universal constant $C>0$.
We fix $N$ such that
$$
\bigg|
\int (u_n(T)-u(T)) g -  \int (u_n^N(T)-u^N(T)) g 
\bigg|
\leq \frac \varepsilon 2, \quad \text{uniformly in $n$.}
$$
But, for fixed $N$, we let $n\to +\infty$, and use step 1 and the proof is
concluded.

\subsection{Well-posedness result for the nonlinear BO equation with potential}
In this subsection, for $0<b<b_0$, $b_0$ small, we consider the IVP
\begin{equation}
\left\{
\begin{array}{l}
v_t=(-\mathcal{H}v_x)_x - (Q(x{-}t) v)_x -\frac b 2 (v^2)_x=0  \quad \text{on $[-T,T]\times \mathbb{R}$},\\[5pt]
v(t=0,x)= v_0(x)  \quad \text{on $\mathbb{R}$}.
\end{array}
\right.\label{BOQ}
\end{equation}
The well-posedness of the Cauchy problem in $L^2$  for this equation is clear from \cite{IK}
since $u(t,x)=Q(x{-}t)+bv(t,x)$ satisfies the BO equation. 
Our main concern is a result of equicontinuity of the map $t\mapsto v(t)$ in $L^2$ 
with respect to $b$. 
To establish such a result we follow the strategy of \cite{IK} on equation
\eqref{BOQ}, using the special form of $Q$ and keeping track of the dependency in $b$.

\begin{theorem}\label{THA}
\begin{description}
\item{{\rm (a)}} Assume $v_0\in H^\infty$. Then, there exists $T=T(Q)>0$ and a unique solution $v=S_b^{\infty}(v_0)$ of \eqref{BOQ} in $[-T,T]$, 
$v\in C([-T,T],H^\infty)$.
\item{{\rm (b)}} There exists a constant $C$, independent of $b$ such that
\begin{equation}\label{eq:b}
	\sup_{t\in [-T,T]} \|v(t)\|_{H^2}\leq C \|v_0\|_{H^2}.
\end{equation}
\item{{\rm (c)}} The mapping $S_b^\infty$ extends uniquely to a continuous mapping
$S_b^0 : L^2 \to C([-T,T],L^2)$, and there exists $C$, independent of $b$ such that
\begin{equation}\label{eq:ci}
\sup_{t\in [-T,T]} \|v(t)\|_{L^2} \leq C \|v_0\|_{L^2}.
\end{equation}
Moreover, given $v_0\in L^2$, $\|v_0\|_{L^2}\leq 2$,
for any $\varepsilon>0$, there exits $\delta=\delta(v_0,\varepsilon)>0$
($\delta$ independent of $b$) such that for any $v_1\in L^2$, $\|v_1\|_{L^2}\leq 2,$
\begin{equation}\label{eq:cii}
\|v_0-v_1\|_{L^2}\leq \delta\quad \Rightarrow \quad
\sup_{t\in [-T,T]} \|S_b^0(v_0)(t)-S_b^0(v_1)(t)\|_{L^2}\leq \varepsilon.
\end{equation}
Finally, there exists $\widetilde \delta=\widetilde \delta(v_0,\varepsilon)>0$ (independent of $b$) such that for any $t,t'\in [-T,T]$,
\begin{equation}\label{eq:ciii}
|t-t'|\leq \widetilde \delta \quad \Rightarrow \quad 
\|S_b^0(v_0)(t)-S_b^0(v_0)(t')\|_{L^2} \leq \varepsilon.
\end{equation}
\end{description}
\end{theorem}

\noindent\emph{Reduction of the proof.}
For $0<\lambda\ll 1$, consider $v_\lambda(t,x)=\lambda v(\lambda^2 t, \lambda x)$.
Then $v_\lambda$ solves
\begin{equation}
\left\{
\begin{array}{l}
(v_\lambda)_t=(-\mathcal{H}(v_\lambda)_x)_x - (\lambda Q(\lambda x{-}\lambda^2 t) v_\lambda)_x -\frac b 2 (v_\lambda^2)_x=0  \quad \text{on $[-T,T]\times \mathbb{R}$},\\[5pt]
v_\lambda(t=0,x)= v_{0,\lambda}(x)  \quad \text{on $\mathbb{R}$},\quad
v_{0,\lambda}(x)=  \lambda v_0( \lambda x).
\end{array}
\right.\label{BOQl}
\end{equation}
Define
$$
Q_\lambda(t,x)=\lambda Q(\lambda x{-}\lambda^2 t).
$$
Then the proof of Theorem \ref{THA}  reduces to prove that for the following
(IVP)
\begin{equation}
\left\{
\begin{array}{l}
v_t=(-\mathcal{H}v_x)_x - (Q_\lambda(t,x)v)_x -\frac b 2 (v^2)_x=0  \quad \text{on $[-1,1]\times \mathbb{R}$},\\[5pt]
v(t=0,x)= v_{0}(x)  \quad \text{on $\mathbb{R}$},\quad
\|v_{0}\|_{L^2}\leq \lambda^{\frac 12},
\end{array}
\right.\label{BOQl2}
\end{equation}
 we have

\begin{theorem}\label{THB}
There exists $b_0$, $\lambda>0$ small enough such that 
if $0<b<b_0$, the following hold
\begin{description}
\item{{\rm (a)}} Assume $v_0\in H^\infty$. Then, there exists  a unique solution $v=S_b^{\infty}(v_0)$ of \eqref{BOQl2} in $[-1,1]$, 
$v\in C([-1,1],H^\infty)$.
\item{{\rm (b)}} There exists a constant $C$, independent of $b$ such that
\begin{equation}\label{eq:b2}
	\sup_{t\in [-1,1]} \|v(t)\|_{H^2}\leq C \|v_0\|_{H^2}.
\end{equation}
\item{{\rm (c)}} The mapping $S_b^\infty$ extends uniquely to a continuous mapping
$S_b^0 : L^2 \to C([-1,1],L^2)$, and there exists $C$, independent of $b$ such that
\begin{equation}\label{eq:ci2}
\sup_{t\in [-1,1]} \|v(t)\|_{L^2} \leq C \|v_0\|_{L^2}.
\end{equation}
Moreover, given $v_0\in L^2$, $\|v_0\|_{L^2}\leq \lambda^{\frac 12}$,
for any $\varepsilon>0$, there exits $\delta=\delta(v_0,\varepsilon)>0$
($\delta$ independent of $b$) such that for any $v_1\in L^2$, 
$\|v_1\|_{L^2}\leq \lambda^{\frac 12}$,
\begin{equation}\label{eq:cii2}
\|v_0-v_1\|_{L^2}\leq \delta\quad \Rightarrow \quad
\sup_{t\in [-1,1]} \|S_b^0(v_0)(t)-S_b^0(v_1)(t)\|_{L^2}\leq \varepsilon.
\end{equation}
Finally, there exists $\widetilde \delta=\widetilde \delta(v_0,\varepsilon)>0$ (independent of $b$) such that for any $t,t'\in [-1,1]$,
\begin{equation}\label{eq:ciii2}
|t-t'|\leq \widetilde \delta \quad \Rightarrow \quad 
\|S_b^0(v_0)(t)-S_b^0(v_0)(t')\|_{L^2} \leq \varepsilon.
\end{equation}
\end{description}
\end{theorem}

The proof of Theorem \ref{THB} is based on the following three propositions.

\begin{proposition}\label{p:1}
Assume $v_0\in H^\infty$, then there exists $T=T(\|v_0\|_{H^2})$
and a unique solution $v$ of \eqref{BOQl2} in $(-T,T)$. Also, for any $\sigma\geq 2$,
\begin{equation}
\sup_{t\in (-T,T)}\|u(t)\|_{H^\sigma} \leq 
C(\sigma, \|v_0\|_{\sigma}, \sup_{t\in (-T,T)} \|v(t)\|_{H^2}).
\end{equation}
In particular, the constant $C$ is independent of $b$, ($0<b<b_0$) and $\lambda<1$.
\end{proposition}
Proposition \ref{p:1} is a consequence of the energy method, taking into account that
$$\|\partial_xQ_\lambda\|_{L^1((-1,1),L^\infty_x)}\leq C.$$
\begin{proposition}\label{p:2}
For $\lambda$ small enough, we have that if $T\in (0,1]$, $\|v_0\|_{L^2}\leq \lambda^{\frac 12}$, $v=S^\infty(v_0)\in C((-T,T),H^\infty)$ is a solution, then
$$
\sup_{t\in [-T,T]} \|v(t)\|_{H^2} \leq C \|v_0\|_{H^2},
$$
where $C$ is independent of $b$ $(0<b<b_0)$.
\end{proposition}
\begin{proposition}\label{p:3}
For $v_0\in H^\infty$, $N\in [1,\infty)$, $\|v_0\|_{L^2}\leq \lambda^{\frac 12}$,
let $\widehat {v_0^N}(\xi)=\mathbf{1}_{[-N,N]}(\xi) \hat v_0(\xi)$, $v_0^N\in H^\infty$.
Then,
$$
\sup_{t\in (-1,1)} \|S_b^\infty(v_0)(t)-S_b^\infty(v_0^N)(t)\|_{L^2}
\leq C \|v_0-v_0^N\|_{L^2},
\quad
\sup_{t\in (-1,1)} \|S_b^\infty(v_0)(t)\|_{L^2}\leq C \|v_0\|_{L^2}.
$$
where $C$ is independent of $b$ $(0<b<b_0)$.
\end{proposition}

\noindent\emph{Proof of Theorem \ref{THB} from Propositions \ref{p:1}, \ref{p:2}
and \ref{p:3}}.
First, note that Propositions \ref{p:1} and \ref{p:2} clearly give (a) and (b) in 
Theorem \ref{THB}.
Let us turn to the proof of (c): it suffices to show first that if $v_{0,n}\in H^\infty$,
$\lim_{n\to +\infty} v_{0,n} = v_0$ in $L^2$, the sequence $S_b^\infty(v_{0,n})$ is 
Cauchy in $C([-1,1],L^2)$. Let $\varepsilon>0$ be given. We want to show that 
there exists $M_\varepsilon$ (independent of $b$) such that
$$
m,\ n \geq M_\varepsilon \quad \Rightarrow \quad
\sup_{t\in [-1,1]} \|S_b^\infty(v_{0,n})(t)-S_b^\infty(v_{0,m})(t)\|_{L^2}\leq \varepsilon.
$$
Observe that 
$$
\|v_{0,n}-v_{0,n}^N\|_{L^2}\leq \|v_{0}-v_{0}^N\|_{L^2}
+\|v_{0}-v_{0,n}\|_{L^2}.
$$
Hence, we can fix $N=N(\varepsilon,v_0)$ large and $M^1_\varepsilon$ large
such that $\|v_{0,n}-v_{0,n}^N\|_{L^2}\leq \frac {\varepsilon}{4C}$, for $n\geq M^1_\varepsilon$,
where $C$ is the constant in Proposition \ref{p:3} ($\|v_{0,n}\|_{L^2} \leq \lambda^{\frac 12}$).
Then, by Proposition \ref{p:3}, for $n\geq M^1_\varepsilon$,
$
\sup_{t\in [-1,1]} \|S_b^\infty(v_{0,n})(t)-S_b^\infty(v_{0,n}^N)(t)\|_{L^2}
\leq \frac \varepsilon 4.
$

It remains to estimate
$
\sup_{t\in [-1,1]} \|S_b^\infty(v_{0,n}^N)(t)-S_b^\infty(v_{0,m}^N)(t)\|_{L^2}.
$
But energy estimates for the difference equation give
\begin{align*}
& \quad \sup_{t\in [-1,1]} \|S_b^\infty(v_{0,n}^N)(t)-S_b^\infty(v_{0,m}^N)(t)\|_{L^2}\\
&\leq \|v_{0,n}^N - v_{0,m}^N\|_{L^2} 
\exp\left(C \int_{-1}^{1} \|\partial_x (S_b^\infty(v_{0,n}^N)(t)\|_{L^\infty_x}+
C\|\partial_x (S_b^\infty(v_{0,m}^N)(t)\|_{L^\infty_x}\right)\\
& \leq \|v_{0,n}  - v_{0,m}\|_{L^2}
\exp\left(C \sup_{t\in (-1,1)} \| S_b^\infty(v_{0,n}^N(t)\|_{H^2}+
+C\sup_{t\in (-1,1)} \| S_b^\infty(v_{0,m}^N(t)\|_{H^2}\right)\\
&\leq \|v_{0,n}  - v_{0,m}\|_{L^2} \exp\left(CN^2 \|v_{0,n}\|_{L^2}+
CN^2 \|v_{0,n}\|_{L^2}\right)
\leq \|v_{0,n}  - v_{0,m}\|_{L^2} \exp(CN^2) \leq \frac \varepsilon 2,
\end{align*}
for $n$, $m$ large
(we have used the estimate of Proposition \ref{p:2}).
Also, by Proposition \ref{p:3}, we have
$\sup_{t\in (-1,1)} \|S_b^\infty(v_{0,n})(t)\|_{L^2}\leq C$.
Thus, we obtain the unique extension $S_b^0$ and \eqref{eq:ci2} holds.

To check \eqref{eq:cii2}, fix $v_0$, $\|v_0\|_{L^2}\leq \lambda^{\frac 12}$,
let $\varepsilon>0$ be given. With $C$ as in Proposition \ref{p:3}, find $N$
($N=N(\varepsilon,v_0)$)
so large that  $\|v_0-v_0^N\|_{L^2} \leq \frac \varepsilon {8C}.$
Now find $\delta_1=\delta_1(\varepsilon,v_0)$ so small that if 
$\|v_0-v_1\|_{L^2}\leq \delta_1$, then
$ \|v_1-v_1^N\|_{L^2} \leq \frac \varepsilon{4C}$. We have
\begin{align*}
   &\sup_{t\in [-1,1]} \|S_b^0(v_0)(t)-S_b^0(v_1)(t)\|_{L^2} \leq 
\sup_{t\in [-1,1]} \|S_b^0(v_0)(t)-S_b^0(v_0^N)(t)\|_{L^2}
\\& \quad +\sup_{t\in [-1,1]} \|S_b^0(v_1)(t)-S_b^0(v_1^N)(t)\|_{L^2}
+\sup_{t\in [-1,1]} \|S_b^0(v_1^N)(t)-S_b^0(v_0^N)(t)\|_{L^2}.
\end{align*}
By Proposition \ref{p:3}, the first two terms are smaller than $\frac \varepsilon 2$.
For the last one, we again use the energy estimate and get, as before
$$
\sup_{t\in [-1,1]} \|S_b^0(v_1^N)(t)-S_b^0(v_0^N)(t)\|_{L^2}\leq
C \|v_1-v_0\|_{L^2} \exp(CN^2),
$$
using Propositions \ref{p:2} and \ref{p:3} and \eqref{eq:cii2} follows.

For \eqref{eq:ciii2}, first find $N=N(\varepsilon,v_0)$ so large that
$\|v_0-v_0^N\|_{L^2} \leq \frac {\varepsilon}{4C}$, where $C$ is as in Proposition \ref{p:3}.
Then,
$\sup_{t\in [-1,1]} \|S_b^0(v_0)(t)-S_b^0(v_0^N)(t)\|_{L^2}\leq \frac \varepsilon 4$
and we are reduced to showing, for $N$ fixed that if $|t-t'|\leq \widetilde \delta$,
then
$
\|S_b^0(v_0^N)(t)-S_b^0(v_0^N)(t')\|_{L^2} \leq \frac \varepsilon 2.
$

Let $f(t)=\|S_b^0(v_0^N)(t)\|_{L^2}^2$. The energy method, combined with Proposition \ref{p:2}
shows that
$|f'(t)|\leq f(0) \exp(CN^2)$. But then, for $|t-t'|\leq \widetilde \delta_1$,
$|f(t)-f(t')|\leq \frac \varepsilon 4$. But
\begin{align*}
& \|S_b^0(v_0^N)(t)-S_b^0(v_0^N)(t')\|_{L^2}^2
=f(t)+f(t')-2 \int S_b^0(v_0^N)(t).S_b^0(v_0^N)(t')dx\\
&= f(t')-f(t)
+ 2 \int S_b^0(v_0^N)(t)[S_b^0(v_0^N)(t)-S_b^0(v_0^N)(t')] dx.
\end{align*}
Let $v^N(t)=S_b^0(v_0^N)(t)$. The second term equals
$$
2 \int v^N(t)\int_{t'}^t \partial_s v^N(s) ds dx=
2 \int^{t}_{t'} \int v_{N}(t)
[-\mathcal{H} \partial_x^2 v^N(s) - (Q_\lambda v^N)_x - \frac b2
((v^{N})^2)_x(s)] dx ds.
$$
But by Proposition \ref{p:2},
$
\sup_{t\in [-1,1]} \|v^N(t)\|_{H^2}\leq
C\|v_0^N\|_{H^2} \leq CN^2.
$
Thus, the second term is controlled by $C|t-t'| N$, and the proof is complete,
provided we prove Propositions \ref{p:2} and \ref{p:3}.

\medskip

\noindent\emph{Proof of Propositions \ref{p:2} and \ref{p:3}.}
\emph{Step 1.} Assume $v_0\in H^\infty$, $\|v_0\|_{H^2}\leq M$ and $0<T\leq 2$,
$v=S_b^\infty(t)$ exists in $[-T,T]$.
Then, there exist $\lambda_0=\lambda_0(M)$, $b_0=b_0(M)$
such that for $0<\lambda<\lambda_0$, $0\leq b<b_0$, we have
\begin{equation}\label{page10}
\sup_{t\in [-T,T]} \|v(t)\|_{H^2} \leq 2 \|v_0\|_{H^2}.
\end{equation}
Proof of \eqref{page10}. Note that $\|\partial_x^k Q\lambda\|_{L^\infty} \leq C_k \lambda^{k+1}$.
Let $f(t)=\|v(t)\|_{H^2}^2$. The standard energy method shows that
$$
|f'(t)|\leq C (\lambda_0^2  + b_0 \|\partial_x v(t)\|_{L^\infty_x}) f(t)
\leq (\lambda_0^2 +b_0 (f(t))^{\frac 12}) f(t).
$$
Integrating the ODE gives the result.

As a corollary, we obtain under the circumstances of Step 1 that $v$ exists in $(-1,1)$ and
$$
\sup_{t\in [-2,2]} \|v(t)\|_{H^2} \leq 2 \|v_0\|_{H^2}.
$$

\emph{Step 2.} From now on, we will follow closely \cite{IK}. Some of the ideas used before
were developed in a forthcoming paper \cite{HIKK}. We have now reduced everything to
\emph{a priori} estimates. We will change notation slightly to match \cite{IK}.
We then study the problem
\begin{equation}\left\{
\begin{aligned}
& u_t + \mathcal{H} u_{xx} + (Q_\lambda u)_x + b (\tfrac 12 u^2)_x=0\quad
(t,x)\in (-1,1)\times \mathbb{R}, \\
& u_{|t=0}=\phi,\quad \|\phi\|_{L^2} \leq \lambda^{\frac 12}, 
\end{aligned}\right.
\end{equation}
We use the notation $P_{\rm low}$, $P_{\rm \pm high}$ as in \cite{IK}:
$$
P_{\rm low} \text{ defined by the Fourier multiplier } 
\xi \to \mathbf{1}_{[-2^{10},2^{10}]}(\xi);
$$
$$
P_{\rm \pm high} \text{ defined by the Fourier multiplier } 
\xi \to \mathbf{1}_{[2^{10},\infty)}(\pm \xi);
$$
$$
P_{\pm} \text{ defined by the Fourier multiplier } 
\xi \to \mathbf{1}_{[0,\infty)}(\pm \xi).
$$
Let $\phi_0=P_{\rm low} \phi\in H^\infty$, real-valued, $\|\phi_0\|_{H^2}\leq 2^{20}=M$.
We choose $\lambda_0$, $b_0$ as in Step 1 and its corollary, so that Proposition \ref{p:1}
and these results gives, with $u_0^{(1)}=S^\infty_b(\phi_0)(t)$ that
$$
\sup_{t\in [-2,2]} \|\partial_t^{\sigma_1} \partial_x^{\sigma_2} u_0^{(1)}\|_{L^2_x}
\leq C_{\sigma_1,\sigma_2} \|\phi\|_{L^2}, \quad \sigma_i \geq 0.
$$
Let $\widetilde u=u-u_0^{(1)}$. The equation for $\widetilde u$ is
\begin{equation}\left\{
\begin{aligned}
&\widetilde u_t + \mathcal{H} \widetilde u_{xx} + (Q_\lambda \widetilde u)_x + b (u_0^{(1)} \widetilde u)_x
+ b (\tfrac 12 \widetilde u^2)_x=0,\\
& \widetilde u_{|t=0} = P_{\rm + high} \phi + P_{\rm - high} \phi.
\end{aligned}\right.
\end{equation}
Let now $u_0(t,x)=Q_\lambda(t,x) + b u_0^{(1)} (t,x)$. Then
$
\sup_{t\in [-2,2]} \|\partial_t^{\sigma_1} \partial_x^{\sigma_2} u_0\|_{L^2_x}
\leq C_{\sigma_1,\sigma_2} (\lambda_0^{\frac 12} + b_0).
$
We now want to construct $U_0$ similarly to \cite{IK}, with the following properties
$\partial_x U_0(t,x)= \frac 12 u_0(t,x)$, $U_0(0,0)=0$
and $\sup_{t\in [-2,2]} \|\partial_t^{\sigma_1} \partial_x^{\sigma_2} U_0(t,.)\|_{L^2_x}
\leq C_{\sigma_1,\sigma_2} (\lambda_0^{\frac 12} + b_0)$
where $\sigma_1,\sigma_2\geq 0$, $(\sigma_1,\sigma_2)\neq (0,0)$.

Since $Q_\lambda(t,x)= \frac {4\lambda}{1+(\lambda x -\lambda^2 t)^2}$, we set
$U_0^{(2)} (t,x)=2 \arctan(\lambda x -\lambda^2 t)$.
We next recall the equation $u_0^{(1)}(t,x)$ verifies:
$$
\partial_t (\tfrac 12 {u_0^{(1)}} ) + \mathcal{H} \partial_x^2(\tfrac 12 {u_0^{(1)}} )
+ \partial_x (Q_\lambda \tfrac 12 {u_0^{(1)}} ) + b \partial_x( (\tfrac 12 {u_0^{(1)}})^2)=0.
$$
We then define first $U_0^{(1)}(t,0)$ by the formula
$$
\partial U_0^{(1)}(t,0) + \mathcal{H} \partial_x(\tfrac 12 {u_0^{(1)}(t,0)} ) 
+ Q_{\lambda}(t,0) \tfrac 12 {u_0^{(1)}(t,0)} + b  (\tfrac 1  2 {u_0^{(1)}(t,0)} )^2=0,
\quad U_0^{(1)} (0,0) = 0.
$$
We then construct $U_0^{(1)}(t,x)$ by $\partial_x U_0^{(1)} (t,x)=\frac 12 u_0^{(1)}(t,x)$.
Notice that $U_0^{(1)}$ is real-valued. Using the equation for $u_0^{(1)}$, we have
$$
\partial_x\left(\partial_t U_0^{(1)} + \mathcal{H} \partial_x^2 U_0^{(1)}
+ Q_\lambda \partial_x U_0^{(1)}+ b (\partial_x U_0^{(1)})^2 \right)=0
\quad \text{on $\mathbb{R}\times [-2,2]$}. 
$$
But then, on $\mathbb{R}\times [-2,2]$, we have
$$
\partial_t U_0^{(1)}(t,x) + \frac 12 \mathcal{H} \partial_x u_0^{(1)}(t,x)
+ Q_\lambda(t,x) \frac 12 u_0^{(1)}(t,x) + 
\frac b4 (u_0^{(1)}(t,x))^2.
$$
We then define $U_0(t,x)=b U_0^{(1)}(t,x) + U_0^{(2)}(t,x)$, and
all our properties hold. We recall that
\begin{equation}\left\{
\begin{aligned}
&\widetilde u_t + \mathcal{H} \widetilde u_{xx} + (u_0 \widetilde u)_x 
+ b (\tfrac 12 \widetilde u^2)_x=0,\\
& \widetilde u_{|t=0} = P_{\rm + high} \phi + P_{\rm - high} \phi.
\end{aligned}\right.
\end{equation}
We now proceed as in Section 2 of \cite{IK}. We define
$P_{+\rm high} \widetilde u= e^{-iU_0} w_+$,
$P_{-\rm high} \widetilde u= e^{iU_0} w_-$ and
$P_{\rm low} \widetilde u = w_0$.
Applying $P_{+\rm high}$, $P_{-\rm high}$, $P_{\rm low}$ to the above equation
and using the definitions above, we have
(we write the equation for $w_+$, the one for $w_-$ is analoguous, the one for $w_0$ will be written later).
Following the argument in \cite{IK}, one gets:
\begin{align*}
& (w_+)_t + \mathcal{H} \partial_x^2 w_+ =
- \frac b2 e^{iU_0} P_{+\rm high} \partial_x((e^{-iU_0} w_+ + e^{iU_0} w_-
+w_0)^2)
\\ &- e^{-iU_0} P_{+\rm high} \partial_x ( u_0(e^{iU_0} w_- + w_0)) + 
e^{iU_0} (P_{-\rm high} +P_{\rm low}) (e^{iU_0} u_0 \partial_x w_+)
+ 2i P_- \partial_x^2 w_+
\\& -e^{iU_0} P_{+\rm high} (\partial_x (u_0e^{-iU_0} w_+)) 
+i w_+ \left[(U_0)_t - i (U_0)_{xx} - ((U_0)_x)^2\right],
\end{align*}
and so after more calculations, we get
\begin{align*}
& (w_+)_t + \mathcal{H} \partial_x^2 w_+ =
- \frac b2 e^{iU_0} P_{+\rm high} \partial_x((e^{-iU_0} w_+ + e^{iU_0} w_-
+w_0)^2)\\&
- e^{-iU_0} P_{+\rm high} \left[\partial_x ( u_0 P_{-\rm high}(e^{iU_0} w_-)+ u_0 P_{\rm low}(w_0)) \right]\\ &
+
e^{iU_0} (P_{-\rm high} 
+P_{\rm low}) \left[\partial_x( u_0 P_{+\rm high}(e^{-iU_0} w_+))\right]\\&
+ 2i P_- \left[\partial_x^2 (e^{iu_0}P_{+\rm high} (e^{-iU_0} w_+))\right]
+i w_+ \left[(U_0)_t + \mathcal{H} \partial_x^2 U_0 + (\partial_x U_0)^2 + iP_+ \partial_x U_0
\right],
\end{align*}
We recall $\partial_x U_0^{(2)} = \frac 12 Q_\lambda$ and that $Q_\lambda$ solves
$\partial_t Q_\lambda + \mathcal{H}\partial_x^2 Q_\lambda + \partial_x(\tfrac 12
Q_\lambda^2)=0$ or
$\partial_t U_0^{(2)} + \mathcal{H} \partial_x^2 U_0^{(2)} =- \tfrac 14 Q_\lambda^2$
and
$\partial_t U_0^{(1)} + \mathcal{H}\partial_x^2 U_0^{(1)} = -Q_\lambda \partial U_0^{(1)} - b (\partial_x U_0^{(1)})^2.$ Hence, $
\partial_t U_0 + \mathcal{H} \partial_x^2 U_0 + (\partial_x U_0)^2=0$
and we get $\partial w_+ + \mathcal{H} w_+ = E_+(w_+,w_-,w_0)$, where $E_+$
is defined as in \cite{IK}, p. 756, except that the first term is multiplied now
by $b$.  The equation for $w_-$ and $E_-$ is similar. The equation for
$w_0$ writes
$$
\partial_t (P_{\rm low} \widetilde u) 
+\mathcal{H} \partial_x^2 P_{\rm low} \widetilde u + P_{\rm low} \partial_x(u_0 \widetilde u)
+ \tfrac b2 P_{\rm low} \partial_x((\widetilde u)^2)=0,
$$
where $\widetilde u= e^{-iU_0}w_+ + e^{iU_0} w_-+ w_0$. Next, we note that, with
$\delta=(\lambda_0^{\frac 12}+b_0)$, the estimates (10.19) in \cite{IK} hold.
Because of this and the form of $E_+$, $E_-$, $E_0$, just
as in Proposition 10.5 in \cite{IK}, we have
\begin{align*}
\|\psi(t) (\mathbf{E}(\mathbf{w}) - \mathbf{E}(\mathbf{w}'))\|_{N^\sigma}
&\leq C b_0 \|\mathbf{w}-\mathbf{w'}\|_{F^\sigma}  
(\| \mathbf {w}\|_{F^0} + \| \mathbf{w}'\|_{F^0})\\
&\quad +C b_0 \|\mathbf{w}-\mathbf{w'}\|_{F^0}  
(\|\mathbf{w}\|_{F^\sigma} + \|\mathbf{w}'\|_{F^\sigma})+
C \delta \|\mathbf{w}-\mathbf{w'}\|_{F^\sigma} .
\end{align*}
Note that $\mathbf{w}=(w_+,w_-,w_0)$ and
$\mathbf{E}(\mathbf{w})=(E_+(w_+,w_-,w_0),E_-(w_+,w_-,w_0),E_0(w_+,w_-,w_0))$
as in \cite{IK}.
The rest of the notation (the norm $\|.\|_{N^\sigma}$ and the function $\psi$)
is also taken from \cite{IK}.
We have a slightly different formula for $E_0$, but (10.27) in \cite{IK}
gives the estimate in our case also.

We then construct a solution to 
\begin{equation*}\left\{\begin{aligned}
& \mathbf{v}_t + \mathcal{H} \mathbf{v}_{xx}=\mathbf{E}(\mathbf{v}) \quad
\text{on } \mathbb{R}\times [-\tfrac 54,\tfrac 54],\\
& \mathbf{v}(0)=\Phi,
\end{aligned}\right.\end{equation*}
as in (10.32)-(10.37) in \cite{IK}. Note that (10.35) 
and $\|v(\Phi)-v(\Phi')\|_{F^0([-\frac 54,\frac 54])} \leq
C \|\Phi-\Phi'\|_{\widetilde H^0}$ hold here too.
Next, with $\Phi=(\phi_+,\phi_-,\phi_0)=
(e^{iU_0(0,.)} P_{+\rm high} \phi,e^{-i U_0(0,.)}P_{-\rm high} \phi,0)$,
$\Phi\in \widetilde H^{20}$, by Lemma 10.1 in \cite{IK}.

We next show $(w_+,w_-,w_0)=\mathbf{v}(\Phi)$ in $\mathbb{R}\times [-1,1]$.
This is as in \cite{IK}. Proposition \ref{p:2}, and the second estimate in Proposition \ref{p:3} now follow from the bounds on $\mathbf{v}(\Phi)$ i.e. (10.35).
For Proposition \ref{p:3}, note that for $N$ large, $U_0$ corresponding to 
$\phi$ and to $\phi_N$ defined by $\hat \phi_N=\mathbf{1}_{[-N,N]}(\xi) \hat \phi(\xi)$
are the same. We then have $u(t,x)=u_0^{(1)} + u - u_0^{(1)}=
u_0^{(1)} + \widetilde u = u_0^{(1)} + e^{-i U_0} w_+ + e^{iU_0} w_- + w_0$
and similarly, $u^N(t,x)=u_0^{(1)}+ u^N - u_0^{(1)}=
u_0^{(1)} + e^{-iU_0} w_+^N + e^{iU_0}Êw_-^{N} + w_0^N$.
Hence, 
\begin{align*}
\sup_{t\in [-1,1]} \|u(t,.)-u^N(t,.)\|_{L^2} &\leq
\sup_{t\in [-1,1]} \|w(t)-w^N(t)\|_{L^2} 
\\ & \leq C \|\psi(t) [w-w^N]\|_{F^0} \leq C \|\phi-\phi^N\|_{L^2}
\end{align*}
as desired, giving Proposition \ref{p:3}.

\appendix

\section{Appendix}\label{sec:A}

First, we recall the following inequalities:
\begin{lemma}\label{COMMUTATOR}
\begin{equation}\label{gn}
\forall 2\leq p< +\infty,\quad \|f\|_{L^p}
 \leq C_p \|f\|_{L^2}^{\frac 2p} \|D^{\frac 12} f\|_{L^2}^{\frac {p-2}p},
\end{equation}
\begin{equation}\label{COMM2}
\|D(fg)- g Df \|_{L^2}
\leq C \|f\|_{L^4} \|D g\|_{L^4},
\end{equation}
\begin{equation}\label{COMM1}
\|D^{\frac 12}(fg)- g D^{\frac 12} f \|_{L^2}
\leq C \|f\|_{L^4} \|D^{\frac 12} g\|_{L^4}.
\end{equation}
\end{lemma} 
Recall that \eqref{gn} is the Gagliardo-Nirenberg inequality, which follows from complex interpolation and Sobolev embedding.

Estimate \eqref{COMM2} is due to Calder\'on \cite{Ca}, see also Coifman and Meyer \cite{CM}, 
formula (1.1).

Estimate \eqref{COMM1} is a consequence of Theorem A.8 in \cite{KPV} for functions depending only on $x$, with the following choice of parameters: $\alpha=\frac 12$, $\alpha_1=0$, $\alpha_2=\frac 12$,
$p=2$, $p_1=p_2=4$.

\subsection{Proof of \eqref{toprove}}\label{secTOPROVE}
We claim that for a function $u(x)$ fixed in $H^2(\mathbb{R})$ 
\begin{equation}
 \int_{y=0} \partial_y (U^2) \Phi 
 = - 2 \iint_{\mathbb{R}^2_+} |\nabla U|^2 \Phi
   +\int_{y=0} U^2 \partial_y \Phi 
\end{equation}
where $U(x,y)$ is the harmonic extension of $u(x)$ in $\mathbb{R}^2_+$ and 
$\Phi(x,y)$ is defined in \eqref{defphicap}.

First, we observe that 
\begin{equation}\label{double}
U, \nabla U \in L^\infty(\mathbb{R}^2_+) \quad  \text{and}\quad 
\sup_{y>0}  |U(x,y)| \to 0 \text{ as $|x|\to +\infty$}.
\end{equation}
Indeed, from \cite{ST}, Theorem 1, p. 62, we have
$\sup_{y>0} |U(x,y)|\leq Mu(x)$, where $Mu(x)$ is the maximal function of $u$
(see \cite{ST} Chapter 1), and similarly,
$\sup_{y>0} |\partial_x U(x,y)|\leq Mu_x(x)$,
$\sup_{y>0} |\partial_y U(x,y)|\leq M(\mathcal{H}u_x)(x)$.
Moreover, from \cite{ST} Theorem 1, p. 5, since
$u,u_{x},\mathcal{H}u_{x}\in H^1\subset L^\infty$, we obtain $Mu, Mu_x, M(Hu_x)\in L^\infty$. Finally, since
$u \in H^1$, we have 
$|u(x)| \to 0$ as $|x|\to +\infty$, which implies by the definition of the maximal function (see \cite{ST}, page 4) that 
$Mu(x) \to 0$ as $|x|\to +\infty$.
Thus \eqref{double} is proved.

\medskip

Let $R>0$. We use the Green formula on $D_R^+=\{(x,y)\in \mathbb{R}_+^2 ~|~ x^2+y^2<R^2 \}$.
Let $\Gamma_R^+=\{(x,y)\in \mathbb{R}_+^2 ~|~ x^2+y^2=R^2 \}$
and $I_R={(x,0)~|~x\in [-R,R]}$. Then:
\begin{equation}\begin{split}
 \int_{\Gamma_R^+\cup I_R} \partial_n (U^2) \Phi   
& = - \iint_{D_R^+} (\Delta U^2) \Phi
+  \iint_{D_R^+}  U^2 \Delta \Phi +   \int_{\Gamma_R^+\cup I_R} U^2 \partial_n \Phi\\ 
& = - 2 \iint_{D_R^+} |\nabla U|^2 \Phi + \int_{\Gamma_R^+\cup I_R} U^2 \partial_n \Phi,
\end{split}
\end{equation}
where $\partial_n$ denotes the inward normal derivative since
$\Delta \Phi=0$ and $\Delta U^2= 2 |\nabla U|^2$.
Therefore, we only have to prove the following convergence results:
\begin{align}
& \lim_{R\to +\infty} \int_{\Gamma_R^+} \partial_n (U^2) \Phi  =
0,\quad 
\lim_{R\to +\infty} \int_{I_R} \partial_n (U^2) \Phi  =
\int_{y=0} \partial_y (U^2) \Phi =2 \int(\mathcal{H} u_x) u \varphi' \label{only1}\\
& \lim_{R\to +\infty} \int_{\Gamma_R^+} U^2 \partial_n \Phi  =
0,\quad
\lim_{R\to +\infty} \int_{I_R} U^2 \partial_n \Phi  =
\int_{y=0}U^2  \partial_y \Phi  =\int u^2 (\mathcal{H} \varphi'') \label{only2}.
\end{align}
The limits $\lim_{R\to +\infty} \int_{-R}^{R} (\mathcal{H} u_x) u \varphi' 
=\int (\mathcal{H} u_x) u \varphi' $ and
$\lim_{R\to +\infty} \int_{-R}^{R}  u^2 (\mathcal{H} \varphi'')
=\int  u^2 (\mathcal{H} \varphi'')$ are clear since $u\in H^1$.
Next, from the expression of $\Phi(x,y)$ in \eqref{defphicap}, we have
$\Phi(x,y)\leq C (1+y) R^{-2}$ on $\Gamma_R^+$. Therefore, from \eqref{double},
($d \sigma$ denotes the unit lenght element on $\Gamma_R^+$)
\begin{equation*}\begin{split}
\int_{\Gamma_R^+} | \partial_n (U^2) \Phi | & \leq 
\frac 1{R^2} \int_{\Gamma_R^+} |\nabla U| |U| (1+y) d\sigma
\\& \leq \frac C{R^2} \int_{\Gamma_R^+\cap \{|x|\leq \sqrt{R}\}} (1+y) d\sigma
+ C \sup_{|x|>\sqrt{R}, y>0} |U(x,y)|
\\ &\leq \frac C {\sqrt{R}}+ C \sup_{|x|>\sqrt{R}, y>0} |U(x,y)|
\end{split}
\end{equation*}
and so \eqref{only1} is proved. Estimate \eqref{only2} is proved similarly and
 is easier since $\partial_y \Phi$ has more decay than $\Phi$.

\subsection{Proof of \eqref{appendix1}}\label{secAPPENDIX1}
In the proof of \eqref{appendix1}, the time $t$ is fixed, so we set
$y_0=x_0+\lambda(t_0-t)$.

Let $\chi:\mathbb{R}\to \mathbb{R}$ be a $C^{\infty}$ function such that
$\chi=1$ on $[0,1]$, $\chi=0$ on $(-\infty,-1]\cap [2,+\infty)$ and
$\chi\leq 1$ on $\mathbb{R}$. Let $\chi_n(x)=\chi(x-n)$.
Then, by the Gagliardo Nirenberg inequality \eqref{gn},
we obtain
\begin{equation*}
\begin{split}
\int |\eta|^3 \varphi'(x-y_0)& \leq
\sum_{n\in \mathbb{Z}} \int_{n}^{n+1} |\eta|^3 \varphi'(x-y_0)
\leq \sum_{n\in \mathbb{Z}} \bigg(\int |\eta|^3 \chi_n^3\bigg) \sup_{[n-y_0,n+1-y_0]} \varphi'  \\
&\leq \sum_{n\in \mathbb{Z}} \bigg(\int |D^{\frac 12} (\eta \chi_n)|^2\bigg)^{\frac 12}\bigg(\int (\eta \chi_n)^2\bigg) \sup_{[n-y_0,n+1-y_0]} \varphi'.
\end{split}
\end{equation*}

By Lemma \ref{COMMUTATOR} and
\eqref{ortho}, we get
\begin{equation*}
\|D^{\frac 12} (\eta \chi_n)\|_{L^2}
\leq C \|(D^{\frac 12}  \eta) \chi_n\|_{L^2}+ C \|\eta\|_{L^4} 
\|D^{\frac 12} \chi_n\|_{L^4}
\leq C \|\eta\|_{H^{\frac 12}}\leq C \alpha_0.
\end{equation*}
Thus,
\begin{equation*}
\begin{split}
\int |\eta|^3 \varphi'(x-y_0)& \leq C\alpha_0
 \sum_{n\in \mathbb{Z}}  \bigg(\int (\eta \chi_n)^2\bigg) \sup_{[n-y_0,n+1-y_0]} \varphi'\leq C \alpha_0 \int \eta^2 \varphi'(x-y_0)
\end{split}
\end{equation*}
by the properties of $\chi$ and the following elementary remark:
\begin{equation}\label{rkonphi}
\forall y\in \mathbb{R},\quad
\sup_{[y,y+4]} \varphi'\leq C \inf_{[y,y+4]} \varphi'.
\end{equation}
Note that the constant $C$ is independent of $A$, for $A>1$.

\subsection{Properties of the operator $\mathcal{L}$}\label{secKERNEL}
We recall from \cite{We85}--\cite{W} and \cite{BB} the following properties of $\mathcal{L}$ (recall $\mathcal{L}\eta=-\mathcal{H}\eta_x + \eta - Q \eta$).

\begin{lemma}\label{76} The operator $\mathcal{L}$ is self-adjoint on $L^2$ and satisfies the following properties.
\begin{description}
\item{\rm (i)} The operator $\mathcal{L}$ has exactly one negative eigenvalue $\lambda_0$ of multiplicity $1$
with corresponding eigenfunction $f_0$, which can be chosen so that $f_0>0$.
\item{\rm (ii)} ${\rm Ker}\, \mathcal{L}= {\rm span}\{Q'\}$.
\item{\rm (iii)} There exists $\lambda>0$ such that, for all $z\in H^{\frac 12}$,
\begin{equation}\label{weinstein}
(z,Q)=(z,Q')=0 \quad \Rightarrow \quad
(\mathcal{L} z,z)\geq \lambda (z,z).
\end{equation}
\end{description}
\end{lemma}

\begin{remark}
Recall from Bennett et al. (\cite{BB}, Appendix B) that the spectrum of $\mathcal{L}$ is completely understood.
Indeed, the operator $\mathcal{L}$ has exactly four eigenvalues,
$\lambda_0=-\frac 12 (1+\sqrt{5})$, $0$, $\frac 12 (-1+\sqrt{5})$, $1$ and a continuous spectrum $[1,+\infty)$. \end{remark}

Now, we sketch a proof of Lemma \ref{76} using general arguments from \cite{We85}--\cite{W}.

\medskip

\noindent\emph{Sketch of proof.}
One easily checks that
$\mathcal{L}Q'=0$ (differentiate the equation of $Q(x+x_0)$ with respect to $x_0$ and
take $x_0=0$), and that 
$\mathcal{L}f_0=-\lambda_0 f_0$, where $f_0=Q+\frac 14 (1+\sqrt{5}) Q^2$ (by \eqref{surT}).
Moreover, the proof of (i) follows from the variational characterization of $Q$, see Proposition~4.2 of \cite{W}.
Recall that $\frac d{dc} \int Q_c^2 = \int Q^2>0$ (subcriticality) implies
that $\inf\{(\mathcal{L}f,f); \ (f,Q)=0,\ \|f\|_{L^2}=1\}=0$ (see proof of Proposition 5.1 in \cite{W}
and Proposition 3.1 in \cite{We86}).

Now, we give a new proof for (ii).
Let  $f\in L^2$ be such that $\mathcal{L} f=0$.
First, we remark that $f\in H^s$, for all $s\geq 0$. Moreover, by similar  estimates as in \cite{AT},
we have $|f(x)|\leq \frac C{1+x^2}.$ 
Integrating  $\mathcal{L} f=0$ on $\mathbb{R}$, we obtain $\int (f - fQ)=0$. But, we also have
$(f,Q)=-(f,\mathcal{L} S)=-(\mathcal{L} f,S)=0$ (see \eqref{surT}). Thus, $\int f=0$ and
we can define $g(x)=\int_{-\infty}^x f(s) ds\in L^2$, which satisfies $\mathcal{L}(g')=0$.
Let now $\widetilde g= g -a Q$ be such that $(\widetilde g,Q)=0$ and $\mathcal{L}(\widetilde g')=0$.
From \eqref{defLtilde} and \eqref{trucdec}, we obtain $\int |D^{\frac 12}\widetilde g|^2 + (\mathcal{L}\widetilde  g,\widetilde g)
\leq 0$. But, since $(\widetilde g,Q)=0$, we have $(\mathcal{L}\widetilde  g,\widetilde g)\geq 0$.
Thus, $\int |D^{\frac 12}\widetilde g|^2=0$ and $\widetilde g\equiv 0$, so that $g=aQ$
and $f=aQ'$.

Finally, we sketch the proof of (iii), which follows from the arguments of  the proof of Proposition 2.9 in \cite{We85}
(see also Section 6, example 4 in \cite{W}). By contradiction, assuming that
$$
\inf\{ (\mathcal{L}f,f); \ (f,Q)=(f,Q')=0,\ \|f\|_{L^2}=1\}=0,
$$
and using compactness arguments as in Proposition 2.9 in \cite{We85}, we obtain the existence of $f\in H^{\frac 12}$,
$\lambda,\beta,\gamma \in \mathbb{R}$ (Lagrange multipliers) such that
$$
(\mathcal{L}f,f)=0,\quad 
(\mathcal{L}-\lambda) f =\beta Q +\gamma Q',
\quad (f,Q)=(f,Q')=0,\quad \|f\|_{L^2}=1.
$$
But, taking the scalar product by $f$, we find $\lambda=0$. Then, taking the scalar product by $Q'$, we find $\gamma=0$.
Taking the scalar product with $S$ (see \eqref{surT}), using
$(S,Q)=\frac 12 (Q,Q)$ and $\mathcal{L}(S)=-Q$, we find $\beta=0$, so that $\mathcal{L} f=0$ and $(f,Q')=0$. This implies $f=0$ by (ii), a contradiction.

\end{document}